\documentclass{amsart}

\usepackage{amsmath, verbatim, amssymb, amsthm}
\usepackage[all]{xypic}

\numberwithin{equation}{section}
\theoremstyle{plain}
\newtheorem{theorem}[equation]{Theorem}
\newtheorem{proposition}[equation]{Proposition}
\newtheorem{corollary}[equation]{Corollary}
\newtheorem{lemma}[equation]{Lemma}

\newtheorem{conjecture}[equation]{Conjecture}

\newtheorem{sublemma}[equation]{Sublemma}
\newtheorem{assnot}[equation]{Assumption-Notation}

\theoremstyle{definition}
\newtheorem{defn}[equation]{Definition}
\newtheorem{notation}[equation]{Notation}

\newtheorem{example}[equation]{Example}

\theoremstyle{remark}

\newcommand{\beq}{\begin{equation}}
\newcommand{\eeq}{\end{equation}}

\newcommand{\ang}[1]{\langle #1 \rangle}

\newcommand{\st}{\left\vert\right.}

\newcommand{\bbar}[1]{\overline{#1}}

\providecommand{\abs}[1]{\lvert#1\rvert}

\DeclareMathOperator{\id}{{Id}}

\DeclareMathOperator{\Ext}{{Ext}}

\DeclareMathOperator{\Aut}{{Aut}}
\DeclareMathOperator{\Proj}{Proj}
\DeclareMathOperator{\Spec}{Spec}

\DeclareMathOperator{\len}{len}

\DeclareMathOperator{\im}{Im}

\DeclareMathOperator{\chrr}{char}
\DeclareMathOperator{\Skl}{Skl}

\DeclareMathOperator{\GK}{GKdim}
\DeclareMathOperator{\gr}{gr}
\DeclareMathOperator{\WDiv}{WDiv}
\DeclareMathOperator{\Bs}{Bs}

\DeclareMathOperator{\Cosupp}{Cosupp}

\newcommand{\sh}{\mathcal}
\newcommand{\mf}{\mathfrak}

\DeclareMathOperator{\shTor}{\mathcal{T}\!\mathit{or}}
\renewcommand{\Lsh}{\mathcal{L}}
\newcommand{\Ish}{\mathcal{I}}
\newcommand{\struct}{\mathcal{O}}
 \newcommand{\kk}{{\Bbbk}}
\newcommand{\dra}{\dashrightarrow}

\newcommand{\ZZ}{{\mathbb Z}}
\newcommand{\zed}{{\mathbb Z}}
\newcommand{\PP}{{\mathbb P}}
\newcommand{\NN}{{\mathbb N}}

\newcommand{\mb}{\mathbb}
\newcommand{\CC}{{\mb C}}

\newcommand{\I}{{\mathbb I}}

\newcommand{\ver}[1]{^{(#1)}}

\DeclareMathOperator{\Supp}{Supp}

\newcommand{\surfdata}{(X, \Lsh, \sigma, \sh{A}, \sh{D}, \sh{C}, \Omega)}

\newcommand{\surfD}{{\mathbb D}}
\newcommand{\Ea}{E_{\mf{a}}}
\newcommand{\Ed}{E_{\mf{d}}}
\newcommand{\Ec}{E_{\mf{c}}}
\renewcommand{\div}{{\rm div}}
\DeclareMathOperator{\NS}{NS}
\newcommand{\adm}{normal}

\title{Classifying birationally commutative projective surfaces}
\author{S. J. Sierra}
\address{Mathematics Department,  Fine Hall, Washington Road, Princeton, NJ 08544}
\email{ssierra@princeton.edu}
\date{\today}

\begin{document}

\begin{abstract}
Let $R = \bigoplus_{n \geq 0} R_n$ be a noetherian connected graded domain of Gelfand-Kirillov dimension 3 over an uncountable algebraically closed field.  
Suppose that the graded quotient ring of $R$ is of the form
\[ Q_{\gr}(R) \cong K[z, z^{-1}; \sigma]\]
where $K$ is a field;  we say that $R$ is a {\em birationally commutative projective surface}.  We classify  birationally commutative projective surfaces and show that they fall into four families, parameterized by geometric data.  This generalizes work of Rogalski and Stafford on birationally commutative projective surfaces generated in degree 1; our proof techniques are quite different.
\end{abstract}

\maketitle
\tableofcontents

\section{Introduction}\label{INTRO}

The  classification of {\em noncommutative projective surfaces}--- connected $\NN$-graded noetherian domains of Gelfand-Kirillov dimension 3---has been one of the major motivating problems in noncommutative ring theory over the past two decades.  In this paper, we resolve an important special case of this problem, classifying {\em birationally commutative projective surfaces}.  

We begin with some terminology.  We work over an  uncountable algebraically closed field $\kk$.  An $\NN$-graded algebra $R$ is {\em connected graded} if $R_0 = \kk$.   Let $R$ be a connected $\NN$-graded noetherian domain.  Its graded quotient ring is then of the form
\[ Q_{\gr} (R) \cong D[z, z^{-1}; \sigma],\]
where $D$ is a division ring and $\sigma$ is an automorphism of $D$.  By abuse of terminology, we refer to $D$ as the {\em function field} of $R$.  If  $D$ is a (commutative) field, we say that $R$ is {\em birationally commutative}.  If $R$ has GK-dimension 3, we say that $R$ is a {\em birationally commutative projective surface}.  

In this paper, we  classify  birationally commutative projective surfaces up to taking Veronese subrings.  We show that they occur in four families, and give geometric data  defining them.  Our results generalize (and give new proofs of) results of Rogalski and Stafford \cite{RS} on birationally commutative projective surfaces generated in degree 1.

Birationally commutative algebras have been fruitful objects of study in noncommutative ring theory for a number of years.  The fundamental example is the {\em twisted homogeneous coordinate ring} construction of Artin and Van den Bergh \cite{AV}.  Let $X$ be a projective variety, let $\sigma \in \Aut(X)$, and let $\Lsh$ be an appropriately positive invertible sheaf on $X$.  (The technical term is that $\Lsh$ is {\em $\sigma$-ample}; for a more precise definition, see Section~\ref{BACKGROUND}.)  For notational purposes, let
\[ \Lsh^{\sigma}:= \sigma^* \Lsh,\]
and let
\[ \Lsh_n := \Lsh \otimes \Lsh^{\sigma} \otimes \cdots \otimes \Lsh^{\sigma^{n-1}}.\]
The twisted homogeneous coordinate ring $B(X, \Lsh, \sigma)$ is defined as the section algebra of the twisted powers $\Lsh_n$ of $\Lsh$.  That is, 
\[ B(X, \Lsh, \sigma): = \bigoplus_{n \geq 0} H^0(X, \Lsh_n).\]
Multiplication on $B(X, \Lsh, \sigma)$ is induced from the product
\[ \Lsh_n \otimes \Lsh_m \to \Lsh_n \otimes \Lsh_m^{\sigma^n} \overset{\cong}{\to} \Lsh_{n+m}.\]
  For a concrete example, let $\tau: \PP^1 \to \PP^1$ be defined by  $\tau([x:y]) = [x: x+y]$.  Then $B(\PP^1, \struct(1), \tau)$ is isomorphic to the {\em Jordan plane}
\[ \kk_J[x,y] := \kk \ang{x,y}/(xy-yx-x^2).\]
By \cite[Theorem~1.2]{Keeler2000}, if $\Lsh$ is ample and $\sigma\in \Aut^o(X)$, then $B(X, \Lsh, \sigma)$ is a birationally commutative noetherian domain of GK-dimension $\dim X + 1$.

Not all noetherian birationally commutative domains are twisted homogeneous coordinate rings.   For example, let   $X, \Lsh, \sigma$ be as above, and let $Z$ be a proper subscheme of $X$, defined by an ideal sheaf $\sh{I}$.  In \cite{S-idealizer}, the author studied the ring
\beq\label{idealizer}
 R:= R(X, \Lsh, \sigma, Z): = \kk \oplus \bigoplus_{n \geq 1} H^0(X, \sh{I}\Lsh_n) \subseteq B(X, \Lsh, \sigma).
\eeq
Note that $R_{\geq 1}$ is a right ideal of $B(X, \Lsh, \sigma)$.  Under relatively mild conditions on $\sigma$ and $Z$, $R$ is the {\em idealizer} of this right ideal:  the maximal subring of $B(X, \Lsh, \sigma)$ in which $R_{\geq 1}$ is a two-sided ideal.  

When $R(X, \Lsh, \sigma, Z)$ is noetherian is understood.  It depends on a geometric condition on the motion of $Z$ under $\sigma$, known as {\em critical transversality}: roughly speaking, $\sigma$ and $Z$ should be in general position.  
\begin{defn}\label{idef-CT}
Let $X$ be a variety, let $\sigma \in \Aut(X)$, and let $Z$ be a closed subscheme of $X$.   Let $A \subseteq \ZZ$.  The set $\{\sigma^n(Z)\}_{n \in A}$ is {\em critically transverse} if, for all closed subschemes $Y$ of $X$, we have for all but finitely many $n \in A$ that 
\[ \shTor^X_i(\struct_Y, \struct_{\sigma^n(Z)}) = 0\]
for $i > 0$.
\end{defn}
If $Z$ is a point, then $\{\sigma^n Z\}$ is critically transverse exactly when it is {\em critically dense}:  it is infinite, and any infinite subset is Zariski-dense.

Recall:
\begin{theorem}\label{ithm-idealizer}
{\em \cite[Theorem~1.8]{S-idealizer})}
If $\{ \sigma^n Z\}$ is critically transverse, then the idealizer $R(X, \Lsh, \sigma, Z)$ is noetherian.\qed
\end{theorem}
For example, let $B:= B(\PP^1, \struct(1), \tau)$ as above, and assume that $\chrr \kk = 0$.  Then the subring
\[ \kk + y B \]
is equal to $R(\PP^1, \struct(1), \tau, [1:0])$.  That it is noetherian follows from  Theorem~\ref{ithm-idealizer}; this  result is originally due to Stafford and Zhang \cite[Example~0.1]{SZ}.  

By Artin and Stafford's remarkable classification of noncommutative projective curves, we have now seen all the classes of noetherian connected graded domains of GK-dimension 2.  
\begin{theorem}\label{ithm-AS}
{\em (\cite{AS})}
Let $R$ be a noetherian connected $\NN$-graded domain of GK-dimension 2.  Then some Veronese subring of $R$
   is either:

$(1)$ a twisted homogeneous coordinate ring $B(X, \Lsh, \sigma)$ for some  projective curve $X$,  automorphism $\sigma$  of $X$, and $\sigma$-ample invertible sheaf  $\Lsh$ on $X$; or

$(2)$ an  idealizer $R(X, \Lsh, \sigma, Z)$, where $X$, $\Lsh, \sigma$ are as above, and $Z$ is a proper subscheme of $X$, supported at points of infinite order. \qed
\end{theorem}

In this paper, we consider noetherian finitely $\NN$-graded domains  of GK-dimension 3, and  classify those that are birationally commutative.  We show that these also depend on geometric data, although, not surprisingly, the data is more complicated than  that of Theorem~\ref{ithm-AS}.
Besides twisted homogeneous coordinate rings and idealizers in them, another large class of birationally commutative projective surfaces are the {\em na\"ive blowup algebras}, as constructed in \cite{KRS} (and generalized in \cite{RS-0}).

To construct a  na\"ive blowup algebra, begin as usual with a projective surface $X$, an automorphism $\sigma$ of $X$, and a $\sigma$-ample invertible sheaf $\Lsh$ on $X$.   Also choose a point $P \in X$ (or more generally, let $P$ be  a 0-dimensional subscheme of $X$).  Let $\sh{I} = \sh{I}_P$ be the ideal sheaf of $P$.  Then we may form a ring
\beq\label{display-NB} S(X, \Lsh, \sigma, P) := \bigoplus_{n \geq 0} H^0(X, \Ish \Ish^\sigma \cdots \Ish^{\sigma^{n-1}} \Lsh_n),
\eeq
which we refer to as a {\em na\"ive blowup of $X$ at $P$}.     The construction of $S(X, \Lsh, \sigma, P)$ mimics the construction of a commutative blowup as a Rees ring:  we are taking (sections of) higher and higher successive powers of the ideal defining $P$, using the multiplication in the twisted homogeneous coordinate ring $B(X, \Lsh, \sigma)$.   

As with idealizers, the properties of the na\"ive blowup $S(X, \Lsh, \sigma, P)$ depend on the geometry of the orbits $\{ \sigma^n P\}_{n \in \ZZ}$.  In particular, the algebra is noetherian exactly when this set is critically transverse.  

\begin{theorem}\label{thm-blowup}
{\em (\cite[Theorem~1.1]{RS-0}, \cite[Proposition~4.12]{S-surfprop})}
Let $X$ be a  projective variety, let $\sigma \in \Aut(X)$, and let  $\Lsh$ be a $\sigma$-ample  invertible sheaf on $X$.  Let $P$ be a 0-dimensional closed subscheme of $ X$.  The ring $S(X, \Lsh, \sigma, P)$ is noetherian if and only if all points of $P$ have critically dense orbits. \qed
\end{theorem}

 Rogalski and Stafford \cite{RS} have  classified all birationally commutative projective surfaces that are generated in degree 1.  It turns out that   twisted homogeneous coordinate rings and na\"ive blowups are the only two types of rings that  occur.  
 \begin{theorem}\label{ithm-RS}
 {\em (\cite[Theorem~1.1]{RS})}
If $R$ is  a birationally  commutative projective surface that   is generated in degree 1, then  up to a finite-dimensional vector space $R$ is  either:

$(1)$ the twisted homogeneous coordinate ring of a projective surface; or 

$(2)$ the  na\"ive blowup  of  a projective surface at a 0-dimensional subscheme supported on points that move in  critically dense orbits. \qed
\end{theorem}

We note that Rogalski and Stafford consider a slightly more general class of rings than our noncommutative projective surfaces.  They study finitely $\NN$-graded noetherian domains $R$  whose graded quotient ring is of the form 
\[ K [z, z^{-1}; \sigma] \]
where $K = \kk(X)$ is the function field of a projective surface $X$ such that $\sigma$ induces an automorphism of $X$.  By \cite[Theorem~1.1]{R-GK}, any such $R$ has GK-dimension 3 or 5, and any birationally commutative domain of GK-dimension 3 that is generated in degree 1 is of the form considered in \cite{RS}.     See Section~\ref{FINAL} for more discussion of the GK-dimension of noncommutative surfaces.  

The hypothesis in Theorem~\ref{ithm-RS} that $R$ be generated in degree 1 seems overly restrictive; note that in contrast with the commutative case, there are many noncommutative noetherian graded rings that have no Veronese subring generated in degree 1.  (For example, the idealizers in \eqref{idealizer} have this property.)
We are able to remove this restriction and give a  complete classification of birationally commutative surfaces, using methods that are quite different from the proof of Theorem~\ref{ithm-RS}. 
Besides idealizers, na\"ive blowups, and twisted homogeneous coordinate rings, one new type of ring arises; we refer to these as {\em ADC rings}.  These are a generalized form of na\"ive blowup, with many similar properties and some important differences --- notably, these algebras may not have any Veronese generated in degree 1.  (For more discussion of ADC rings, see \cite{S-surfprop}.)

We obtain:

\begin{theorem}\label{ithm-main}
Let $R$ be a connected $\NN$-graded birationally commutative noetherian domain of GK-dimension 3.  
Then  
 some Veronese subring of $R$ is either:

$(1)$ the twisted homogeneous coordinate ring of a projective surface;

$(2)$ a na\"ive blowup or ADC ring on a projective surface; or

 $(1')$, $(2')$  an idealizer inside a ring of type $(1)$ or $(2)$.  

Further, all defining data is critically transverse.
\end{theorem}



We wish to state Theorem~\ref{ithm-main} more precisely.  To do so, we recall some terminology from \cite{S-surfprop}.  We say that the tuple $\surfD = \surfdata$ is {\em surface data} if:
\begin{itemize}
\item $X$ is  a projective  surface;
\item $\sigma$ is an automorphism of $X$;
\item $\Lsh$ is an invertible sheaf on $X$;
\item $\Omega$ is a curve on $X$;
\item $\sh{D}$ is the ideal sheaf of a 0-dimensional subscheme of $X$ such that all points in the cosupport of $\sh{D}$ have distinct infinite $\sigma$-orbits; and
\item $\sh{A}$ and $\sh{C}$ are ideal sheaves on $X$, cosupported on  0-dimensional subschemes consisting of points of infinite order, such that 
\beq\label{i-ADCeq}
\sh{I}_{\Omega} \sh{A} \sh{C} \subseteq \sh{D}. 
\eeq
\end{itemize}
Let $Z$ be the scheme defined by $ \sh{D}$, let $\Lambda$ be the scheme defined by $ \sh{A}$, and let $\Lambda'$ be the scheme defined by $\sh{C}$.  The surface data $\surfD$ is {\em transverse} if:
\begin{itemize}
\item $\Lsh$ is $\sigma$-ample; and
\item $\{\sigma^n Z\}_{n \in \ZZ}$, $\{ \sigma^n \Lambda\}_{n \geq 0}$, $\{ \sigma^n \Lambda'\}_{n \leq 0}$, and  $\{\sigma^n \Omega\}_{n \in \ZZ}$ are critically transverse.
\end{itemize}

Given surface data $\surfD= \surfdata$, we define sheaves $\sh{T}_n$ by setting $\sh{T}_0 := \struct_X$ and 
\[ \sh{T}_n := \sh{I}_{\Omega} \sh{A} \sh{D}^{\sigma} \cdots \sh{D}^{\sigma^{n-1}} \sh{C}^{\sigma^n}  \Lsh_n\]
for $n \geq 1$.  
We define an algebra $T(\surfD)$ associated to $\surfD$ by
\[ T(\surfD):= \bigoplus_{n\geq 0} H^0(X, \sh{T}_n).\]

Our main result is:
\begin{theorem}\label{ithm-surfclass}
Let $\kk$ be an uncountable algebraically closed field and let $R$ be a   birationally commutative projective surface over $\kk$.  Then  there is  transverse surface data $\surfD$ so that some Veronese of $R$ is isomorphic to $T(\surfD)$.  
\end{theorem}

We compare Theorem~\ref{ithm-main} and Theorem~\ref{ithm-surfclass}.  Let $\surfD= \surfdata$ be transverse surface data.  Let $\sh{A}' \supseteq \sh{I}_{\Omega}\sh{A}$ and $\sh{C}' \supseteq \sh{C}$ be maximal so that
\[ \sh{A}'\sh{C}' \subseteq \sh{D}.\]
Let 
\[\mb{E}:=(X, \Lsh, \sigma, \sh{A}', \sh{D}, \sh{C}', \emptyset).\]
The algebra $T(\mb{E})$ is called an {\em ADC ring}; of course, if $\sh{D} = \struct_X$ then $T(\mb{E})$ is just the twisted homogeneous coordinate ring $B(X, \Lsh, \sigma)$.  It turns out that $T(\surfD)$ is a left idealizer inside a right idealizer inside $T(\mb{E})$.

By enumerating the possible types of birationally commutative surfaces, Theorem~\ref{ithm-surfclass} has important applications.  For example, there are strong restrictions on the behavior of the Artin-Zhang $\chi$ conditions (see Section~\ref{FINAL}) for birationally commutative projective surfaces.  
\begin{theorem}\label{ithm-chi2}
Let $\kk$ be an uncountable algebraically closed field and let $R$ be a birationally commutative projective surface over $\kk$ that satisfies left or right $\chi_2$.  Then some Veronese subring of  $R$ is a twisted homogeneous coordinate ring,  and $R$ satisfies  $\chi$.
\end{theorem}
  
The main difficulty in proving both Theorem~\ref{ithm-surfclass} and Theorem~\ref{ithm-RS} is constructing the classical projective surface  $X$ that is associated to a given birationally commutative projective surface $R$.   Rogalski and Stafford prove Theorem~\ref{ithm-RS} through  a delicate analysis of a certain class of modules, called {\em point modules}  over  $R$.  This is quite difficult because for na\"ive blowups such modules are parameterized by an infinite series of projective schemes but not by any individual projective scheme; see \cite[Theorem~1.1]{KRS}.  In contrast, in the proof of Theorem~\ref{ithm-surfclass}  we construct the surface $X$ much more directly, through a method of successive approximations of the ``correct'' surface.  While there are technical issues involved in this proof, most of them are involved with showing that this method does, in fact, lead to an appropriate projective surface, and the actual construction is relatively straightforward.  

As mentioned, the classification given in Theorem~\ref{ithm-main} is  more complex  than the classification of curves in Theorem~\ref{ithm-AS}.  Note that a  connected $\NN$-graded domain of GK 2 is automatically, by Theorem~\ref{ithm-AS}, birationally commutative. For surfaces, this requires an additional assumption.  In fact, the full birational classification of noncommutative projective surfaces is still unknown, although there is a long-standing conjecture due to Artin.  Artin's conjecture may be phrased as:
\begin{conjecture} \label{conj-A}
{\em (\cite[Conjecture~4.1]{Artin1995})}
If $R$ is a noncommutative projective surface, then its function field is either:

$(1)$ a  field of transcendence degree 2 (birationally commutative);

$(2)$ a division ring finite-dimensional over a central field of transcendence degree 2 (birationally PI);

$(3)$ the full quotient division ring of  an Ore extension $K[x; \sigma, \delta]$,  where  $K$  is a field of transcendence degree 1  (a ``quantum ruled surface''); or

$(4)$ $D(E, \sigma)$, the function field of the Sklyanin algebra $\Skl_3(E, \sigma)$ for some elliptic curve $E$ and automorphism $\sigma$ of $E$ (a ``quantum rational surface'').
\end{conjecture}

Theorem~\ref{ithm-surfclass}  resolves the classification of noncommutative surfaces falling into case (1) of Artin's conjecture.   
We note that   Artin's original formulation of Conjecture~\ref{conj-A} assumed  much stronger technical  conditions on the rings under study; it is  interesting that these assumptions are not necessary to understand the birationally commutative case. 

Since Theorem~\ref{ithm-surfclass} is our main result, its proof takes up much of this paper.  We briefly summarize the organization of the paper here.  One of our main techniques will be to work, not with an algebra $R$, but with an associated quasicoherent sheaf $\sh{R}$ on some projective surface.  This sheaf is known as a {\em bimodule algebra}.  We give necessary background on bimodule algebras and some other technical terms in Section~\ref{BACKGROUND}.  In Sections~\ref{CODIM1}--\ref{DIM0} we construct surface data $\surfD=\surfdata$ associated to $R$.  This gives an approximation to $R$, and in fact we show that $T(\surfD)$ is a finite left and right module over a Veronese of $R$.  The surface $X$ is always normal, and $T(\surfD)$ is best thought of as some sort of  normalization of $R$.  In Section~\ref{DATA-CT}, we show that $\surfD$ is actually transverse.  In Section~\ref{FINDY}, we construct a finite birational morphism $\theta:  X \to Y$, an action of $\phi$ on $Y$ conjugate to $\sigma$, and a $\phi$-ample invertible sheaf on $Y$ so that (some Veronese of) $R$ is contained in $B(Y, \sh{M}, \phi)$.  In Section~\ref{SURJ} we prove a general result giving sufficient conditions for a subring of some algebra $T(\mb{E})$ to be equal to $T(\mb{E})$ in large degree.  Finally, in Section~\ref{FINAL} we show that there are ideal sheaves $\sh{E}$, $\sh{F}$, and $\sh{G}$ and a curve $\Phi$ on $Y$ so that 
 $\surfD'= (Y, \sh{M}, \phi, \sh{E}, \sh{F}, \sh{G}, \Phi)$ is transverse surface data and so that some Veronese of $R$ is isomorphic to $T(\surfD')$.  This proves  Theorem~\ref{ithm-surfclass} and Theorem~\ref{ithm-main};  we also prove Theorem~\ref{ithm-chi2}.

{\bf Acknowledgements}.  The results in this paper were part of my Ph. D. thesis at the University of Michigan, under the direction of J. T. Stafford.  During the writing of this paper, I was supported by NSF grants   DMS-0555750 and DMS-0802935.  I would like to thank Mel Hochster, Dan Rogalski, and Karen Smith for many informative conversations.

\section{Notation and background}\label{BACKGROUND}

Throughout, we let $\kk$ be a fixed uncountably algebraically closed field; all schemes are of finite type over $\kk$.

Given a birationally commutative projective surface $R$, we seek to produce associated geometric data.  One of our key techniques will be to work, not with a graded ring, but with a {\em bimodule algebra} on some projective scheme $X$.  
 In this section, we give the relevant definitions and notation for bimodule algebras.  Most of the material in this section was developed in \cite{VdB1996} and  \cite{AV}. Our presentation  follows that in \cite{KRS} and \cite{S-idealizer}.

A bimodule algebra on a variety $X$ is, roughly speaking, a quasicoherent sheaf with a multiplicative structure.  

\begin{defn}\label{def-bimod}
Let $X$ be a projective variety; that is, a projective integral separated scheme (of finite type over $\kk$).  An {\em $\struct_X$-bimodule}  is  a quasicoherent $\struct_{X \times X}$-module $\sh{F}$, such that for every coherent $\sh{F}' \subseteq \sh{F}$,  the projection maps $p_1, p_2: \Supp \sh{F}' \to X$ are both finite morphisms.    The left and right $\struct_X$-module structures associated to an $\struct_X$-bimodule $\sh{F}$ are defined respectively as $(p_1)_* \sh{F}$ and $ (p_2)_* \sh{F}$.  
We make the notational convention that when we refer to an $\struct_X$-bimodule simply as an $\struct_X$-module, we are using the left-handed structure (for example, when we refer to the global sections or higher cohomology
 of an $\struct_X$-bimodule).

There is a tensor product operation on the category of bimodules that has the expected properties; see \cite[Section~2]{VdB1996}.
\end{defn}

All the bimodules that we consider will be constructed from bimodules of the following form:
\begin{defn}\label{def-LR-structure}
Let $X$ be a variety and let $\sigma, \tau \in \Aut(X)$. Let $(\sigma, \tau)$ denote the map
\begin{align*}
X & \to X \times X \\
x & \mapsto (\sigma(x), \tau(x)).
\end{align*}
If  $\sh{F}$ is a quasicoherent sheaf on $X$, we define the $\struct_X$-bimodule ${}_{\sigma} \sh{F}_{\tau}$ to be 
\[ {}_{\sigma} \sh{F}_{\tau} = (\sigma, \tau)_* \sh{F}.\]
If $\sigma = 1$ is the identity, we will often omit it; thus we write $\sh{F}_{\tau}$ for ${}_1 \sh{F}_{\tau}$ and $\sh{F}$ for  the  $\struct_X$-bimodule ${}_1 \sh{F}_1 = \Delta_* \sh{F}$, where $\Delta: X \to X\times X$ is the diagonal.
\end{defn}

\begin{defn}\label{def-BMA}
Let $X$ be a variety and let $\sigma \in \Aut(X)$.  
An $\struct_X$-bimodule $\sh{R}$ is an {\em $\struct_X$-bimodule algebra} or  {\em bimodule algebra} if it is an algebra object in the category of bimodules.  That is, there are a unit map $1: \struct_X \to \sh{R}$ and a product map $\mu: \sh{R} \otimes \sh{R} \to \sh{R}$ that have the usual properties.
\end{defn}

We follow \cite{KRS} and define
\begin{defn}\label{def-gradedBMA}
Let $X$ be a variety and let $\sigma \in \Aut(X)$.  
A bimodule algebra $\sh{R}$ is a {\em graded $(\struct_X, \sigma)$-bimodule algebra} if:

(1) There are coherent sheaves $\sh{R}_n$ on $X$ such that
\[ \sh{R} = \bigoplus_{n \in \ZZ} {}_1(\sh{R}_n)_{\sigma^n};\]

(2)  $\sh{R}_0 = \struct_X $;

(3) the multiplication map $\mu$ is given by $\struct_X$-module maps
$\sh{R}_n \otimes \sh{R}_m^{\sigma^n} \to \sh{R}_{n+m}$, satisfying the obvious associativity conditions.  
\end{defn}

\begin{defn}\label{def-module}
Let $X$ be a variety and let $\sigma \in \Aut(X)$.  
Let $\sh{R}$ be a graded $(\struct_X, \sigma)$-bimodule algebra.  A {\em right $\sh{R}$-module} $\sh{M}$ is a an  $\struct_X$-bimodule $\sh{M}$ together with a right $\struct_X$-module map $\mu: \sh{M} \otimes \sh{R} \to \sh{M}$ satisfying the usual axioms.  We say that $\sh{M}$ is {\em graded} if there is a direct sum decomposition 
\[\sh{M} = \bigoplus_{n \in \ZZ}  (\sh{M}_n)_{\sigma^n}\]
 with multiplication giving a family of $\struct_X$-module maps $\sh{M}_n \otimes \sh{R}_m^{\sigma^n} \to \sh{M}_{n+m}$, obeying the appropriate axioms.
 
  We say that $\sh{M}$ is {\em coherent} if there are a coherent $\struct_X$-module $\sh{M}'$ and a surjective map $\sh{M}' \otimes \sh{R} \to \sh{M}$ of ungraded $\struct_X$-modules.  
We make similar definitions for left $\sh{R}$-modules.  The bimodule algebra $\sh{R}$ is {\em right (left) noetherian}
if every right (left) ideal of $\sh{R}$ is coherent.  By standard arguments, a graded $(\struct_X, \sigma)$-bimodule algebra is right (left) noetherian if and only if every graded right (left) ideal is coherent.
\end{defn}

If $\sh{R}$ is a graded $(\struct_X, \sigma)$-bimodule algebra, we may form its {\em section algebra}
\[ H^0(X, \sh{R})  = \bigoplus_{n\geq 0}H^0(X, \sh{R}_n).\]
Multiplication on $H^0 (X, \sh{R})$ is induced from the multiplication map $\mu$ on $\sh{R}$; that is, from the maps
\[ \xymatrix{
H^0(X, \sh{R}_n) \otimes H^0(X, \sh{R}_m) \ar[r]^{1\otimes \sigma^n} 
&  H^0(X, \sh{R}_n) \otimes H^0(X, \sh{R}_m^{\sigma^n}) \ar[r]^<<<<<{\mu} 
& H^0(X, \sh{R}_{n+m}).
}\]

We recall that ampleness is a key concept when working with bimodule algebras.  

\begin{defn}\label{def-ample}
Let $X$ be a projective variety, let $\sigma \in \Aut(X)$, and let  $\{\sh{R}_n\}_{n \in \NN}$ be a sequence of coherent sheaves  on $X$.  The sequence of bimodules $\{(\sh{R}_n)_{\sigma^n}\}_{n \in \NN}$ is {\em right ample} if for any coherent $\struct_X$-module $\sh{F}$, the following properties hold:
 
 (i)   $\sh{F} \otimes \sh{R}_n$ is globally generated for $n \gg 0$;
 
 (ii) $H^q(X, \sh{F} \otimes \sh{R}_n) = 0$  for $n \gg 0$ and all $q \geq 1$.  
 
\noindent The sequence $\{ (\sh{R}_n)_{\sigma^n}\}_{n \in \NN}$  is {\em left ample}  if for any coherent $\struct_X$-module $\sh{F}$, the following properties hold:
 
 (i)  $\sh{R}_n \otimes \sh{F}^{\sigma^n}$ is globally generated for $n \gg 0$;
 
 (ii) $H^q(X, \sh{R}_n \otimes \sh{F}^{\sigma^n}) = 0$  for $n \gg 0$ and all $q \geq 1$.

We say that an invertible  sheaf $\Lsh$ is {\em $\sigma$-ample} if the $\struct_X$-bimodules 
\[\{(\Lsh_n)_{\sigma^n}\}_{n \in \NN} = \{ \Lsh_\sigma^{\otimes n}\}_{n \in \NN}\]
 form a right ample sequence.  By  \cite[Theorem~1.2]{Keeler2000}, this is true if and only if the $\struct_X$-bimodules $\{(\Lsh_n)_{\sigma^n}\}_{n \in \NN}$ form a left ample sequence.
\end{defn}

The fundamental example of a bimodule algebra is as follows.  
\begin{example}
\label{eg-twist}
Let $X$ be a projective scheme, let $\sigma \in \Aut(X)$, and let $\Lsh$ be an invertible sheaf on $X$.  We define the {\em twisted bimodule algebra of $X$, $\Lsh$, and $\sigma$} to be
\[ \sh{B} := \sh{B}(X, \Lsh, \sigma) := \bigoplus_{n \geq 0} (\Lsh_n)_{\sigma^n}.\]
Then $\sh{B}$ is an $(\struct_X, \sigma)$-graded bimodule algebra.  
Taking global sections of $\sh{B}(X, \Lsh, \sigma)$ gives the twisted homogeneous coordinate ring $B(X, \Lsh, \sigma)$.
\end{example}
 All the bimodule algebras that we consider in this paper will be sub-bimodule algebras of twisted bimodule algebras.  

Let $\surfD = \surfdata$ be surface data, as defined in the Introduction.  Then one may form a bimodule algebra
\[ \sh{T}(\surfD):= \bigoplus_{n \geq 0} (\sh{T}_n)_{\sigma^n}\]
where $\sh{T}_0:=\struct_X$ and 
\[ \sh{T}_n:=\sh{I}_{\Omega} \sh{A} \sh{D}^{\sigma} \cdots \sh{D}^{\sigma^{n-1}} \sh{C}^{\sigma^n}  \Lsh_n.\]
The ring $T(\surfD)$ from the Introduction is just the section ring of $\sh{T}(\surfD)$.

The algebras $T(\surfD)$ are studied in the companion paper \cite{S-surfprop}, in slightly greater generality.  The surface data defined here is known in that paper as surface data of exponent 1.  We also caution that in \cite{S-surfprop}, we allow $\surfD$ to be transverse surface data if $\Lsh$ is not $\sigma$-ample.  We note that if $\surfD$ is transverse surface data, then $T(\surfD)$ may have GK-dimension 3 or 5, since twisted homogeneous coordinate rings on surfaces may have GK-dimension 3 or 5. 

 We recall the main result of \cite{S-surfprop}:
\begin{theorem}\label{thm-surfprop-main}
{\em (\cite[Theorem~4.13]{S-surfprop})}
Let $\surfD = \surfdata$ be  surface data, where $\Lsh$ is $\sigma$-ample.  Then the following are equivalent: 

$(1)$ $T(\surfD)$ is  noetherian;

$(2)$ $\sh{T}(\surfD)$ is noetherian;

$(3)$  $\surfD$ is transverse. \qed
\end{theorem}

We will also use the following technical result from \cite{S-surfprop}.
\begin{lemma}\label{lem-ample}
{\em (\cite[Lemma~3.8(2)]{S-surfprop})}
Let $\surfD = \surfdata$ be surface data, and let $\sh{T}:= \sh{T}(\surfD)$.  Suppose that $\Lsh$ is $\sigma$-ample, $\{\sigma^n \Omega\}$ is critically transverse, and that all points in the cosupport of $\sh{ADC}$ have dense orbits.  Then the sequence of bimodules $\{(\sh{T}_n)_{\sigma^n}\}$ is left and right ample. \qed
\end{lemma}

If $R = \bigoplus R_n$ is  a graded ring,  we denote the $k$'th Veronese of $R$ by $R\ver{k} := \bigoplus R_{nk}$; we use similar notation for graded bimodule algebras.

\section{Approximating birationally commutative surfaces in codimension 1}\label{CODIM1}

Let $R$ be a birationally commutative projective surface with function field $K$.   To prove Theorem~\ref{ithm-surfclass}, we must construct surface data $\surfD$
 that will correspond to $R$.     
We will approach the  data $\surfD$ through successive approximations.  In this section, we  construct a smooth surface $X$, an automorphism $\sigma$ of $X$, an effective curve $\Omega$ on $X$, and a line bundle $\Lsh$ on $X$ that capture some broad properties of the ring $R$.  We refer to this data as {\em approximating $R$ in codimension 1}.

 Since by assumption, the graded quotient ring of $R$ is isomorphic to $K [z, z^{-1}; \sigma]$, where $\sigma$ is a  $\kk$-automorphism  of $K$ and $K$ necessarily has transcendence degree 2, we begin with a birational class of surfaces and a  birational self-map of each surface. The central problem in constructing the surface data $\surfD$ is to find the correct surface $X$ within this birational class. This will occupy us for virtually all of the paper. 

 We say that $\sigma \in \Aut_{\kk}(K)$ is {\em geometric} if 
   there is a projective surface $X$ with $K \cong \kk(X)$ such that $\sigma$ is induced by an automorphism of $X$.   We call such a pair $(X, \sigma)$ a {\em model} for $R$.  Not all automorphisms of fields of transcendence degree 2 have models; for example, by \cite[Remark~7.3]{DF2001}, the automorphism $(x,y)\mapsto (x, xy)$ of $\mathbb{C}(x,y)$ is not geometric. 
  Fortunately, in our situation a  result of Rogalski ensures that  there is a model $(X, \sigma)$ for $R$.  Results of Artin and Van den Bergh then allow us to get precise information on the numeric action of  $\sigma$.   

We define some terminology.    
  If $X$ is a projective scheme, we denote the group of Cartier divisors on $X$ modulo numerical equivalence by $\NS(X)$.  We say that the automorphism $\sigma$ of $X$ is {\em numerically trivial} if the induced action of $\sigma$ on $\NS(X)$ is trivial; that is, if $\sigma D \equiv D$ for any Cartier divisor $D$ on $X$, where $\equiv$ denotes numerical equivalence.  
  We will say that an $\sigma$  is {\em quasi-trivial} if there is some integer $r > 0$ so that  $\sigma^r$  is numerically trivial.

  \begin{theorem}\label{thm-R-model}
  {\em (Rogalski, Artin-Van den Bergh)}
   Let $K/\kk$ be a finitely generated field extension of transcendence degree 2, and let 
   $\sigma \in \Aut_{\kk}(K)$.   Then every locally finite $\NN$-graded domain $R$  such that $Q_{\gr}(R) \cong  K [z, z^{-1};\sigma]$  has the same 
GK-dimension $d \in \{3, 4, 5, \infty\}$.   Moreover, 
$d \in  \{ 3, 5 \}$ if and only if $\sigma$ is geometric, and  $d = 3$ if and only if for any model $(X, \sigma)$ for $R$, the automorphism $\sigma$ is quasi-trivial.
\end{theorem} 

\begin{proof}
The first and second statements are \cite[Theorem~1.1]{R-GK}. Now suppose that $\sigma$ is geometric, and let $(X, \sigma)$ be a model for $R$.    Let $P \in O(\NS(X))$ be the matrix giving the numeric action of $\sigma$ on $\NS(X)$.  By \cite[Theorem~7.1]{R-GK} and \cite[Lemma~2.12]{R-GK}, all eigenvalues of $P$ have modulus 1;  now by \cite[Lemma~5.3]{AV}, the eigenvalues of $P$ are all roots of unity.     Let $\Lsh$ be an ample invertible sheaf on $X$.   Then \cite[Theorem~1.7]{AV} implies that $\Lsh $ is $\sigma$-ample, and that the GK-dimension of $B(X, \Lsh, \sigma)$, which is equal to $d$, is  3 if and only if $\sigma$ is quasi-trivial.  
\end{proof}

Suppose that $X$ and $X'$ are birationally equivalent surfaces; let $\sigma$, respectively $\sigma'$, be an automorphism of $X$, respectively $X'$.  We say that $\sigma$ and $\sigma'$ are {\em conjugate} if  there is a birational map $\pi: X' \to X$ so that $\pi \sigma' = \sigma \pi$ (as birational maps).
      Rogalski and Stafford note that it is an easy consequence of the existence of resolutions of singularities for surfaces (see \cite{Lipman1969}) that any geometric automorphism of a field of transcendence degree 2 is conjugate to an automorphism of a nonsingular surface.
   
   \begin{lemma}\label{lem-model}
   {\em (\cite[Lemma~6.2]{RS})}
If $K$ is a finitely generated field of transcendence degree 2 over $\kk$ and $\sigma \in \Aut_{\kk}(K)$ is a geometric automorphism of $K$, then there is a nonsingular surface $X$ with $\kk(X) = K$ such that $\sigma$ is induced from an automorphism of $X$.  In particular, if  a birationally commutative projective surface has a model, it has a nonsingular model.    
 \qed
\end{lemma}
  
  We now establish notation for some geometric data determined by $R$.  
  If $X$ is a projective variety and $V \subset K = \kk(X)$ is a finite-dimensional $\kk$-vector space, we will denote the coherent subsheaf of the constant sheaf $K$ on $X$ generated by the elements of $V$ by 
  \[ V \cdot \struct_X.  \]

Since $R$ has GK-dimension 3, Theorem~\ref{thm-R-model} implies that there is a model $(X, \sigma)$ for $R$.    By Lemma~\ref{lem-model}, we may also, if we choose, assume that $X$ is nonsingular.  
Replacing $R$ by an appropriate Veronese subring, we may assume that $R_1 \neq 0$.  

\begin{assnot}\label{assnot1}
We assume that  $R$ is a birationally commutative projective surface with $R_1 \neq 0$.    
 Let $K$ be the function field of $R$ and let $(X, \sigma)$ be a model for $R$. Fix $z \neq 0 \in R_1$.  
For all $n \geq 0$, we define $\bbar{R_n} := R_n \cdot z^{-n} \subset K$, so that
\[ R = \bigoplus_{n \geq 0} \bbar{R}_n z^n \subset K[z, z^{-1}; \sigma] .\]
 Let $\sh{R}_n(X) := \bbar{R}_n \cdot \struct_X$.
\end{assnot}

\begin{example}\label{eg-extended}
Before beginning to work with our noncommutative ring $R$, suppose for a moment that $R = \kk[x,y,z]$.  We know, of course, that $R \cong B(\PP^2, \struct(1), 1) = B(\PP^2, \struct(1))$ and that $\PP^2 = \Proj R$.  However, we cannot construct $\Proj R$ directly using noncommutative techniques.  Instead, we will   construct the defining data $(\PP^2, \struct(1))$ from the graded pieces of $R$.

The function field of $R$ is $K = \kk(x/z, y/z)$.  Consider the model  $X = \PP^1 \times \PP^1$  for $K$, where we think of $X$ as $\Proj$ of the bigraded ring $\kk[s,t][u,v]$.  We will let $s/t = x/z$ and $u/v = y/z$ in $K$.  

Let $\bbar{R}_1 := R_1 z^{-1} \subset K$.  Then 
\[ \bbar{R}_1 =  \kk \cdot  \frac{x}{z} \oplus \kk \cdot  \frac{y}{z} \oplus \kk 
 = \kk \cdot  \frac{sv}{tv} \oplus \kk \cdot \frac{tu}{tv} \oplus \kk \cdot \frac{tv}{tv}.  \]
Let $D$ be the divisor on $X$ defined by the equation $tv = 0$.  
On $X$, the rational functions in $\bbar{R}_1$ correspond to sections of $\struct_X(D) \cong \struct(1,1)$, and they generate
\[ \bbar{R}_1 \cdot \struct_X = \sh{I}_{[1:0] \times [1:0]} \struct_X(D).\]
We will modify $X$ by blowing up the base locus of $\bbar{R}_1 \subset H^0(X, \struct_X(D))$.

Let $\pi: \widetilde{X} \to X$ be the blowup of $X$ at $[1:0] \times [1:0]$.  Let $E = \pi^{-1}([1:0]\times[1:0]) $ be the exceptional locus of $\pi$, and let $F_1$ and $F_2$ be the strict transforms of $[1:0] \times \PP^1$ and $\PP^1 \times [1:0]$ respectively.  Then $F_1$, $F_2$ and $E$ are the three (-1) curves on $\widetilde{X}$, and on $\widetilde{X}$, the rational functions in $\bbar{R}_1$ generate the invertible sheaf 
\[
\Lsh := \bbar{R}_1\cdot \struct_{\widetilde{X}} = \struct_{\widetilde{X}}(F_1+F_2+E)\cong \sh{I}_E \struct_{\widetilde{X}}(\pi^* D).
\]  

One may check that $R$ and the section ring  $B(\widetilde{X}, \Lsh)$ are isomorphic.  However, $\Lsh$ is not ample.  By the Nakai-Moishezon criterion \cite[Theorem~V.1.10]{Ha}, the failure of ampleness of $\Lsh$ is equivalent to the existence of an effective curve $C$ so that $(F_1 +F_2+E).C \leq 0$.  One checks that $(F_1+F_2+E). F_1 = (F_1+F_2+E). F_2 = 0$.  That is, the curves $F_1$ and $F_2$ are contracted by the morphism defined by the base point free linear system $\bbar{R}_1 \subseteq H^0(\widetilde{X}, \struct_{\widetilde{X}}(F_1+F_2+E))$ on $\widetilde{X}$.  The image of $\widetilde{X}$ under this morphism is, of course, $\PP^2$, the ``correct'' model for $R$. 
\end{example}

We now return to the setting of a noncommutative projective surface $R$.  We assume Assumption-Notation~\ref{assnot1}.  It is immediate that the bimodule 
\[ \sh{R}(X) := \bigoplus_{n \geq 0} (\sh{R}_n(X))_{\sigma^n}\]
is in fact a graded $(\struct_X, \sigma)$-bimodule algebra, and of course $R \subseteq H^0(X, \sh{R}(X))$.  While ultimately we wish to understand $R$, our fundamental technique will be to approach $R$ by analyzing the bimodule algebra $\sh{R}(X)$ on a suitable model $(X, \sigma)$ for $R$; to construct $X$, we will mimic the steps carried out in Example~\ref{eg-extended}.

For all $n \geq 0$, let $\sh{R}_n := \sh{R}_n(X)$.  
We have
\beq\label{one}
\sh{R}_n \sh{R}_m^{\sigma^n} \subseteq \sh{R}_{n+m}
\eeq
for all $n,  m \geq 0$.  Since $R$ is an affine $\kk$-algebra, there is some $r \geq 1$ such that 
 for all $n > r$ we have:
\beq\label{two}
\sh{R}_n = \sum_{i=1}^r \sh{R}_i \sh{R}_{n-i}^{\sigma^i}.
\eeq

  We introduce some notation and terminology on divisors associated to finite-dimensional vector spaces of rational functions; see \cite[Chapter~1]{Laz} for a more detailed discussion.  

\begin{defn}\label{def-baselocus}
 If $X$ is a normal projective variety and $f$ is a rational function on $X$, we will denote its associated Weil  divisor by $\div_X (f)$, or $\div(f)$ if $X$ is understood.  We note that if $\sigma \in \Aut(X)$ and $f \in \kk(X)$, then $\div_X(f^\sigma) = \sigma^{-1} \div_X (f)$.  For any finite dimensional $\kk$-vector space $V \subseteq K$, and for any normal projective model $X$ for $K$, let $D^X(V)$ 
 be the minimal Weil divisor $D$ on $X$ such that $\div_X(f) + D \geq 0$ for all $f \in V$.   That is, $\struct_X(D^X(V))$ is canonically isomorphic to the double dual $ (V \cdot \struct_X)^{**}$.
 
 Now suppose that $X$ is an arbitrary projective variety and let $K = \kk(X)$.  Let $D$ be a Cartier divisor on $X$. 
We will denote  the  complete linear system $H^0(X, \struct_X(D))$ associated to $D$ by $\abs{D}$.  
Since $\struct_X(D)$ is by definition a subsheaf of the constant sheaf $K$, elements of $\abs{D}$ are rational functions on $X$.

Let $V \subset K$ be a finite-dimensional $\kk$-vector space.  Note that $V$ may be contained in many complete linear systems.  If $V \subseteq \abs{D}$ for some Cartier divisor $D$, we define the image of the natural map
\[ V \otimes \struct_X(-D) \to \struct_X\]
to be the {\em base ideal of $V$ with respect to $D$}.   The closed subscheme of $X$ that it defines is called the {\em base locus of $V$ with respect to $D$}.  We write it $\Bs^D(V)$.  If $(V \cdot \struct_X)^{**}$ is an invertible sheaf, then it corresponds to an effective Cartier divisor $D$ with $V \subseteq \abs{D}$.  This is the minimal such $D$, and in this situation we refer to the base ideal (respectively base locus) of $V$ with respect to $D$ simply as the {\em base ideal of $V$} (respectively, the {\em base locus of $V$}).   We write the base locus of $V$ as $\Bs(V)$.

If the base locus of the complete linear system $\abs{D}$ is empty, so $\struct_X(D)$ is globally generated, we say that $D$ and $\abs{D}$ are {\em base point free}.
\end{defn}

If $X$ is nonsingular or $X$ is normal and $D^X(V)$ is Cartier, then the base ideal and base locus of $V$ are always defined.  Note that in either situation, the base locus of $V$ has codimension at least 2.

  \begin{lemma}\label{lem-divisorlem}
  Let $X$ be a normal surface and let $K := \kk(X)$.  Let $\sigma \in \Aut(X)$, and let 
 $V, W \subseteq K$ be finite-dimensional $\kk$-vector spaces. 
 
 $(1)$ $D^X(VW)  = D^X(V) + D^X(W).$
 
 $(2)$  For every $n$, $D^X(V^{\sigma^n}) = \sigma^{-n} (D^X(V))$.
 \end{lemma}
 \begin{proof}
 (1) For any $f \in V$ and $g \in W$, we have
 \[ \div_X(fg) + D^X(V) + D^X(W) = \div_X(f) + \div_X(g) + D^X(V) + D^X(W) \geq 0,\]
 and so 
 \beq \label{a}
 D^X(V) + D^X(W) \geq D^X(VW).
 \eeq
 Now fix $f \in V$.  For any $g \in W$, we have 
 $ D^X(VW) + \div_X(f) + \div_X(g) \geq 0$.  
Thus
  $D^X(VW) + \div_X(f) \geq D^X(W)$.  As this holds for any $f \in V$, we obtain that
 \beq\label{b}
 D^X(VW) - D^X(W) \geq D^X(V).
 \eeq
 Combining \eqref{a} and \eqref{b}, we have proved (1).
 
 (2) is a consequence of the equality $\div_X(f^{\sigma}) = \sigma^{-1} \div_X(f)$.
 \end{proof}

Let $(X,\sigma)$ be a normal model for $R$.  We can now define data that determine $R$ in codimension 1, as follows.
 
 \begin{assnot}\label{assnot1a}
 Assume that  $R$ is a birationally commutative projective surface  with $R_1 \neq 0$.  
 Let $K$ be the function field of $R$ and let $(X, \sigma)$ be a normal model for $R$.  Fix $z \neq 0 \in R_1$.   Let  $\bbar{R}_n := R_n \cdot z^{-n}$  and  let $\sh{R}_n := \sh{R}_n(X) := \bbar{R}_n \cdot \struct_X$ for all $n \geq 0$.  
 
 For all $n \geq 0$, let $D_n := D^X(\sh{R}_n)$.  If $n <0$, let $D_n := 0$.   
  If $D_n$ is Cartier for all $n \geq 1$ (for example, if $X$ is nonsingular), then for $n \geq 1$ we further let $\sh{I}_n$ be the base ideal of $\bbar{R}_n$  and let $W_n$ be the base locus of $\bbar{R}_n$.  
 \end{assnot}


The following purely combinatorial lemma is a  restatement of results of  Artin and Stafford on the combinatorics of divisors on smooth curves.

\begin{lemma}\label{lem-AS}
{\em (Artin-Stafford)}
Let $A = \zed/(k)$ for some $k \in \zed$ (possibly $k= 0$).  Let $M$ be the free abelian group on the generating set $\{ P_i \st i \in A\}$; define a partial order $\geq$ on $M$ by saying that $E \geq 0$ if $E = \sum  n_i P_i$ where $n_i \geq 0$ for all $i$.  Define an automorphism $\sigma$ of $M$ by $\sigma(P_i) = P_{i+1}$.   

Suppose there is a sequence of elements $\{ E_i \st i \in \zed \}$ in $M$ satisfying:

{\rm (i)} $E_i \geq 0$ for all $i \geq 0$, and $E_i = 0$ if $i < 0$.

{\rm(ii)} There exists an integer $r$ such that 
\[E_n = \sup_{i=1}^r (E_i + \sigma^{-i} E_{n-i})\]
 for all $n \geq 1$.

Then:

$(1)$ If $k=0$, so $A = \zed$, there is an element $\Psi \geq 0 \in M$ and an integer $t \geq 0$ such that
\[
E_{m+n} = E_m + \sigma^{-m} (E_n) + \sigma^{-m}(\Psi)\]
for all $m, n \geq t$.

$(2)$  If  $k = 1$, then there is an integer $\ell$ so that $E_{n\ell} = n E_{\ell}$ for all $n \geq 1$.
\end{lemma}
\begin{proof}
(1) is \cite[Corollary~2.12]{AS}.  (2) is \cite[Lemma~2.7]{AS}.
\end{proof}

\begin{lemma} \label{lem-codim1}
Assume Assumption-Notation~\ref{assnot1a}.    
 Then 
 there are Weil divisors $0 \leq \Omega \leq D$ on $X$ and an integer $k \geq 1$ such that for all $n \geq 1$ we have 
 \beq\label{ge1}
 D_{kn} = D + \sigma^{-k} D + \cdots + \sigma^{-k(n-1)} D - \Omega,
 \eeq
  and so that no irreducible component of $\Omega$ is fixed by any power of $\sigma$. 
  Furthermore,
  \beq\label{ge2}
   D_{(n+m)k} = D_{nk} + \sigma^{-nk} (D_{mk}) + \sigma^{-nk}(\Omega)
   \eeq
  for all $n, m \geq 1$.
\end{lemma}

\begin{proof} 
We note that it suffices to prove the lemma for a Veronese subring of $R$; it then holds for $R$ by changing $k$ and $D$.  

   We claim that for all $n, m \geq 0$ we have 
\beq\label{foo1}
D_{n+m} \geq D_n + \sigma^{-n} D_m
\eeq
 and that there is $r \geq 1$ such that for all $n \geq 1$, we have   
 \beq\label{foo2}
 D_{n} = \sup_{i=1}^r \Bigl( D_i + \sigma^{-i} (D_{n-i}) \Bigr).
 \eeq

To see this, fix $m, n \geq 0$.  Let $D' := D^X(\bbar{R}_n (\bbar{R}_m)^{\sigma^n})$.   
By Lemma~\ref{lem-divisorlem},  $D' = D^X(\bbar{R}_n) + \sigma^{-n} D^X(\bbar{R}_m)$.
%
Because $\bbar{R}_n (\bbar{R}_m^{\sigma^n}) \subseteq \bbar{R}_{n+m}$, we have that $D_{n+m} \geq D'$.  This gives \eqref{foo1}.   Because $1 \in \bbar{R}_1$, we have $D_{n+1} \geq D_n$ for all $n$.  
Let $r \geq 1$ be such that for all $n \geq r$, we have  $R_n = \sum_{i=1}^r R_i R_{n-i}$.  Then \eqref{foo2} follows.

Let $\WDiv(X)$ denote the group of Weil divisors on  $X$.   Equation~\eqref{foo2} implies that there are only finitely many $\sigma$-orbits of prime divisors in $\WDiv(X)$ on which some $D_n$ is nonzero.   In particular, there are only finitely many such $\sigma$-orbits that are finite.  Thus for some $\ell$, 
 each $\sigma^{\ell}$-orbit of $\WDiv(X)$ on which some  $D_{n\ell}$ is nonzero is  either infinite or consists of one element.    Without loss of generality, we may replace $R$ by $R\ver{\ell}$ (which is still finitely generated) and assume that all curves of finite order that appear in some $D_n$ are $\sigma$-invariant. 

Let $E \in \WDiv(X)$ be a $\sigma$-invariant irreducible curve  such that some $D_n \geq E$.  There are only finitely many such $E$. 
Let $E_n := D_n |_E$.  Equations~\eqref{foo1} and \eqref{foo2} imply that $\{ E_n\}$ satisfies the hypotheses of Lemma~\ref{lem-AS}, with $k=1$.  Thus, 
by Lemma~\ref{lem-AS}(2), there is an  integer $m \geq 1$ such that for all $n \geq 1$, we have
\[
 D_{nm} |_E = n (D_m |_E).
\]
If $E \in \WDiv(X)$ is of finite order under $\sigma$ but not $\sigma$-invariant, then 
\[ D_{nm}|_E = 0 = n (D_m|_E)\]
for all $n$.  Thus, by replacing $R$ by $R \ver{m}$,  we may assume that 
\[ D_n |_E = n (D_1 |_E)\] for all  irreducible curves $E$ that are of finite order under $\sigma$ and for all $n \in \NN$. 

Let $\{P^1, \ldots, P^s\}$ be irreducible generators of the finitely many distinct infinite $\sigma$-orbits  in $\WDiv(X)$ on which some $D_n$ is nonzero.  Fix $1 \leq i \leq s$, and let $M$ be the subgroup of $\WDiv(X)$ generated by $\{\sigma^n(P^i)\}_{n \in \zed}$.  Let $E_n := D_n |_M$.
As before,  $\{E_n\}$ satisfies the hypotheses of Lemma~\ref{lem-AS}.  
Thus there exist $t$ and $\Psi$ as in the statement of Lemma~\ref{lem-AS}(1).  By varying $i$, we obtain integers $t^1, \ldots, t^s$ and divisors $\Psi^1, \ldots, \Psi^s$, with $\Psi^i$ supported on $\{ \sigma^n P^i\}_{n \in \ZZ}$.   Let $k = \max\{t^i\}$ and let $\Omega = \Psi^1 + \cdots  + \Psi^s$.  By construction, $\Omega$ contains no components of finite order under $\sigma$, and
\beq\label{keith}
 D_{m+n} = D_m + \sigma^{-m} (D_n) + \sigma^{-m} (\Omega)
 \eeq
 for all $n, m \geq k$.  
Define $D := D_k + \Omega$.  Then \eqref{ge2} holds for all $n, m \geq 1$, and an easy induction shows that
\eqref{ge1} holds 
 for all $n \geq 1$. 
  \end{proof}

 \begin{defn}\label{def-gap}
 Assume Assumption-Notation~\ref{assnot1}; in particular, fix $0 \neq z \in R_1$.  Let $(X, \sigma)$ be a normal model for $R$.  Let $D_n := D^X(\bbar{R}_n)$.  If there are effective Weil divisors $D$ and $\Omega$ on $X$ and an integer $k$ so that
 \eqref{ge1} and \eqref{ge2} hold for all $n, m \gg 0$, we follow the terminology of \cite{AS} and  say that $\Omega$ is {\em a gap divisor for $R$ on $X$ associated to $z$} (or more briefly a {\em gap divisor for $R$ on $X$}).  We say that $D$ is  {\em a coordinate divisor for $R$ on $X$ (associated to $z$)}.
     \end{defn}
 
 Note that $\Omega$ is a gap divisor for $R$ associated to $z$ if and only if it is a gap divisor associated to $z^n$ for some (or all) $R \ver{n}$.
Further, this gap divisor is unique (at least up to choice of $z$).  
 
 \begin{lemma}\label{lem-gap-unique}
 Assume Assumption-Notation~\ref{assnot1a}.   For a fixed $z \neq 0 \in \bbar{R}_1$, there is exactly one Weil divisor $\Omega$ that is a gap divisor for $R$ associated to $z$.
 \end{lemma}
 \begin{proof}
 By Lemma~\ref{lem-codim1}, there is a gap divisor $\Omega$ for $R$ on $X$ associated to $z$.  
 Suppose that there are Weil divisors $\Omega$ and $\Omega'$ so that for some $k,k' \geq 1$ we have
 \[ D_{k(n+m)} - D_{kn} - \sigma^{-kn}(D_{km}) = \sigma^{-kn}( \Omega)\]
 and
  \[ D_{k'(n+m)} - D_{k'n} - \sigma^{-k'n}(D_{k'm}) = \sigma^{-k'n}( \Omega')\]
for all $n, m \geq 1$.  Then 
\[ \sigma^{-kk'} (\Omega) = D_{2kk'} - D_{kk'} - \sigma^{-kk'}( D_{kk'}) = \sigma^{-k'k} (\Omega').\]
Thus $\Omega = \Omega'$.
\end{proof}

Initially, it will be more convenient to work on a nonsingular model for $R$.  
By Lemma~\ref{lem-codim1} and Theorem~\ref{thm-R-model}, we may replace $R$ by a Veronese subring to assume without loss of generality that we are in the following situation:
 
 \begin{assnot}\label{assnot2} 
  Assume that  $R$ is a birationally commutative projective surface  with $R_1 \neq 0$.  
 Let $K$ be the function field of $R$ and assume that there is a nonsingular model  $(X, \sigma)$  for $R$ so that $\sigma$ is numerically trivial.  As usual, we will identify Weil and Cartier divisors.  Fix $z \neq 0 \in R_1$.  
  Let  $\bbar{R}_n := R_n \cdot z^{-n}$  and  let $\sh{R}_n := \sh{R}_n(X) := \bbar{R}_n \cdot \struct_X$ for all $n \geq 0$.     
 For all $n \geq 0$, let $D_n := D^X(\sh{R}_n)$.  If $n <0$, let $D_n := 0$.   
Let $\sh{I}_n$ be the base ideal of $\bbar{R}_n$  and let $W_n$ be the base locus of $\bbar{R}_n$.  By definition, $W_n$ is 0-dimensional.

 Further assume that there are a gap divisor $\Omega$ and a coordinate divisor $D$ associated to $z$ so that \eqref{ge1} and \eqref{ge2} hold with $k =1$ for all $n \geq 1$, and that 
   $\Omega \cap \sigma^{k} \Omega$ is finite for all $k \neq 0$.

  Let $\Lsh := \struct_X(D)$.   Recall that $\Lsh_n = \Lsh \otimes \Lsh^{\sigma} \otimes \cdots \otimes \Lsh^{\sigma^{n-1}}$.     For $n \geq 0$, define 
  \[ \Delta_n := D + \cdots + \sigma^{-(n-1)}D,\]
  so $\Lsh_n  = \struct_X(\Delta_n)$.
  
  Our assumptions imply that  
  \[ \sh{R}_n= \sh{I}_n \sh{I}_{\Omega} \Lsh_n = \sh{I}_n(\Delta_n - \Omega)\]
  for all $n \geq 1$.  
    \end{assnot}
 
 We note that if Assumption-Notation~\ref{assnot2} holds for $R$, then it holds for any Veronese $R\ver{n}$ of $R$, by replacing $\sigma$ by $\sigma^n$.  (The effect of this change is to also replace $D$ by $\Delta_n$.)  Also note that in the setting of Assumption-Notation~\ref{assnot2}, we may regard $R$ as a subring of $B(X, \Lsh, \sigma) \cong \bigoplus H^0(X, \Lsh_n) z^n$, even if $\Lsh$ is not ample or $\sigma$-ample.  
 That is, elements of $\bbar{R}_n$ correspond to global sections of $\Lsh_n$.  We will make this identification throughout the rest of the paper.

\section{Points of finite order}\label{FINITEPOINTS}

The model $X$ that we chose for $R$ was picked  arbitrarily within its birational class.  To begin improving our approximation, we show in this section that we can  blow up  $X$ finitely many times to remove any points of finite order in the  base loci of the rational functions $\bbar{R}_n$.

\begin{lemma}\label{lem-finitedata}
Assume Assumption-Notation~\ref{assnot2}.  Then there is a finite set $V$ so that $W_n$ is supported on $ V \cup \cdots \sigma^{-(n-1)}(V)$ for all $n \geq 1$.  In fact, we may take 
 \beq\label{eqV}
 V = W_1 \cup W_2 \cup (\sigma^{-1} \Omega \cap \sigma^{-2} \Omega).
 \eeq
\end{lemma}

\begin{proof}
Recall that $\sh{I}_n$ is the base ideal of the vector space $\bbar{R}_n$ of rational functions.   
For all $m, n \geq 0$, the equation 
\[\sh{R}_n \sh{R}_m^{\sigma^n}  =  \sh{I}_{\Omega} \sh{I}_{n} \sh{I}_{\Omega}^{\sigma^n} \sh{I}_{m}^{\sigma^n} \Lsh_{n+m} \subseteq \sh{R}_{n+m} = \sh{I}_{\Omega} \sh{I}_{{n+m}}\Lsh_{n+m}\]
 gives  a set-theoretic containment 
\beq \label{kiri}
W_{n+m} \subseteq W_n \cup \sigma^{-n} (W_m) \cup \sigma^{-n} (\Omega).
\eeq
Define $V$ as in \eqref{eqV}.  Our assumptions imply that   $V$ is finite.

   Assume that for all $j \leq n$, we have $W_j \subseteq V \cup \cdots \cup \sigma^{-(j-1)} (V)$.  By construction, this is true for $n=2$.  For $n > 2$, it follows from \eqref{kiri}  that 
\[ W_{n+1}   \subseteq (W_1 \cup \sigma ^{-1} W_n \cup \sigma^{-1} \Omega) \cap
 (W_2 \cup \sigma ^{-2} W_{n-1} \cup \sigma^{-2} \Omega). \]
 By induction,  therefore, 
 \[ W_{n+1}  \subseteq V \cup \cdots \cup \sigma^{-n}V \cup (\sigma^{-1} \Omega \cap \sigma^{-2} \Omega)  = V \cup \cdots \cup \sigma ^{-n} V. \]
  \end{proof}
  
We give an elementary lemma on how base ideals transform under birational projective morphisms.
    
  \begin{lemma}\label{lem-baseideal}
  Let $\pi: X' \to X$ be a birational morphism of projective varieties, and let $D$ be an effective (Cartier) divisor on $X$.  Let $V \subseteq \abs{D}$.  Then the base ideal of $V$ on $X'$ with respect to $\pi^*D$ is the expansion to $X'$ of the base ideal of $V$ on $X$ with respect to $D$.  If $X$ and $X'$ are normal, $D^X(V)$ is Cartier, and the indeterminacy locus of $\pi^{-1}$ consists of smooth points of $X$, then  
  \[ \pi^* D^X(V) - D^{X'} (V)\]
  is effective and supported on the exceptional locus of $\pi$.  
  \end{lemma}
  \begin{proof}
  Note that the elements of $V$ are also elements of the linear system  $\abs{\pi^*D}$.  
  Let $\sh{I}$ be the base ideal of $V$ with respect to $D$.
  Let $\sh{J}$ be the base ideal of $ V$ on $X'$ with respect to $\pi^* D$; that is, the image of the natural map
  \[  V \otimes \struct_{X'}(- \pi^* D) \to \struct_{X'}.\]
 Pulling  back the surjection
  \[ V \otimes \struct_X(-D) \twoheadrightarrow \sh{I}\]
  to $X'$, we obtain a surjection
  \[  V \otimes \struct_{X'}(-\pi^* D) \twoheadrightarrow \pi^* \sh{I}.\]
 Composing this with the natural map from $\pi^* \sh{I} \to \pi^* \struct_X = \struct_{X'}$, we obtain the map
 \[  V \otimes \struct_{X'}(-\pi^* D) \to \struct_{X'}\]
 defining $\sh{J}$.  The image of $\pi^* \sh{I} \to \struct_{X'}$ is precisely $\sh{I} \struct_{X'}$; that is, $\sh{J}$ is the expansion of $\sh{I}$ to $\struct_{X'}$.
 
 Suppose now that $X$ and $X'$ are  normal, $D^X(V)$ is Cartier, and the indeterminacy locus of $\pi^{-1}$ consists of smooth points of $X$.  Let $F := D^X(V)$, and let $\sh{I}$ be the base ideal of $V$ on $X$.   Then by the above, 
 \[ V \cdot \struct_{X'} = \sh{I} \struct_{X'} (\pi^* F),\]
 and so $D^{X'}(V) = \pi^* F - C$ for some effective Weil divisor $C$ contained in the subscheme of $X'$ defined by $\sh{I} \struct_{X'}$.  Thus $C$ is supported on the exceptional locus of $\pi$.  
 \end{proof}

  We now begin the process of modifying $X$ to remove points of finite order from the base loci $W_n$.   We will do this through a series of blowups at finite orbits, and we begin by studying the effect of blowing up on the gap divisor $\Omega$.  
  
  \begin{lemma}\label{lem-blowup-omega}
  Assume Assumption-Notation~\ref{assnot2}; in particular, fix $0 \neq z \in R_1$, which we will use to calculate gap divisors, and let $\Omega$ be the gap divisor of $R$ on $X$.

$  (1)$ Let $\pi: \widetilde{X} \to X$ be the blowup of $X$ at a finite $\sigma$-orbit.  Then   the gap divisor of $R$ on $\widetilde{X}$ is the strict transform of $\Omega$.
  
  $(2) $ There are a nonsingular projective variety $X'$ and a birational morphism $\pi: X' \to X$ so that there is an automorphism $\sigma'$ of $X'$ with $\pi \sigma' = \sigma \pi$, and so that the gap divisor of $R$ on $X'$ contains no points of finite order under $\sigma'$.
   That is, by changing our smooth model $X$, without loss of generality we may assume that $\Omega$ contains no points of finite order.
  \end{lemma}
  \begin{proof}
  (1)  
By assumption, 
 \beq\label{gapeq1}
 D_n + \sigma^{-n} D_m + \sigma^{-n} \Omega = D_{n+m}
 \eeq
 for all $n, m \geq 1$.  
 For all $n \geq 1$, let $F_n := D^{\widetilde{X}}(\bbar{R}_n)$; let $\sh{J}_n$ be the base ideal of $\bbar{R}_n$ on $\widetilde{X}$, so 
 \[\sh{J}_n = \sh{R}_n(\widetilde{X})(-F_n) \subseteq \struct_{\widetilde{X}}.\] 
Note $\widetilde{X}$ is also nonsingular.  Let $\widetilde{\sigma}$ be the automorphism of $\widetilde{X}$ that lifts $\sigma$.   By Lemma~\ref{lem-codim1}, let $\widetilde{\Omega}$ be the gap divisor of $R$ on $\widetilde{X}$.   All components of $\widetilde{\Omega}$ are of infinite order under $\widetilde{\sigma}$ by Lemma~\ref{lem-codim1}, and there is some $k$ so that 
 \beq\label{gapeq2}
  F_{nk} + \widetilde{\sigma}^{-nk} F_{mk} + \widetilde{\sigma}^{-nk} \widetilde{\Omega} = F_{(n+m)k}
  \eeq
 for all $n, m \gg 0$.  
 
 For all $n \geq 0$, let $E_n := \pi^* D_n - F_n$.  By Lemma~\ref{lem-baseideal}, $E_n$ is effective and supported on the exceptional locus of $\pi$.   
Pulling back \eqref{gapeq1} to $\widetilde{X}$, we obtain that 
\[ \pi^* D_n + \widetilde{\sigma}^{-n} (\pi^* D_m) + \widetilde{\sigma}^{-n} (\pi^* \Omega) = \pi^* (D_{n+m})\]
for all $n, m \geq 1$.  Comparing this to \eqref{gapeq2}, we see that
\[ E_{nk} + \widetilde{\sigma}^{-nk} (E_{mk}) + \widetilde{\sigma}^{-nk}(\pi^* \Omega - \widetilde{\Omega}) = E_{(n+m)k}\]
for all $n, m \geq 0$.  Thus $\pi^* \Omega - \widetilde{\Omega}$ is supported on the exceptional locus of $\pi$.  All its components are thus of finite order under $\widetilde{\sigma}$; as $\widetilde{\Omega}$ contains no components of finite order under $\widetilde{\sigma}$, we see that $ \widetilde{\Omega}$ is the strict transform of $\Omega$.

(2) Suppose that $\Omega$ contains a point $p$ of finite order, and let $\pi: \widetilde{X} \to X$ be the blowup of $X$ at the orbit of $p$.   Let $\widetilde{\sigma}$ be the automorphism of $\widetilde{X}$ conjugate to $\sigma$.   Let  $\widetilde{\Omega}$ be the gap divisor of $R$ on $\widetilde{X}$.  By (1), $\widetilde{\Omega}$ is the strict transform of $\Omega$.    

Note that $\widetilde{\sigma}$ is still quasi-trivial.    Thus we may choose  $k$ so that $\widetilde{\sigma}^k$ is numerically trivial.  By assumption,  $\widetilde{\Omega} \cap \widetilde{\sigma}^{k} \widetilde{\Omega}$ is finite.  Then we have:
\beq\label{wombat}
 \Omega^2 > (\widetilde{\Omega})^2 = \widetilde{\Omega} . \widetilde{\sigma}^{k} (\widetilde{\Omega}) \geq 0.
 \eeq
If $\widetilde{\Omega}$ contains any points of finite order, we may repeat this process and reduce $(\widetilde{\Omega})^2$ further.  Since \eqref{wombat} shows that the gap divisor always has non-negative self-intersection, this process must terminate after finitely many steps.  That is, after finitely many steps we  obtain a gap divisor containing no points of finite order.
\end{proof}

 We are ready to prove that there is some model of $R$ for which all the sets $W_n$ --- that is, all the base loci of  $\bbar{R}_n$ --- consist of points of infinite order.  Before doing so, we recall some terminology from commutative algebra.  
  Let $(S, M)$ be a regular local ring of dimension 2, and let $I$ be an $M$-primary ideal of $S$.   Recall \cite[Section~12.1]{Eis} that the {\em Hilbert-Samuel function} of $S$ with respect to $I$ is defined as
  \[ H_I(n) = \len I^n/I^{n+1}.\]
  Recall further \cite[Exercise~12.6]{Eis} that the {\em multiplicity} of $I$, written $e(I)$, is defined as
  \[e(I) = (2!) \cdot ( \mbox{the leading coefficient of $H_I$)}.\]  
  This is  a positive integer that may be defined more geometrically as follows:  let $a, b \in  I$ be a regular sequence.  Then $e(I)$ is the intersection multiplicity of two general members of the ideal $aS+bS$.  
  
  Now let $X$ be a nonsingular surface and let $Z$ be a 0-dimensional subscheme of $X$.   We define the {\em multiplicity} $e(Z)$ of $Z$ to be the sum of the multiplicities of the defining ideal of $Z$ at all points in $\Supp (Z)$.  By definition, $e(Z) \geq 0$, and $e(Z) = 0$ if and only if $Z = \emptyset$.  
    Let $p \in Z$ and let $\pi: X' \to X$ be the blowup of $X$ at $p$; let $Z' \subset X'$ be the strict transform of $Z$.   The identification of $e(Z)$ with an intersection multiplicity shows that  $e(Z')$ is strictly less than $e(Z)$.

  \begin{proposition}\label{prop-blowup}
  Let $R$ be  a birationally commutative surface with $R_1 \neq 0$.  There is a smooth model $(X, \sigma)$ for some Veronese $R \ver{r}$ of $R$ so that  $\sigma$ is numerically trivial,  the gap divisor of $R \ver{r}$ on $X$ contains no points of finite order, and so that for all $n \geq1$ the base locus of $\bbar{R}_{nr}$ on $X$ is  supported on points of infinite order.  
    \end{proposition}
  \begin{proof}
  Choose a smooth model $(X, \sigma)$ for $R$.   By Lemma~\ref{lem-codim1}, by replacing $R$ by a Veronese subring, we may assume that we are in the situation of Assumption-Notation~\ref{assnot2}.  By Lemma~\ref{lem-blowup-omega}, by changing $X$ and possibly replacing $R$ by a further Veronese (to ensure that Assumption-Notation~\ref{assnot2} still holds), we may further assume that $\Omega$ contains no points of finite order.
    
Let $M$ be such that $R$ and $\sh{R} = \sh{R}(X)$ are generated in degrees $\leq M$.   If there is some  $1 \leq i \leq M$ such that $W_i$  contains a  point $p$ of finite order under $\sigma$, replace $X$ by the blowup of $X$ at the orbit of $p$.  As $e(W_i)$ is reduced each time, continuing finitely many times, we may assume that  there is a surface $\widetilde{X}$ with a morphism $\pi: \widetilde{X} \to X$ and an automorphism $\widetilde{\sigma}$ of $\widetilde{X}$,  so that for $i = 1 \ldots M$ the base locus of $\bbar{R}_i$ on $\widetilde{X}$ contains no points of finite order.  

One would hope that for $i>M$ the base locus of $\bbar{R}_i$ on $\widetilde{X}$ also contains no points of finite order.  This is unfortunately not necessarily true; however, by appealing to commutative algebra we can show that we can remove points of finite order by a further blowup.  We establish some notation.  
For all $n \geq 1$, let $F_n := D^{\widetilde{X}}(\bbar{R}_n)$, and let $\sh{J}_n := \sh{R}_n(\widetilde{X})(-F_n)$ be the base ideal of $\bbar{R}_n$ on $\widetilde{X}$.    We caution that \eqref{ge1} and \eqref{ge2} may not hold for the $F_i$ with $k = 1$, although they do, of course, hold for some $k$. 
By Lemma~\ref{lem-baseideal}, for all $n \geq 1$ the divisor   
\[E_n := \pi^* D_n - F_n = \pi^* \Delta_n - \pi^* \Omega - F_n\]
 is effective and supported on the exceptional locus of $\pi$.  In particular, all components of any $E_n$ are of finite order under $\sigma$.  
  By Lemma~\ref{lem-blowup-omega}(1), the gap divisor of $R$ on $\widetilde{X}$ is equal to $\pi^* \Omega$. 

The bimodule algebra $\sh{R}(\widetilde{X})$ on $\widetilde{X}$ is still generated in degrees $\leq M$.  That is, 
\beq\label{bar1}
\sh{J}_n(F_n) = \sh{R}_n(\widetilde{X}) =  \sum_{i=1}^M \sh{R}_i(\widetilde{X}) \sh{R}_{n-i}(\widetilde{X})^{\widetilde{\sigma}^i} 
 \eeq
for all $n \geq M$.    
We may rewrite \eqref{bar1} as
\begin{multline*}
 \sh{J}_n \struct_{\widetilde{X}}( - \pi^* \Omega -  E_n + \pi^* \Delta_n) 
= \\
\sum_{i=1}^M \sh{J}_i \struct_{\widetilde{X}} (- \pi^* \Omega - E_i + \pi^* \Delta_i) \cdot
 \sh{J}_{n-i}^{\widetilde{\sigma}^i} \struct_{\widetilde{X}}( \widetilde{\sigma}^{-i} (-\pi^* \Omega - E_{n-i} + \pi^* \Delta_{n-i}))
 \end{multline*}
 for all $n \geq M$.
 Subtracting $\pi^*(\Delta_n - \Omega)$ from both sides, we obtain 
 \beq \label{bar2}
 \sh{J}_n \sh{I}_{E_n}  = 
 \sh{J}_n \struct_{\widetilde{X}} (-E_n) = 
 \sum_{i=1}^M \sh{J}_i \sh{J}_{n-i}^{\widetilde{\sigma}^i} \struct_{\widetilde{X}} (-E_i - \widetilde{\sigma}^{-i}(E_{n-i} + \pi^* \Omega)).
 \eeq
 for all $n \geq M$.

For all $n$, let $\sh{K}_n$ be the minimal ideal sheaf on $\widetilde{X}$ that contains $\sh{J}_n$ and is cosupported at points on the exceptional locus of $\pi$; this exists because $\sh{J}_n$ is coartinian.   
Since $\Omega$ does not contain any fundamental points of $\pi^{-1}$, 
 by restricting \eqref{bar2} to the exceptional locus of $\pi$ we obtain 
\[ 
\sh{K}_n \sh{I}_{E_n} = 
\sum_{i=1}^M \sh{K}_i \sh{K}_{n-i}^{\widetilde{\sigma}^{i}} \struct_{\widetilde{X}}(-E_i - \widetilde{\sigma}^{-i} E_{n-i} )= 
\sum_{i=1}^M \sh{K}_i \sh{K}_{n-i}^{\widetilde{\sigma}^i} \sh{I}_{E_i} \sh{I}^{ \widetilde{\sigma}^{i}}_{E_{n-i}}
\]
for all $n \geq M$.  

For all $n$, let $\hat{\sh{K}}_n$ be the minimal ideal sheaf on $\widetilde{X}$ containing $\sh{K}_n$ and cosupported at points of finite order.   Now, there is some $k $ so that the ideal sheaf $\sum_{i=1}^M \sh{I}_{E_i} \sh{I}^{\widetilde{\sigma}^i}_{E_{n-i}}$ is  $\widetilde{\sigma}^{k}$-invariant,  as all $E_i$ are of finite order under $\widetilde{\sigma}$. This implies that
\[ 
\hat{\sh{K}}_n \sh{I}_{E_n} = \sum_{i=1}^M \hat{\sh{K}}_i \hat{\sh{K}}_{n-i}^{\widetilde{\sigma}^i} \sh{I}_{E_i } \sh{I}^{\widetilde{\sigma}^i}_{E_{n-i}}
\]
for all $n \geq M$.  But for $1 \leq i \leq M$, the base locus of $\bbar{R}_i$ on $\widetilde{X}$ contains no points of finite order, and so $\hat{\sh{K}}_i = \struct_{\widetilde{X}}$.  Thus
\beq\label{bar3}
\hat{\sh{K}}_n \sh{I}_{E_n} = \sum_{i=1}^M \hat{\sh{K}}_{n-i}^{\widetilde{\sigma}^i} \sh{I}_{E_i }\sh{I}^{\widetilde{\sigma}^i}_{E_{n-i}}\eeq
for $n \geq M$.  
Let $\ell$ be such that $\widetilde{\sigma}^{\ell}$ fixes all irreducible exceptional curves and so all components of $E_1, \ldots, E_M$.
   Then \eqref{bar3} and an easy induction imply that for all $n$, 
$\hat{\sh{K}}_{n}$ is $\widetilde{\sigma}^{\ell}$-invariant.  
%

Let $\sh{S}$ be the graded $(\struct_{\widetilde{X}}, \widetilde{\sigma}^{\ell})$-bimodule algebra defined by 
\[ \sh{S} := \bigoplus_{n \geq 0} (\sh{S}_n)_{\widetilde{\sigma}^{n\ell}},\]
where
\[ \sh{S}_n := \hat{\sh{K}}_{n\ell} \sh{I}_{E_{n\ell}}.\]
As all $\sh{S}_{n}$ are $\widetilde{\sigma}^{\ell}$-invariant, $ \sh{S}$ may be given the structure of a sheaf of (commutative) graded algebras (i.e., a commutative bimodule algebra) on $\widetilde{X}$. Call this bimodule algebra $\sh{S}_c$.  
 Since $R\ver{\ell}$ is finitely generated,  there is $N \geq 1$ so that
\[ \sh{R}_{n\ell}(\widetilde{X}) =\sum_{i=1}^N \sh{R}_{i\ell} (\widetilde{X}) \sh{R}_{(n-i) \ell}(\widetilde{X})^{\widetilde{\sigma}^{i\ell}}\]
for all $n \geq N$.  
Restricting to the exceptional locus of $\pi$ and to finite orbits, we obtain that
\[ \sh{S}_{n} =\sum_{i=1}^N \sh{S}_{i} \sh{S}_{n-i}^{\widetilde{\sigma}^{i\ell}}
= \sum_{i=1}^N \sh{S}_{i} \sh{S}_{n-i}
\]
for all $n \geq N$, and so $\sh{S} $ and $\sh{S}_c$ are finitely generated.  
An easy globalization of  \cite[Section~III.1.3, Proposition~3]{Bour}) shows that there is some $e \geq 1$ so that 
$\sh{S}_{c}\ver{e}$ is generated in degree 1.  Therefore, $\sh{S}\ver{e}$ is also generated in degree 1.

That is, 
\[ \hat{\sh{K}}_{n\ell e} \sh{I}_{E_{n \ell e}}=  \sh{S}_{ne} = (\sh{S}_e)^n= (\hat{\sh{K}}_{\ell e} \sh{I}_{E_{\ell e}})^n\]
for all $n \geq 1$.   As $\hat{\sh{K}}_{\ell e}$ is $\widetilde{\sigma}^{\ell}$-invariant, we may resolve it by a sequence of point blowups at finite $\widetilde{\sigma}$-orbits.  We obtain  a nonsingular surface $X'$ with a birational morphism $\psi: X' \to X$ so that $\widetilde{\sigma}$ is conjugate to an automorphism $\sigma'$ of $X'$ and so that the expansion of $\hat{\sh{K}}_{\ell e}$ to $X'$ is invertible.  Thus the expansion of $\hat{\sh{K}}_{n \ell e}$ to $X'$ is invertible for all $n \geq 1$. 

Recall that $\sh{J}_n$ is the base ideal of $\bbar{R}_n$ on $\widetilde{X}$.  For all $n$, there is an ideal sheaf $\sh{C}_n $  so that 
\[\sh{J}_n = \hat{\sh{K}}_n \sh{C}_n.\]
Necessarily, $\sh{C}_n$ is cosupported at points of infinite order.   Let $Z_n$ be the subscheme of $\widetilde{X}$ defined by $\sh{C}_n$.  
  Lemma~\ref{lem-baseideal} implies that 
 the base locus of $\bbar{R}_{n \ell e}$ on $X'$ is $\psi^{-1} (Z_{n \ell e})$, as the expansion of $\hat{\sh{K}}_{n \ell e}$ to $X'$ is invertible.   
 This 
  contains no points of finite order for any $n \geq 1$.  By Lemma~\ref{lem-blowup-omega}, the gap divisor of $R$ on $X'$ contains no points of finite order under $\sigma'$.  
  \end{proof}

 \section{An ample model for $R$}\label{AMPLEMODEL}

Let $(X, \sigma)$ be a   model for $R$.  
If a coordinate divisor  of $R$ on $X$ is $\sigma$-ample, we refer to $X$ or to the pair $(X,\sigma)$ as an {\em ample model} for $R$.   The goal of this section is to show that a normal  ample model for $R$ exists.

We begin by  giving the $\sigma$-twisted versions of some results about big and nef divisors.   Recall that we denote linear equivalence of divisors by $\sim$ and  numerical equivalence by $\equiv$.  If $D$ is  a divisor on $X$ and $m \geq 1$, let $\Delta_m := D + \sigma^{-1} D + \cdots + \sigma^{-(m-1)} D$.

Recall that a divisor $D$ on a projective surface $X$ is {\em big} if 
\[h^0(X, \struct_X(nD)) = \dim H^0(X, \struct_X(nD))\]
 grows as $O(n^2)$, and $D$ is {\em nef} if $D.C \geq 0$ for any curve $C$ on $X$.  We refer the reader to \cite{Laz} for the basic properties of big and nef divisors.  

\begin{defn}\label{def-big}
Let $\sigma$ be a quasi-trivial automorphism of the projective surface $X$.  We say that a divisor $D$ is {\em $\sigma$-big} if $h^0(X, \struct_X(\Delta_n))$ grows at least as $O(n^2)$.   
\end{defn}

We note that  if  $(X, \sigma)$ is a normal model for $R$ and $D$ is  a coordinate divisor for $R$ on $X$, then   $D$ is $\sigma$-big by assumption on the GK-dimension of $R$.  

\begin{lemma}\label{lem-kodaira}
{\em (Kodaira's Lemma; cf. \cite[Proposition~2.2.6]{Laz})}
Let $\sigma$ be a quasi-trivial automorphism of a smooth projective surface $X$, and let $D$ be a $\sigma$-big divisor on $X$.  Let $F$ be an effective divisor on $X$.  Then $H^0(X, \struct_X(\Delta_m - F)) \neq 0$ for all sufficiently large $m$.
\end{lemma}
\begin{proof}
We consider the exact sequence
\[  0 \to H^0(X, \struct_X(\Delta_m - F)) \to H^0(X, \struct_X(\Delta_m)) \stackrel{\phi_m}{\to} H^0(X, \struct_F(\Delta_m)). \]
As a  divisor on a smooth surface, $F$ is a local complete intersection.  
Thus Riemann-Roch \cite[Ex~IV.1.9]{Ha} implies that 
 there are constants $n$ and $c$ such that if $E$ is a divisor on $X$ with $E . F \geq n$, then
$h^0(X, \struct_F(E)) = E.F + c$.  Since $\sigma$ is quasi-trivial, $\Delta_m .F$ grows  no faster than  $O(m)$, and thus $h^0(X, \struct_F(\Delta_m))$ grows no faster than $O(m)$.  Since $D$ is $\sigma$-big, for $m \gg 0$ we have that $h^0(X, \struct_X(\Delta_m)) > h^0(X, \struct_F(\Delta_m))$ and therefore the map 
$\phi_m: H^0(X, \struct_X(\Delta_m)) \to H^0(X, \struct_F(\Delta_m))$ must have a kernel.  This gives a section of $\struct_X(\Delta_m -F)$.
\end{proof}

\begin{corollary}\label{cor-amp+eff}
{\em (cf. \cite[Corollary~2.2.7]{Laz})}
Let $\sigma$ be a quasi-trivial automorphism of the smooth projective surface $X$, and let $D$ be a $\sigma$-big divisor on $X$.  Let $A$ be an ample divisor on $X$.  Then there is some $m>0$ and some effective divisor $N$  on $X$ such that $\Delta_m \sim A +N$.
\end{corollary}
\begin{proof}
Choose $r$ such that $(r+1) A$ and $rA$ are both effective.  Using Lemma~\ref{lem-kodaira}, choose $m$ such that $H^0(X, \struct_X(\Delta_m - (r+1)A)) \neq 0$.  Thus there is some effective $N'$ with
\[\Delta_m - (r+1)A \sim N'.\]
  That is, $\Delta_m \sim A + (rA + N')$.  Since $rA$ and $N'$ are both effective, the theorem is proved for $N = rA + N'$.
\end{proof}

\begin{lemma}\label{lem-Wilson}
{\em (Wilson's Theorem; cf. \cite[Theorem~2.3.9]{Laz})}
Let $\sigma$ be a quasi-trivial automorphism of the smooth projective surface $X$, and let $D$ be a $\sigma$-big and nef divisor on $X$.
Then there are an effective divisor $N$ and an integer $m_0>0$ such that for every $m \geq m_0$, both $\Delta_m - N$ and  $\Delta_m - \sigma^{-(m -m_0)}(N)$ are base point free.
\end{lemma}
\begin{proof}
Fujita's Theorem \cite[Theorem~1.4.35, Remark~1.4.36]{Laz} gives a very ample divisor  $B$ such  that $H^i(X, \struct_X(B + P)) = 0$ for all nef $P$ and for all $i \geq 1$.  Note that the same property holds for all $\sigma^n B$.   Corollary~\ref{cor-amp+eff} implies that there is some $m_0>0$ so that $\Delta_{m_0} \sim 3B + N$ for $N$ effective.  Then for $m \geq m_0$, 
\[\Delta_m - N \sim 3B + \sigma^{-m_0} \Delta_{m - m_0}.\]
Now the cone of nef divisors is $\sigma$-invariant.  As $D$ is nef,  $\sigma^{-m_0}\Delta_{m-m_0}$ is also nef.  Thus
\[ \Delta_m - N  \sim B + 2B + \rm{nef}.\]
  Our assumption on $B$ implies  that 
\[ H^i(X, \struct_X(\Delta_m -N -iB) ) = 0\]
for $i = 1, 2$.   By \cite[Theorem~1.85]{Laz},  $\struct_X(\Delta_m - N)$ is globally generated.

Similarly, 
\[\Delta_m - \sigma^{-(m-m_0)} N \sim \Delta_{m-m_0} + 3 \sigma^{-(m-m_0)} B,\]
 and so $H^i(X, \struct_X(\Delta_m - \sigma^{-(m-m_0)}N - i \sigma^{-(m-m_0)} B) ) = 0 $ for $i = 1, 2$.  That is, $\Delta_m - \sigma^{-(m-m_0)}N$ is 0-regular with respect to the very ample divisor $\sigma^{-(m-m_0)} B$.  Applying \cite[Theorem~1.85]{Laz} again, $\struct_X(\Delta_m - \sigma^{-(m-m_0)}N)$ is globally generated.
\end{proof}

 We now begin to  consider the rational maps to projective space defined by the rational functions in $\bbar{R}_n$ and $\abs{\Delta_n}$.  We record here the elementary result that these are birational onto their image for $n \gg 0$.  
 
\begin{lemma}\label{lem-fieldgenerate}
Let $R$ be a birationally commutative projective surface with function field $K := \kk(X)$.  Assume that $R_1 \neq 0$ and fix $0 \neq z \in R_1$.  
For some $n$, the rational functions in $\bbar{R}_n$ generate $K$ as a field and so induce a birational map of $X$ onto its image.
\end{lemma}
\begin{proof}
Let $f_1, \ldots, f_k$ be rational functions that generate $K$.  For each $i$, there are homogeneous elements $a_i, b_i$ of some $R_{n_i}$ so that $f_i = a_i b_i^{-1}$.  By putting all the $f_i$ over a common denominator, we may assume that there are some $c_1, \ldots, c_k, b \in R_n$ with $f_i = c_i b^{-1}$ for all $i$.  Thus $\bbar{R}_n$ generates  $K$.
\end{proof}

By Proposition~\ref{prop-blowup} and Lemma~\ref{lem-fieldgenerate}, we may pass to a further Veronese subring to strengthen our assumptions on $R$.  

\begin{assnot}\label{assnot3a}
  Assume that  $R$ is a birationally commutative projective surface  with function field $K$ so that  $R_1\neq 0$.  Fix $0 \neq z \in R_1$, and define $\bbar{R}_n := R_n z^{-n}$.  Assume that $\bbar{R}_1$ generates $K$ as a field.
    
 Assume also  that there is a nonsingular model   $(X, \sigma)$  for $R$ so that $\sigma$ is numerically trivial.  Define $\sh{R}_n(X)$, $D_n$, $\sh{I}_n$, and $W_n$  as in Assumption-Notation~\ref{assnot1a}.
 Further assume that there are a gap divisor $\Omega$ and a coordinate divisor $D$ associated to $z$ so that \eqref{ge1} and \eqref{ge2} hold with $k =1$ for all $n, m \geq 1$, and that 
   $\Omega \cap \sigma^{k} \Omega$ is finite for all $k \neq 0$.  We further assume that $\Omega$ and all $W_n$ are disjoint from finite $\sigma$-orbits.
 
 We continue to define $\Delta_n := D   + \cdots + \sigma^{-(n-1)} D$ and $\Lsh := \struct_X(D)$.  
\end{assnot}

We remark that if Assumption-Notation~\ref{assnot3a} holds for $R$, it  holds for any Veronese $R \ver{k}$ of $R$, by replacing $\sigma$ by $\sigma^k$ and $D$ by $\Delta_k$.  

\begin{lemma} \label{lem-basics}
Assume Assumption-Notation~\ref{assnot3a}, so $D$ is the coordinate divisor of $R$ on $X$ and $\sigma$ is numerically trivial.  Then  $\Delta_n$ is big and nef for all $n \geq 1$.
\end{lemma}
\begin{proof}
It is enough to show that $D$ is big and nef.
  Suppose that $D$ is not nef.  Then there is some  effective curve $C$ such that $C.D < 0$.  For  $n \gg 0$ we have $(\Delta_n - \Omega) .C = nD.C - \Omega.C <0$.  This implies that if $\Gamma \sim \Delta_n - \Omega$ is effective, then $C \leq \Gamma$; that is, $C$ is contained in the base locus of $\abs{\Delta_n - \Omega}$.   But $\Bs(\abs{\Delta_n - \Omega}) \subseteq \Bs (\bbar{R}_n) = W_n$, and this is 0-dimensional.  
  Thus $D$ is nef.

%
 By assumption on $R$, we know that $D$ is $\sigma$-big.    By Corollary~\ref{cor-amp+eff}, for some $n \geq 1$ we have 
$\Delta_n \sim A + F$, where $A$ is ample and $F$ is effective. 
Thus $\Delta_n$ is big by \cite[Corollary 2.2.7]{Laz}.  Since $\sigma$ is numerically trivial and bigness is numeric \cite[Corollary~2.2.8]{Laz}, we see that $nD$ and therefore $D$ are big.
\end{proof}

\begin{theorem}\label{thm-BPF}
Assume Assumption-Notation~\ref{assnot3a}.  Then there is some $k$ so that $\Delta_{nk}$ is base point free for $n \gg 0$.
\end{theorem}

We note that if $R$ is commutative (so $\Omega = 0$ and $\sigma = \id_X$), then this follows from Zariski's result \cite[Theorem~6.2]{Zariski1962} that if $\Lsh$ is a line bundle on a projective variety with a 0-dimensional base locus, then some tensor power of $\Lsh$ is globally generated.

\begin{proof}
For all $n$, let $Z_n := \Bs^{\Delta_n}(\abs{\Delta_n})$.  We want to show that for some $k$, $Z_{nk} = \emptyset$ for $n \gg 0$.

We first show that $Z_n$ is 0-dimensional for $n \gg 0$.  Let $C_n$ be the 1-dimensional component of $Z_n$.  The coordinate  divisor $D$ is $\sigma$-big by assumption, and nef by Lemma~\ref{lem-basics}.  By Lemma~\ref{lem-Wilson}, we know that there is some effective $N$ such that for all $m \gg 0$, both $\Delta_m - N$ and $\Delta_m - \sigma^{-(m -m_0)} N$ are base point free.  Thus $C_m \subseteq N \cap \sigma^{m-m_0} N$ for all $m \gg 0$, and so for all  $m \gg 0$, $C_m$ is a union of components of $N$ that are of finite order under $\sigma$. 
Now, 
\[C_m \subseteq \Bs^{\Delta_m}(\abs{\Delta_m}) \subseteq \Omega \cup \Bs(\abs{\Delta_m - \Omega}).\]
 As  $\Bs(\abs{\Delta_m - \Omega}) \subseteq \Bs(\bbar{R}_n) = W_n$ is 0-dimensional, we also have that  $C_m \leq \Omega$.  Since $\Omega$ has no components of finite order,  $C_m = 0$ for $m \gg 0$.

By passing to a Veronese subring, and replacing $D $ by some $\Delta_k$ and $\sigma$ by $\sigma^k$, we may assume that $Z_n$ is 0-dimensional for all $n \geq 1$.  
Let $\phi = \phi_{\abs{D}}$ be the rational map from $X$ to some $\PP^N$ defined by the complete linear system $\abs{D}$.  Let $Y$ be the closure of $\phi(X)$; we will abuse notation and refer to $\phi$ as a rational map from $X$ to $Y$.  
Note that $\phi$ is birational by assumption, as $\bbar{R}_1 \subseteq H^0(X, \struct_X(D))$ generates $K$.

By blowing up the finite base locus of $\abs{D}$, we obtain a surface $X'$ and a diagram of birational maps
\[ \xymatrix{
	&X' \ar[ld]_{\pi} \ar[rd]^{\phi'}	& \\
	X \ar@{-->}[rr]_{\phi} && Y} \]
such that $\pi$ and $\phi'$ are morphisms.  
Let $C$ be a reduced and  irreducible hyperplane section of $Y$ that avoids the finitely many  points with positive-dimensional preimage in  $X$ or $X'$ and does not contain any component of the singular locus of $Y$.  Such $C$ exist by Bertini's theorem and \cite[Remark~III.7.9.1]{Ha}.   Then $D':=\pi (\phi')^{-1}(C) $ is a reduced and irreducible curve that is linearly equivalent to $D$.  Without loss of generality, we may replace $D$ by $D'$ and assume that $D$ is reduced and irreducible.

  We will show that $\struct_X(\Delta_m)$ is globally generated  for all $m \gg 0$.  The proof is based on repeated applications of the following long exact cohomology sequence.
  Let  $B$ be an effective divisor on $X$ and let $A$ and $A'$ be divisors such that $A' \sim A - B$.  Then the exact sequence
\[ 0 \to \struct_X(A') \to \struct_X(A) \to \struct_B(A) \to 0\]
induces  a long exact cohomology sequence
\begin{multline}\label{theo} 
0 \to H^0(X, \struct_X(A')) \to H^0(X, \struct_X(A)) \to H^0(X, \struct_B(A)) \\
\to H^1(X, \struct_X(A')) \to H^1(X, \struct_X(A)) \to H^1(X, \struct_B(A)).\end{multline}

In particular, for all $m \geq 0$ there are homomorphisms
\[H^1(X, \struct_X(\sigma^{-1}(\Delta_m))) \to H^1(X, \struct_X(\Delta_{m+1})) \to H^1(X, \struct_D(\Delta_{m+1})).\]
Now $D$ is irreducible and $D. \Delta_m = m D^2$, as $\sigma$ is numerically trivial.  Since $D$ is big and nef, $D^2 > 0$ by \cite[Theorem~2.2.16]{Laz}.  By Riemann-Roch, there is an integer $m_0$ such that if   $m \geq m_0$, then $H^1(D, \struct_D(\Delta_{m})) = 0$.  Thus if $m \geq m_0$ we have that $h^1(X, \struct_X(\Delta_m)) =  h^1(X, \struct_X(\sigma^{-1} \Delta_m) ) \geq h^1(X, \struct_X(\Delta_{m+1}))$.  Therefore, there is some $m_1 \geq m_0$ and some non-negative integer $a$ such that if $m \geq m_1$, we have that 
\[h^1(X, \struct_X(\Delta_m)) = a.\]
  Increasing $m_1$ further if necessary,  we may also assume that if $H$ is any divisor on $X$ with $D.H \geq m_1 D^2$, then   $\struct_{\sigma^jD}(H)$ is globally generated for any $j$ and $H^1(X, \struct_{\sigma^jD}(H)) = 0$.

Suppose that $m \geq 2m_1$.  We claim that $\struct_X(\Delta_m)$ is globally generated; that is, $\Bs(\abs{\Delta_m}) = 0$.  Since $\Bs(\abs{\Delta_m}) \subseteq D \cup \sigma^{-1}(D) \cup \cdots \cup \sigma^{-(m-1)}(D)$, it is enough to show that $\struct_X(\Delta_m)$ is globally generated at every point in $\sigma^{-i}(D)$ for $i =0 \ldots m-1$.

We claim that for any such $i$, we have that 
\beq\label{wannaget}
h^1(X, \struct_X(\Delta_m - \sigma^{-i}(D))) =a.
\eeq
   We will do the case when  $m = 2m'$ for $m' \geq m_1$ and $i \leq m'-1$; similar arguments work for other choices for $m$ and $i$.    For $j = 0 \ldots i$, let 
   \[C_j := \Delta_m - \sigma^{-i}(D) - \cdots - \sigma^{-j}(D).\]
      Define $C_{i+1} := \Delta_m$. Thus for $j = 0 \ldots i$, we have $C_{j} = C_{j+1} - \sigma^{-j}(D)$, 
and  $C_{j+1} \geq \sigma^{-m'} \Delta_{m'}$.  Thus $\sigma^{-j}D \cdot C_{j+1} \geq m_1D^2$ and so by the choice of  $m_1$ we have that 
   $H^1(X, \struct_{\sigma^{-j} D} (C_{j+1})) = 0$. 
     Thus the long exact cohomology sequence \eqref{theo} gives an exact sequence
\[  H^1(X, \struct_X(C_j)) \to H^1(X, \struct_X(C_{j+1})) \to H^1(X, \struct_{\sigma^{-j} D} (C_{j+1})) = 0. \]   
We obtain that
\[h^1(X, \struct_X(C_0)) \geq h^1(X, \struct_X(C_1)) \geq \cdots \geq h^1(X, \struct_X(C_i)) \geq h^1(X, \Delta_m) = a.\] 
Since $C_0 = \Delta_m - \Delta_{i+1} = \sigma^{-(i+1)}(\Delta_{m-i-1})$ and $m-i-1 \geq m' \geq  m_1$, we have that $h^1(X, \struct_X(C_0)) =a$.  Thus  $h^1(X, \struct_X(C_i)) = a$.  The claim \eqref{wannaget} is proved.

Now let $0 \leq i \leq m-1$ be arbitrary.  As a special case of \eqref{theo}, we obtain  the long exact sequence
\begin{multline*}
0 \to H^0(X, \struct_X (\Delta_m - \sigma^{-i}(D))) \to H^0(X, \struct_X (\Delta_m)) \stackrel{\phi}{ \to} H^0(X, \struct_{\sigma^{-i}(D)}(\Delta_m)) \to \\
H^1(X, \struct_X (\Delta_m - \sigma^{-i}(D))) \to H^1(X, \struct_X (\Delta_m)) \to H^1(X, \struct_{\sigma^{-i}(D)}(\Delta_m)).
\end{multline*}
By assumption on $m$, we have $H^1(X, \struct_{\sigma^{-i} D} (\Delta_m)) = 0$, and we have seen that 
 \[ h^1(\struct_X(\Delta_m - \sigma^{-i}(D))) = h^1(\struct_X(\Delta_m)) =a.\]
   Thus  the map
\[ \phi: H^0(X, \struct_X(\Delta_m)) \to H^0(X, \struct_{\sigma^{-i}(D)}(\Delta_m))\]
is surjective.  Since  we have taken $m$ sufficiently large so that $\struct_{\sigma^{-i}(D)}(\Delta_m)$ is globally generated,  $\Bs(\abs{\Delta_m})$ must be disjoint from $\sigma^{-i}(D)$.  Since this holds for all $i$, we see that $\abs{\Delta_m}$ is base point free.
\end{proof}

We are almost ready to construct the ample model for $R$.  We first prove two lemmas about birational maps.  

\begin{lemma}\label{lem-contractedcurves}
Let $D$ be a Cartier divisor on a possibly singular projective variety $X$ and let $V \subseteq \abs{D}$ be a subspace of dimension $d \geq 2$.  Let $\phi = \phi_V$ be the rational map to $\PP^{d-1}$ defined by $V$, and let ${\Gamma}$ be an irreducible curve on $X$ that is disjoint from the base locus of $V$ with respect to $D$.    Then $\phi$ contracts ${\Gamma}$ if and only if $D \cdot {\Gamma} = 0$.  Further, if $\phi$ contracts ${\Gamma}$, then for any $v \in V$, either $v$ never vanishes on ${\Gamma}$ or $v |_{\Gamma} \equiv 0$.
\end{lemma}
\begin{proof}
Suppose that $\phi$ contracts ${\Gamma}$ to a point.  By making a linear change of coordinates, without loss of generality we may assume that  $\phi({\Gamma}) = [1:0:\cdots:0]$.  This is the same as choosing a basis $\{v_1, \ldots, v_d\}$ of $V$ such that $v_1 |_{\Gamma}$ is never 0  and that $v_i |_{\Gamma} \equiv 0$ for all $i \geq 2$.  In particular, the divisor of zeroes of $v_1$ is disjoint from ${\Gamma}$ and so $D.{\Gamma} = 0$.

Conversely, suppose that $D . {\Gamma} = 0$.  Then choose distinct points $x, y \in {\Gamma}$ and $v \in V$.  If $v(x) \neq 0$ but $v(y) =0$ then we have that ${\Gamma}.D > 0$; thus $v(y) = 0$  if and only if $v|_{\Gamma} \equiv 0$.  Now, since ${\Gamma}$ does not meet $\Bs^{D}(V)$, there is some $v \in V$ such that $v(x) \neq 0$.  We may choose a basis $\{v, v_2 \ldots, v_d\}$ for $V$ such that $v_i(x) = 0$ for all $i \geq 2$.  By the above, in these coordinates  $\phi({\Gamma}) = [1:0:\cdots:0]$.
\end{proof}

We obtain as a corollary that any curve $\Gamma$ such that $\Gamma. \Delta_n = 0$ must be disjoint from the gap divisor $\Omega$ and from the base loci $W_m$.  
\begin{corollary}\label{cor-noprob}
Assume Assumption-Notation~\ref{assnot3a}.  
Suppose that $\abs{\Delta_n}$ is base point free.  Let $\phi_n$ be the morphism to projective space defined by $\abs{\Delta_n}$.  If $\phi_n$ contracts an irreducible curve ${\Gamma}$, then 
there is some $f \in \bbar{R}_n$ so that 
\[(\div_X(f) + \Delta_n) \cap {\Gamma} = (\div_X(f) + \Omega + (\Delta_n - \Omega)) \cap {\Gamma} = \emptyset. \]
 In particular, 
 $W_n \cup \Omega$ is disjoint from ${\Gamma}$.
\end{corollary}
\begin{proof}
As $X$ is nonsingular, we may identify Cartier and Weil divisors.  
Lemma~\ref{lem-contractedcurves} implies that the set of irreducible curves contracted by $\phi_n$ is precisely the set of irreducible curves ${\Gamma}$ with ${\Gamma} . \Delta_n = 0$.  As $\sigma$ is numerically trivial,  ${\Gamma} . \Delta_n = 0 $ if and only if $\sigma {\Gamma} . \Delta_n = 0$.  Thus the  set of curves contracted by the morphism $\phi_n$ 
is $\sigma$-invariant.  

By assumption $\phi_n$ is birational onto its image.  Thus  there are finitely many such curves and so all are of finite order under $\sigma$.  In particular, if ${\Gamma}$ is such a curve, then ${\Gamma} \not\leq \Omega$.

Fix an irreducible curve ${\Gamma}$ with $\Delta_n . {\Gamma} = 0$.  As ${\Gamma} \not \subseteq \Omega \cup W_n$, we have some $f \in \bbar{R}_n $ so that ${\Gamma} \not \leq \div_X(f) + \Delta_n$.  Thus ${\Gamma} \cap \div_X(f) + \Delta_n = \emptyset$ by Lemma~\ref{lem-contractedcurves}. 
Now, set-theoretically we have
\[
\Omega \cup W_n  = \bigcap_{g \in \bbar{R}_n} \div_X(g) + \Delta_n \subset \div(f) + \Delta_n.
\]
Thus  $\Omega \cap {\Gamma} = W_n \cap {\Gamma} = \emptyset$.    
\end{proof}

\begin{lemma}\label{lem-definemap}
{\em (Compare \cite[Lemma~3.2]{AS}.)}
Let $X$ be a variety, and let $G_1$, $G_2$, and $G_3$ be effective (Cartier) divisors on $X$; let $E := G_3 - G_1 - G_2$.   
  For $i = 1 \ldots 3$, let $U_i \subseteq \abs{G_i}$ be a vector space of dimension at least 2, and suppose that $U_1 U_2 \subseteq U_3$.   Let $\phi_i: X \to \PP^{N_i}$ be the rational map defined by the sections $U_i$ of $G_i$  and let $Y_i$ be the closure of $\im \phi_i$ in $\PP^{N_i}$.  Further assume that $\phi_3: X \to Y_3$ is birational.     Then there is an induced rational map $\pi: Y_3 \to Y_1$ so that $\pi \phi_3 = \phi_1$ and so that if $x \not \in \Bs^{G_i}(U_i)$ for $i =1 \ldots 3$  and $x \not \in  E$, then $\pi$ is defined at $\phi_3(x)$.  
 \end{lemma}
\begin{proof}
We repeat the proof of \cite[Lemma~3.2]{AS}, to note that it works in our situation as well.  
As rational maps, $\pi = \phi_1(\phi_3)^{-1}$.  
Let $x \in X \smallsetminus  ( E \cup \Bs^{G_1}(U_1) \cup \Bs^{G_2}(U_2) \cup \Bs^{G_3}(U_3))$; then all the maps $\phi_i$ are defined at $x$. 
We may thus choose elements $u_0 \in U_1$ and $v \in U_2$ so that, locally at $x$, $G_1 = - \div_X(u_0)$ and $G_2 = -\div_X(v)$.  Our assumptions imply that, locally at $x$, $G_3 = -\div_X(u_0v)$.  Let $\{ u_0, \ldots, u_r\}$ be a basis for $U_1$.   Locally at $x$, $\phi_1$ is defined by $[ u_0: \cdots : u_r]$; we may also define it by $[ u_0v: \cdots :u_r v]$.  Then if
$\{ u_0v, \ldots, u_r v, w_{r+1}, \ldots, w_s\}$ is a basis for $U_3$, then the rational map $\pi$ is given by projection onto the first $r+1$ coordinates.  This is defined locally at $\phi_3(x)$ by construction.  
\end{proof}

We have seen that the maps $X \to \PP^{N_n}$ defined by $\abs{\Delta_n}$ are birational morphisms onto their image for $n \gg 0$.  We now show that for $n \gg 0$, the image of this map gives a normal ample model for $R$, and that our assumptions on $X$ continue to hold.  

\begin{theorem}\label{thm-amplemodel}
Assume Assumption-Notation~\ref{assnot3a}.    Then there are a normal projective surface $X'$, a birational morphism $\theta: X \to X'$, and an ample invertible sheaf $\Lsh'$ on $X'$ such that for some $k \geq 1$,  $\sigma^k$ is conjugate to  a numerically trivial automorphism $\sigma'$ of $X'$ and  $\theta^* (\Lsh') \cong \Lsh_k$.  In particular, $\Lsh'$ is $\sigma'$-ample.  

Further, the gap divisor of $R\ver{k}$ on $X'$ is Cartier and contains no points or components of finite order.  For all $n \geq 1$, the base locus of $\bbar{R}_{nk}$ on $X'$ contains no points of finite order.
\end{theorem}

\begin{proof}
For all $n$, let $\alpha_n$ be the rational map  from $X$ to some projective space given by $\abs{\Delta_n}$; let $X_n$ be the closure of the image of $X$ under $\alpha_n$.    By Theorem~\ref{thm-BPF}, we may replace $R$ by a Veronese subring to assume that $\abs{\Delta_n}$ is base point free for all $n \geq 1$, so $\alpha_n$ is a birational morphism for all $n \geq 1$.  Assumption-Notation~\ref{assnot3a} continues to hold.

For all $n$, we have $\Delta_n + \sigma^{-n} D = \Delta_{n+1}$ and $\abs{\Delta_n} \cdot \abs{\sigma^{-n} D} \subseteq \abs{\Delta_{n+1}}$.  Using Lemma~\ref{lem-definemap} with $E=0$, for each $n \geq 1$ we obtain a birational morphism  $\pi_n: X_{n+1} \to X_n$ so that the diagram
\[ \xymatrix{
X \ar[r]^{\alpha_{n+1}} \ar[rd]_{\alpha_n}	& X_{n+1} \ar[d]^{\pi_n} \\
	& X_n }\]
commutes.  Likewise, the equation $D + \sigma^{-1} \Delta_n = \Delta_{n+1}$ gives a birational morphism $\rho_n: X_{n+1} \to X_n$ so that
\[ \xymatrix{
X \ar[r]^{\alpha_{n+1}} \ar[rd]_{\alpha_n \circ \sigma}	& X_{n+1} \ar[d]^{\rho_n} \\
	& X_n }\]
commutes.   

Let $\Gamma$ be an irreducible curve on $X$.  Then, as $\sigma$ is numerically trivial, 
\[ \Delta_{n+1} . \Gamma = (n+1) D.\Gamma = \frac{n+1}{n} \Delta_n. \Gamma,\]
so $\Delta_{n+1}. \Gamma = 0$ if and only if $\Delta_n . \Gamma =0$.  By Lemma~\ref{lem-contractedcurves}, $\alpha_{n+1}$ and $\alpha_n = \pi_n \circ \alpha_{n+1}$ contract the same curves; thus $\pi_n: X_{n+1} \to X_n$ does not contract any curves and is a finite morphism.  Likewise, $\rho_n$ is a finite morphism.  
By finiteness of the integral closure,  there is some $k$ such that if $n \geq k$, then both $\pi_n$ and $\rho_n$ are isomorphisms.

Let $\bbar{X} := X_{k}$, and let $\alpha = \alpha_{k}: X \to \bbar{X}$.   Define  $\bbar{\sigma} = (\rho_{k} \pi^{-1}_{k})^k$.   Then $\bbar{\sigma}$ is an automorphism of $\bbar{X}$, and we have that 
\[ \bbar{\sigma} \circ \alpha = \alpha \circ \sigma^k.\]
Clearly $\bbar{\sigma}$ is numerically trivial.

   Let $\pi: X' \to \bbar{X}$ be the normalization of $\bbar{X}$.   Since $X$ is normal by assumption, the morphism $\alpha$ factors through $\pi$ --- that is, there is a birational morphism
 $\theta: X \to X'$ such that the diagram
\[ \xymatrix{
X \ar[rr]^{\theta} \ar[rd]_{\alpha} && X' \ar[ld]^{\pi} \\
 & \bbar{X}} \]
 commutes.    Note that if $\theta$ is finite at $x \in X$, then $\theta$ is a local isomorphism at $x$.  
By the universal property of normalizations, $\bbar{\sigma}$ lifts uniquely to an automorphism $\sigma'$ of $X'$, which is also  numerically trivial. 

 By construction, $\bbar{X}$ carries a very ample line bundle $\bbar{\Lsh}$ such that 
 \[ \Lsh_k \cong \struct_X(\Delta_k) \cong \alpha^* \bbar{\Lsh} \cong \theta^* \pi^* \bbar{\Lsh}.\]  
  Let $\Lsh' = \pi^* \bbar{\Lsh}$.  Then $\Lsh'$ is the pullback of an ample line bundle by a finite map and so is ample by \cite[Exercise~III.5.7(d)]{Ha}.  Further, $\Lsh'$ is $\sigma'$-ample by  \cite[Theorem~1.7]{AV}.  By the projection formula \cite[Exercise II.5.1.(d)]{Ha}, $\theta_*(\Lsh_k) \cong \Lsh'$. 

Let $C$ be the union of the finitely many curves in $X$ that are contracted by $\theta$.  Note that $\theta$ is an isomorphism from the open subset $X \smallsetminus C$ of $X$ onto an open subset of $X'$.  Note also that, by Corollary~\ref{cor-noprob}, $\Omega \subseteq X \smallsetminus C$.  Let $\Omega' := \theta(\Omega)$ be the scheme-theoretic image of $\Omega$.  Thus $\Omega'$ is a Cartier divisor on $X'$.     

Let $D'$ be the Cartier divisor on $X'$ corresponding to the invertible sheaf $\Lsh'$.   The singular locus of $X'$ consists of finitely many points, as $X'$ is normal.  Fix $n \geq 1$.  By restricting the Weil divisor $D^{X'}(\bbar{R}_{nk})$ to the open set where $X'$ is smooth, we obtain that
\[ D^{X'}(\bbar{R}_{nk}) = D' + (\sigma')^{-1} (D') + \cdots + (\sigma')^{-(n-1)}(D') - \Omega'.\]
By Lemma~\ref{lem-gap-unique}, $\Omega'$ is the gap divisor of $R\ver{k}$ on $X'$ associated to $z$. That  $\Omega'$ contains no points  or components of of finite order  follows directly from the corresponding properties for $\Omega$.

Fix $n \geq 1$.  We have seen that $D^{X'}(\bbar{R}_{nk})$ is Cartier.  Let $Z_n$ be the base locus of $\bbar{R}_{nk}$ on $X'$.  
Let $x \in X'$ be a point of finite order under $\sigma'$, and let $\Gamma$ be an irreducible component of $\theta^{-1}(x)$.  If $\Gamma$ is a curve, then by Corollary~\ref{cor-noprob}, there is some $f \in \bbar{R}_{nk}$ so that $\Gamma \cap (\div_X(f) + (\Delta_{nk}-\Omega)) = \emptyset$.   If $\Gamma = \{ p\}$ is a point, then it is of finite order, and so  by assumption $p \not \in  W_n$.  Again, there is an $f \in \bbar{R}_{nk}$ so that $p \not \in \div_X(f) + (\Delta_{nk}- \Omega)$.  In either case, $f$ gives a section of $\Lsh'_n$ that does not vanish at $x$, so $x \not \in Z_n$.    Thus $Z_n$ contains no points of finite order.  
\end{proof}

We comment that in the commutative setting,  $\bbar{X}$ would be normal automatically; see \cite[Theorem~2.1.27, Example~2.1.15]{Laz}.  We do not know if this is true for our construction.

\section{Stabilizing 0-dimensional data}\label{DIM0}

We are now ready to start working with the infinite order 0-dimensional data defining $R$.  In this section, we begin with a normal ample model $(X, \sigma)$ for $R$, and give surface data $\surfD = \surfdata$ so that the bimodule algebras $\sh{R}(X)$ and $\sh{T} (\surfD)$ are equal up to taking Veronese subrings.  The key arguments in this section are  essentially combinatorial.

By Theorem~\ref{thm-amplemodel}, by replacing $R$ by a further Veronese, we may without loss of generality make the following assumptions:

\begin{assnot}\label{assnot4}
We assume that $R$ is a birationally commutative projective surface with function field $K$ and fix $0 \neq z \in R_1$.  Let $\bbar{R}_n := R_n z^{-n}$, and  assume that  $\bbar{R}_1$ generates $K$ as a field.  Let  $(X, \sigma)$ be a normal model for $R$ with $\sigma$ numerically trivial, and let $\sh{R}_n := \bbar{R}_n \cdot \struct_X$. Assume also  that there are an ample and   $\sigma$-ample invertible sheaf $\Lsh$ on $X$, an effective locally principal Weil divisor $\Omega$ on $X$ containing no points or components of finite order under $\sigma$, and 0-dimensional subschemes $W_n$ of $X$, disjoint from finite $\sigma$-orbits, such that for all $n \geq 1$, 
$\sh{R}_n = \sh{I}_{W_n} \sh{I}_{\Omega} \Lsh_n$.  
\end{assnot}

We first show that  $\Omega$ can meet $\sigma$-orbits only finitely often.  

\begin{proposition}\label{prop-omega-finite}
Assume Assumption-Notation~\ref{assnot4}.  
Let $p \in X$ be a point of infinite order under $\sigma$; let $O(p)$ denote the $\sigma$-orbit of $p$.  Then $\Omega$ intersects $O(p)$ only finitely often.  
\end{proposition}
\begin{proof}
Suppose  that $O(p) \cap \Omega$ is infinite.  We will show that $R$ is not noetherian.

First suppose that for infinitely many  $d \leq 0$, we have $\sigma^{d}(p) \in \Omega$.     By  Lemma~\ref{lem-finitedata} there is a finite set $V$ 
  such that, for all $n \geq 1$, we have $W_n \subseteq V  \cup \cdots \cup \sigma^{-(n-1)} V$.  
We define a point $q$ as follows:  if $O(p) \cap V = \emptyset$, let $q := p$.  If $O(p)$ meets $V$, let $c := \min \{ d \st \sigma^d(p) \in V \}$ and let $q := \sigma^{c-1}(p)$.  In either case,  for all $n \geq 1$ and $1 \leq m \leq n$, we have $\sigma^{-n}(q) \not\in W_m$.  

Define a left ideal $J$ of $R$ by letting $J := \bigoplus \bbar{J}_n z^n$, where
\[\bbar{J}_n := H^0(X, \Lsh_n \cdot \sh{I}_q^{\sigma^n}) \cap \bbar{R}_n.\]
   If $\sigma^{-n}(q) \in \Omega$, then $\sh{R}_n \subseteq \Lsh_n \sh{I}_q^{\sigma^n}$ and so $J_n = R_n$.  On the other hand, since $\sigma^{-n}(q) \not\in W_n = \Bs(\bbar{R}_n)$ by construction, if $\sigma^{-n}(q) \not\in \Omega$ then there is some section of $\Lsh_n$ in $\bbar{R}_n$ that does not vanish at $\sigma^{-n}(q)$.   That is,  $J_n = R_n$ if and only if $\sigma^{-n}(q) \in \Omega$.

For all $i< n$ we have 
$R_{n-i} J_i \subseteq H^0(X, \sh{I}_\Omega \cdot \sh{I}_q^{\sigma^n} \Lsh_n) z^n$.  
Fix $m \geq 1$ and $n > m$ such that $\sigma^{-n}(q) \in \Omega$.  Then 
\[(R \cdot J_{\leq m})_n \subseteq H^0(X, \sh{I}_\Omega \cdot \sh{I}_q^{\sigma^n} \Lsh_n) z^n.\]
  As $\sigma^{-n}(q) \not \in W_n = \Bs(\bbar{R}_n)$, we have that
\[ \bigl( H^0(X, \sh{I}_\Omega \cdot \sh{I}_q^{\sigma^n} \Lsh_n) \cap  \bbar{R}_n  \bigr) \neq  \bbar{J}_n = \bbar{R}_n.\]
  Thus $J$ is not finitely generated.

Now suppose that for infinitely many  $d \geq 0$, we have $\sigma^d(p) \in \Omega$.  
Let 
\[ \Lsh' :=  \Lsh(-\Omega + \sigma^{-1} \Omega).\]
Then $\sh{R}_n = \sh{I}_{W_n} (\sh{I}_\Omega)^{\sigma^n} (\Lsh')_n$.  That is, $R$ is also contained in a {\em left} idealizer at $\Omega$ inside $B' = B(X, \Lsh', \sigma)$.
An argument symmetric to the one above constructs a right ideal of $R$ that is not finitely generated. 
\end{proof}

We now analyze the 0-dimensional schemes $(\Omega \cup W_n) \cap O(p)$.  To simplify our computations, we will pass to a Veronese subring so that our data may be presented in a standard form.   That we may do so is the content of the following elementary lemma.  

 For any $k \geq 1$, and for any $p \in X$, we will let $O_k(p)$ denote the $\sigma^k$-orbit of $p$.

\begin{lemma}\label{lem-2pts}
Assume Assumption-Notation~\ref{assnot4}.   Then there is a positive integer $k$ such that, for any $p \in X$,  either $O_k(p)$ is disjoint from all $W_n$ or 
there is a point $q \in O_k(p)$ so that $O_k(p) \cap \Omega \subseteq \{q\}$ and 
\[ \{q \} \subseteq  \bigl( O_k(p) \cap (\Omega \cup W_k) \bigr) \subseteq \{ q, \sigma^{-k}( q) \}.\]
\end{lemma}

We first note:  

 \begin{sublemma}\label{sub-easy}
 Suppose that $q$ is a point of infinite order and that 
 \[(\Omega \cup W_1) \cap O(q) \subseteq \{ q, \ldots, \sigma^{-s}(q)\}.\]
 Then 
 \[ (\Omega \cup W_n) \cap O(q) \subseteq \{ q, \sigma^{-1} (q), \ldots, \sigma^{-(n+s-1)} (q) \}\]
 for all $n \geq 1$.
 \end{sublemma}
 \begin{proof}
 This follows from \eqref{kiri}.
 \end{proof}

\begin{proof}[Proof of Lemma~\ref{lem-2pts}]
By Lemma~\ref{lem-finitedata} we know that 
\[\bigcup_{n \geq 1} W_n\]
 is contained in finitely many infinite $\sigma$-orbits.  By Proposition~\ref{prop-omega-finite} each of those orbits meets $\Omega$ only finitely often.  Thus there is some $s \geq 1$ such that for any $p \in \bigcup_{n \geq 1} W_n$, we have that
\[ (\Omega \cup W_1) \cap O(p)  \subseteq \{ \sigma^{-i} (p), \sigma^{-(i+1)}(p), \ldots \sigma^{-(i+s)} (p) \}\]
for some $i \in \ZZ$.

Let $p$ be a point of $\bigcup_{n \geq 1} W_n$. Let 
\[m := \max \{ n \in \ZZ \st \sigma^n(p) \in \Omega \cup W_1\},\]
 and let $q := \sigma^m(p)$.    Then the hypotheses of Sublemma~\ref{sub-easy} hold, and therefore 
 \[ (\Omega \cup W_n) \cap O(q) \subseteq \{ q, \ldots, \sigma^{-(n+s-1)}(q)\} \]
 for all $n \geq 1$.   
 Thus,
  for any $n \geq s$ and any $0 \leq i \leq n-1$, we have that 
\[  (\Omega \cup W_n) \cap O_n(\sigma^{-i} (q)) \subseteq \{\sigma^{-i} (q), \sigma^{-(i+n)} (q)\}.\]

For $0 \leq i \leq  s-1$, 
we have 
\begin{multline*}
 \Omega \cap O_{2s}(\sigma^{-i}(q)) \subseteq 
\bigl( (\Omega \cup W_{2s}) \cap O_{2s}(\sigma^{-i}(q)) \bigr) \cap
\bigl( (\Omega \cup W_s ) \cap O_s(\sigma^{-i}(q)) \bigr) \\
\subseteq 
\{ \sigma^{-i}(q), \sigma^{-(i+2s)}(q) \} \cap \{\sigma^{-i}(q), \sigma^{-(i+s)}(q)\} = \{ \sigma^{-i}(q)\}.
\end{multline*}
For $s \leq i \leq  2s-1$, we have
\begin{multline*}
 \Omega \cap O_{2s}(\sigma^{-i}(q)) \subseteq 
\bigl( (\Omega \cup W_{2s}) \cap O_{2s}(\sigma^{-i}(q)) \bigr) \cap
\bigl( (\Omega \cup W_s ) \cap O_s(\sigma^{-(i-s)}(q)) \bigr) \\
\subseteq 
\{ \sigma^{-i}(q), \sigma^{-(i+2s)}(q) \} \cap \{\sigma^{-(i-s)}(q), \sigma^{-i}(q)\} = \{ \sigma^{-i}(q)\}.
\end{multline*}
Thus the lemma holds for $k = 2s$.
\end{proof}

Lemma~\ref{lem-2pts} allows us to replace $R$ by a Veronese subring so that without loss of generality we may make the following assumptions:

\begin{assnot}\label{assnot5}
We assume that $R$ is a birationally commutative projective surface with function field $K$ and fix $0 \neq z \in R_1$.  Let $\bbar{R}_n := R_n z^{-n}$, and  assume that  $\bbar{R}_1$ generates $K$ as a field.  Let  $(X, \sigma)$ be a normal model for $R$ with $\sigma$ numerically trivial, and let $\sh{R}_n(X) := \bbar{R}_n \cdot \struct_X$. Assume also  that there are an ample and   $\sigma$-ample invertible sheaf $\Lsh$ on $X$, an effective locally principal Weil divisor $\Omega$ on $X$ containing no points or components of finite order under $\sigma$, and 0-dimensional subschemes $W_n$ of $X$, disjoint from finite $\sigma$-orbits, such that for all $n \geq 1$, 
$\sh{R}_n(X) = \sh{I}_{W_n} \sh{I}_{\Omega} \Lsh_n$.  

 In addition, we assume that for any orbit $O(p)$ that meets $\bigcup_{n \geq 1} W_n$, there is some $q \in O(p)$ such that 
\[ \{ q\} \subseteq O(p) \cap (W_1 \cup \Omega) \subseteq \{ q, \sigma^{-1}(q)\},\]
and  $O(p) \cap \Omega \subseteq \{ q\}$.
\end{assnot}

We note here a combinatorial result on the behavior of certain posets, which we state in generality because we will use it here and in Section~\ref{FINDY}.  

\begin{lemma}\label{lem-M}
Let $(\mb{M}, \preceq)$ be a poset satisfying the ascending chain condition, and let
\[ \{ M^n_j \st  1 \leq j \leq n-1 \}\]
be a collection of elements of $\mb{M}$ so that for all $n, j$ we have:

$(1)$  $M^{n+1}_j \succeq M^n_j$, and

$(2)$ $M^{n+1}_{j+1} \succeq M^n_j$.

Then there are some $K \geq 2a \in \NN$ so that if $n \geq K$, then:
\[ M^n_i = M^K_i \mbox{ for } 1 \leq i \leq a.\]
\[ M^n_i = M^K_a \mbox{ for } a \leq i \leq n-a.\]
\[ M^n_i = M^K_{i-n+K} \mbox{ for } n-a \leq i \leq n-1.\]
Further, we have
\[ M^K_1 \preceq M^K_2 \preceq \cdots \preceq M^K_a = \cdots = M^K_{K-a}
 \succeq M^K_{K-a+1} \succeq \cdots \succeq M^K_{K-1}.\]
\end{lemma}
\begin{proof}
For each $n\geq i$, the chain
\[ M^n_i \preceq M_i^{n+1} \preceq \cdots \]
stabilizes; define $M^{\infty}_i$ to be this limiting value.  Likewise, the chains
\[ M^{n}_{n-j} \preceq M^{n+1}_{n+1-j} \preceq \cdots \]
stabilize to a limiting value $N^{\infty}_j$.   Taking the limit of $(1)$ and $(2)$ as $j \to \infty$, we obtain that
\[ M^{\infty}_i \preceq M^{\infty}_{i+1}\]
and
\[ N^{\infty}_j \preceq N^{\infty}_{j+1}\]
for all $i, j \geq 1$.  Define 
\[ M^{\infty}_{\infty} := \lim_{i \to \infty} M^{\infty}_i\]
and
\[ N^{\infty}_{\infty} := \lim_{j \to \infty} N^{\infty}_j.\]

Since the sets
$ \{ M^{\infty}_j \}$ and $ \{N^{\infty}_j\}$ are cofinal subsets of $\{M^n_j\}$, we have equality of the limits
\[ M^{\infty}_\infty = N^{\infty}_{\infty}.\]
Thus we may choose $a \geq 1$ and $r \geq 2a$ so that
\[ M^{\infty}_{\infty} = M^{\infty}_a = M^r_a = M^r_{r-a} = N^{\infty}_a.\]
Possibly increasing $r$, we may assume that
$ M^r_i = M^{\infty}_i$ and $N^r_{r-i} = N^{\infty}_i$ for all $1 \leq i \leq a$.

Note that if $e \geq c$ and $b \leq d \leq b+e-c$, then 
\[ M^c_b \preceq M^e_{d}.\]
Let $K := 2r-2a$.  Let $n \geq K$.  If  $ a \leq i \leq n/2 \leq n+a-r$, then 
\[ M^{\infty}_{\infty} = M^r_a \preceq M^n_i \preceq M^{\infty}_{\infty}.\]
In particular, $M^K_a = M^{\infty}_{\infty}$.  
Likewise, if $n-a \geq i \geq n/2 \geq r-a$, then
\[ M^K_a = M^{\infty}_{\infty} = M^r_{r-a} \preceq M^n_i \preceq M^{\infty}_{\infty}.\]
By choice of $r$,  if $1 \leq i \leq a$ and $n \geq K$, then $M^n_i = M^{\infty}_i = M^r_i$, and
$M^n_{n-i} =  N^{\infty}_i= M^r_{r-i}$.  
\end{proof}

\begin{lemma}\label{lem-stabledata}
Assume Assumption-Notation~\ref{assnot5}.  
Let $p \in  \bigcup_{n \geq 1} W_n$, and let $\struct := \struct_{X,p}$ be the local ring of $X$ at $p$, with maximal ideal $\mf{p}$.  For all $j \geq 1$ and all $ i \in \zed$, define
$\mf{m}^j_i$ to be the stalk of the ideal sheaf $\sh{R}_j \Lsh_j^{-1} = \sh{I}_\Omega \sh{I}_{W_n}$ at $\sigma^{-i}(p)$, considered as an ideal in $\struct$ via the isomorphism $\sigma^i:  \struct_{X, \sigma^{-i} p} \to \struct$.   Our assumptions imply that by reindexing the orbit of $p$ if necessary, we may assume that $\mf{m}^1_i = \struct$ for all $i < 0$ and $i > 1$, that $\mf{m}^1_0 \neq \struct$, and that $\Omega \cap O(p) \subseteq \{ p\}$.

Then there are integers $ N \geq t  \geq 1$,  ideals $\mf{a}_1, \ldots \mf{a}_{t-1}, \mf{d}, \mf{c}_{t-1}, \ldots \mf{c}_0$ of $\struct$ that are either $\mf{p}$-primary or equal to $\struct$, and an ideal $\mf{a}_0$ of $\struct$ so that for all $n \geq N$, we have:
\begin{align*} \mf{m}^n_i  = \mf{a}_i & \mbox{ for $0 \leq i < t$} \\
\mf{m}^n_i = \mf{d} & \mbox{ if $t \leq  i \leq n-t$} \\
\mf{m}^n_i = \mf{c}_{n-i} & \mbox{ if $n - t < i \leq n$} \\
\mf{m}^n_i = \struct & \mbox{ if $i < 0$ or $i > n$}.
\end{align*}
Further, we have $\mf{a}_0 \mf{c}_0 \subseteq \mf{d}$ and 
\[ \mf{a}_1 \subseteq \mf{a}_2 \subseteq \cdots \subseteq \mf{a}_{t-1} \subseteq \mf{d} \supseteq \mf{c}_{t-1} \supseteq \cdots \supseteq \mf{c}_{1}.\]
\end{lemma}
\begin{proof}
The multiplication $\sh{R}_n \sh{R}_m^{\sigma^n} \subseteq \sh{R}_{n+m}$ implies that
\[ \mf{m}^n_i \mf{m}^m_{i-n} \subseteq \mf{m}^{n+m}_i\]
for all $i, n, m$.  In particular, since if $k < 0 $ and $\ell > 1$ then $\mf{m}^1_k = \mf{m}^1_{\ell} = \struct$, we have that
\[ \mf{m}^n_i = \mf{m}^n_i \mf{m}^1_{i-n} \subseteq \mf{m}^{n+1}_i\]
if $i \leq n-1$, and
\[ \mf{m}^n_i = \mf{m}^1_{i+1} \mf{m}^n_i \subseteq \mf{m}^{n+1}_{i+1}\]
if $i \geq 1$.  Thus the poset
\[ \{ \mf{m}^n_i \st 1 \leq i \leq n-1 \},\]
of ideals of $\struct$, partially ordered by inclusion, satisfies the hypotheses of Lemma~\ref{lem-M}. The existence of
$N, t$, and 
\[ \mf{a}_1 \subseteq \mf{a}_2 \subseteq \cdots \subseteq \mf{a}_{t-1}   \subseteq \mf{a}_t = \mf{d} = \mf{c}_t  \supseteq \mf{c}_{t-1} \supseteq \cdots \supseteq \mf{c}_{1},\]
as claimed, follows directly from Lemma~\ref{lem-M}.   Since the chains $\{ \mf{m}^n_0\}$ and $\{ \mf{m}^n_n\}$ have maximal elements, by possibly increasing $N$ we may also obtain $\mf{a}_0 = \mf{m}^N_0$ and $\mf{c}_0 = \mf{m}^N_N$ as described.  Then
\[ \mf{a}_0 \mf{c}_0 = \mf{m}^N_N \mf{m}^N_0 \subseteq \mf{m}^{2N}_N = \mf{d},\]
as claimed.
\end{proof}

The ideal $\mf{d}$ constructed in Lemma~\ref{lem-stabledata} is called the {\em central stalk} of $R$ at the orbit $O(p)$.  
\begin{corollary}\label{cor-stabledata}
Let $t$ be the integer constructed in Lemma~\ref{lem-stabledata}.  
For all $s \geq t$ there is an integer $M \geq 1$ so that for all $n \geq M$ we have:
\begin{enumerate}
\item If $ 0 \leq k\leq s$ then
$\mf{m}^{2ns}_{2is+k} = \mf{d}$ if $1 \leq i \leq n-1$, and $\mf{m}^{2ns}_{2is+k} = \struct$ if $i <0 $ or $i > n$.
\item  If $s+1 \leq k \leq 2s-1$ then
$\mf{m}^{2ns}_{2is+k} = \mf{d}$ if $0 \leq i \leq n-2$, and $\mf{m}^{2ns}_{2is+k} = \struct$ if $i< 0 $ or $i \geq n$.
\end{enumerate}
\end{corollary}
\begin{proof}
This follows immediately from Lemma~\ref{lem-stabledata}.
\end{proof}

We may now give the defining data for  the bimodule algebra $\sh{R} = \sh{R}(X)$.  

\begin{defn}\label{def-admissible}
We will say that the surface data $\surfD = \surfdata$ is {\em  \adm} if:
\begin{itemize}
\item $X$ is normal;
\item $\sigma$ is numerically trivial;
\item $\Lsh$ is ample (and therefore $\sigma$-ample);  
\item $\Omega$ contains no points or components of finite order under $\sigma$; 
\item the cosupport of $\sh{ADC}$ meets $\sigma$-orbits at most once; and
\item if $w \in \Cosupp \sh{ADC}$, then $\Omega \cap O(w) \subseteq \{w\}$.
\end{itemize}
\end{defn}

  Note that  if $\surfD$ is \adm, then $\Omega$ is locally principal.

\begin{theorem}\label{thm-data}
Let $R$ be   a birationally commutative projective surface.  Then there are  \adm\ surface data $\surfD = \surfdata$  and an integer $ k \geq 1$ so that 
\[\sh{R} (X) \ver{k} = 
\sh{T} (\surfD). \]
 \end{theorem}
 \begin{proof}
By Proposition~\ref{prop-blowup}, after replacing $R$ by a Veronese subring we may assume that we are in the situation of Assumption-Notation~\ref{assnot3a}.  By Theorem~\ref{thm-amplemodel}, by replacing $R$ by a further Veronese subring and possibly changing $X$, we may assume that $R$ and $X$ satisfy Assumption-Notation~\ref{assnot4}.  
%
  
 By Lemma~\ref{lem-2pts}, we may replace $R$ by a further Veronese subring to assume that for all $p$ such that $O(p)$ meets $W_1$, there is a $q \in O(p)$ so that 
 \[ \{q\} \subseteq O(p) \cap (\Omega \cup W_1) \subseteq \{q, \sigma^{-1}(q)\}\]
 and $O(p) \cap \Omega \subseteq \{q\}$.  By Sublemma~\ref{sub-easy}, 
 \[ O(p) \cap (\Omega \cup W_n) \subseteq \{ q, \ldots, \sigma^{-n} (q)\}\]
 for all $n \geq 1$.  By Lemma~\ref{lem-finitedata}, there are  points
 $q^1, \ldots, q^r$ with orbits $O^j := O(q^j)$, so that 
 \[ \bigcup_{n \geq 1} W_n \subseteq \bigcup_{j =1}^r O^j.\]

For $j = 1 \ldots r$, apply Lemma~\ref{lem-stabledata} to $O^j$ to obtain positive integers  $t^j$ and $N^j$ and  ideals   $\mf{d}^j$, $\mf{a}^j_i$, and $\mf{c}^j_i$ in $\struct_{X, q^j}$, for $0 \leq i \leq t^j-1$.  Let $s := \max\{t^j, N^j\}$.  For integers $i$ with  $t^j \leq i \leq s-1$, define $\mf{a}^j_i := \mf{c}^j_i := \mf{d}^j$.  We have
\[ \mf{a}^j_1 \subseteq \cdots \subseteq \mf{a}^j_{s-1} \subseteq \mf{d} \supseteq \mf{c}^j_{s-1} \supseteq \cdots \supseteq \mf{c}^j_1\]
and $\mf{a}^j_0 \mf{c}^j_0 \subseteq \mf{d}^j$ 
for all $j$.  Corollary~\ref{cor-stabledata} implies that, replacing $R$ by $R\ver{2s}$ (and so replacing $\sigma$ by $\sigma^{2s}$, replacing $\Lsh$ by $\Lsh_{2s}$, and replacing $r$ by $2sr$) and reindexing orbits as needed, we may presume that the conclusions of Lemma~\ref{lem-stabledata} hold on all orbits, with $N = t = 1$.

We define an ideal sheaf  $\sh{A}' \subseteq \sh{I}_{\Omega}$ so that $\sh{I}_{\Omega}/\sh{A}'$ is supported on 
$\{ q^1, \ldots, q^r \}$ 
 by setting the stalk of $\sh{A}'$ at $q^j$ to be isomorphic to $\mf{a}^j_0$.   Similarly, we define an ideal sheaf $\sh{C}$, cosupported on 
 $\{q^1, \ldots, q^r \}$,
 by setting the stalk of $\sh{C}$ at $q^j$ to be isomorphic to $\mf{c}^j_0$.  
Let $\sh{D}$ be the ideal sheaf cosupported at $\{q^1, \ldots, q^r\}$ so that 
\[ \sh{D}_{q^j} \cong \mf{d}^j.\]
Then $\sh{A}'\sh{C} \subseteq \sh{D}$.  Let $\sh{A}:=\sh{I}_{\Omega}^{-1}\sh{A}'$.

Now, $\surfD:=\surfdata$ is surface data, and we have seen it is in addition normal.
By construction,
\[ \sh{R}_n  = \sh{A}' \sh{D}^{\sigma} \cdots \sh{D}^{\sigma^{n-1}} \sh{C}^{\sigma^n} \Lsh_n 
= \sh{I}_{\Omega}  \sh{A} \sh{D}^{\sigma} \cdots \sh{D}^{\sigma^{n-1}} \sh{C}^{\sigma^n} \Lsh_n
\]
if $n \geq 1$.  Thus $\sh{R} = \sh{T}(\surfD)$.    
\end{proof}

We record for future reference  the effect of taking Veronese subrings upon surface data. 
 
 \begin{lemma}\label{lem-Veronese-data}
 Suppose that 
 \[\surfD' = \surfdata\]
  is surface data.  Let $n \geq 1$, and let 
 \[ \surfD := (X, \Lsh_n, \sigma^n, \sh{A} \sh{D}^{\sigma} \cdots \sh{D}^{\sigma^{n-1}}, 
 \sh{D} \sh{D}^{\sigma} \cdots \sh{D}^{\sigma^{n-1}}, \sh{C}, \Omega).\]
 Then $\surfD$ is surface data, and 
 \[ \sh{T}(\surfD') \ver{n} = \sh{T}(\surfD).\]
 Furthermore, if the surface data  $\surfD'$ is \adm, respectively transverse, then  the surface data $\surfD$ is \adm, respectively transverse.
 \end{lemma}
 \begin{proof}
 This is an elementary computation, which we leave to the reader.
 \end{proof}

\section{Transversality of the defining data}\label{DATA-CT}
In Section~\ref{DIM0}, we constructed  \adm\ surface data $ \surfD = \surfdata$ such that (up to a Veronese) we have  $\sh{R}(X) = \sh{T}(\surfD)$.  In this section, we show that  the data $\surfD$ is in fact transverse, and that $T(\surfD)$ is a finite module over (a Veronese of) $R$.   This allows us to consider $T(\surfD)$ a kind  of normalization of $R$, and further justifies the term ``normal surface data.''

By Theorem~\ref{thm-data}, we may assume without loss of generality that we are in the following situation:
\begin{assnot}\label{assnot6}
We assume that $R$ is a birationally commutative projective surface with $R_1 \neq 0$ and fix $0 \neq z \in R_1$.   As usual, we define $\bbar{R}_n:= R_n z^{-n}$.   In addition, we assume that $\bbar{R}_1$ generates $K$ as a field.  We assume also that there is   surface data $\surfD = \surfdata$, \adm\ in the sense of Definition~\ref{def-admissible}, so that  if $\sh{R}_n(X) := \bbar{R}_n \cdot \struct_X$, then 
\[ \sh{R}(X) = \sh{T} (\surfD).\]
  We will continue to let $W_n$ be the base locus of $\bbar{R}_n$ for $n \geq 1$, so that $W_n$ is defined by 
  \[  \sh{A} \sh{D}^{\sigma } \cdots \sh{D}^{\sigma^{n-1}} \sh{C}^{\sigma^n}\]
  for all $n \geq 1$.

Let $Z$  be the subscheme defined by $\sh{D}$, let $\Lambda$ be the subscheme defined by $\sh{A}$, and let $\Lambda'$ be the subscheme defined by $\sh{C}$.
\end{assnot}


We first prove the unsurprising result that in this situation $\Omega$ has good transversality properties.  
\begin{lemma}\label{lem-omega-transverse}
The set $\{\sigma^n \Omega\}_{n \in \zed}$ is critically transverse.
\end{lemma}
\begin{proof}
By assumption, $\Omega$ contains no points or components of finite order under $\sigma$.  Since $X$ is normal, $\Omega$ is locally principal.  By Proposition~\ref{prop-omega-finite}, $\Omega$ meets infinite orbits at only finitely many points.  Thus by \cite[Lemma~3.1]{S-surfprop}, $\{ \sigma^n \Omega \}$ is critically transverse.
\end{proof}

We next prove two lemmas that will, in many cases, allow us to work with the full algebra $T(\surfD)$ instead of the subalgebra $R$.  The first is an easy generalization of a lemma of Rogalski and Stafford.

\begin{lemma}\label{lem-surjective}
{\em(Compare \cite[Lemma~9.3]{RS})}
Let $X$ be a projective scheme with automorphism $\sigma$.  Let $\{(\sh{R}_n)_{\sigma^n}\}$ be a left and right ample of sequence of bimodules on $X$ such that for each $n$, the set where  $\sh{R}_n$ is not locally free has dimension $\leq 0$.  Let $\sh{F}$ be a globally generated coherent sheaf on $X$ and let $V \subseteq H^0(X, \sh{F})$ be a vector space that generates $\sh{F}$.   Let $i \in \ZZ$.  Then for $n \gg 0$, the natural homomorphism
\[ \alpha: V \otimes H^0(X, \sh{R}_n^{\sigma^i}) \to H^0(X, \sh{F} \otimes \sh{R}_n^{\sigma^i})\]
is surjective.
\end{lemma}
\begin{proof}
By assumption, there is an exact sequence
\[ \xymatrix{
0 \ar[r] 
&\sh{H} \ar[r] 
&V \otimes \struct_X \ar[r] 
& \sh{F} \ar[r] 
&0.}\]
Tensoring with $\sh{R}_n^{\sigma^i}$, we obtain an exact sequence
\[ \xymatrix{
0 \ar[r]
& \shTor^X_1(\sh{F}, \sh{R}_n^{\sigma^i}) \ar[r]
& \sh{H} \otimes \sh{R}_n^{\sigma^i} \ar[r]^{\theta_n}
& V \otimes \sh{R}_n^{\sigma^i} \ar[r]
& \sh{F} \otimes \sh{R}_n^{\sigma^i} \ar[r]
& 0.}\]

Let $\sh{K}_n := \im \theta_n$.  Our assumptions on $\sh{R}$ imply that $\shTor^X_1(\sh{F}, \sh{R}_n^{\sigma^i})$ is supported on a set of dimension 0, and so $H^i(X,  \shTor^X_1(\sh{F}, \sh{R}_n^{\sigma^i}))= 0$ for all $n$ and for all $i \geq 1$.  Thus $H^1(X, \sh{K}_n) \cong H^1(X, \sh{H}\otimes \sh{R}_n^{\sigma^i})$.  By ampleness of $\{ (\sh{R}_n)_{\sigma^n}\}$, this vanishes for $n \gg 0$.  From the long exact cohomology sequence
\[ \xymatrix{
  H^0(X, \sh{K}_n) \to V \otimes_k H^0(X, \sh{R}_n^{\sigma^i}) \ar[r]^>>>>>{\alpha}
 & H^0 (X, \sh{F} \otimes \sh{R}_n^{\sigma^i}) \ar[r]
 & H^1(X, \sh{K}_n), }\]
 we deduce that $\alpha$ is surjective for $n \gg 0$.
 \end{proof}

 \begin{lemma}\label{lem-fg}
  Let $X$ be a projective scheme, let $\sigma$ be an automorphism of $X$, and let $\Lsh$ be a $\sigma$-ample invertible sheaf on $X$.  Suppose that $R$ is a finitely generated graded subalgebra of $B(X, \Lsh, \sigma) = \bigoplus H^0(X, \Lsh_n) z^n$.  For all $n \geq 1$ let $\sh{R}_n \subseteq \Lsh_n$ be the sheaf generated by the sections in $\bbar{R}_n$.    Let $T := \bigoplus_{n \geq 0} H^0(X, \sh{R}_n) z^n$.  
  
  Suppose that for all $n$, the set where $\sh{R}_n$ is not locally free has dimension $\leq 0$ and that the sequence of bimodules $\{ (\sh{R}_n)_{\sigma^n}\}$ is left and right ample.  
Then $T$ is finitely generated as a left and right $R$-module. 
  \end{lemma}
  \begin{proof}
  By symmetry, it suffices to prove that $_{R} T$ is finitely generated.  
  
  Let $k$ be such that 
  \[ R_n = \sum_{i =1}^k R_i R_{n-i}\]
  for all $n > k$.   Then 
  \[ \sh{R}_n = \sum_{i=1}^k \sh{R}_i \sh{R}_{n-i}^{\sigma^i}\]
  for all $n > k$; taking global sections we have
  \beq\label{widget}
  T_n = \sum_{i=1}^k H^0(X, \sh{R}_i \sh{R}_{n-i}^{\sigma^i}) z^n
  \eeq
  for all $n >k$.
  
    For each $1 \leq i \leq k$, the sections in $\bbar{R}_i$ generate $\sh{R}_i$.  Applying Lemma~\ref{lem-surjective}, we obtain that there is some $n_0$, which we may take to be greater than $k$, so that the multiplication map
  \[ R_i \otimes T_{n-i} \to H^0(X, \sh{R}_i \otimes \sh{R}_{n-i}^{\sigma^i})z^n \]
  is surjective for $n \geq n_0$ and $1 \leq i \leq k$.     

Now consider the exact sequence
\[ 0 \to  \sh{J}_{i,n} \to \sh{R}_i \otimes \sh{R}_{n-i}^{\sigma^i} \to \sh{R}_i \sh{R}_{n-i}^{\sigma^i} \to 0. \]
The kernel $\sh{J}_{i,n}$ is supported at finitely many points, and so $H^1(X, \sh{J}_{i,n}) = 0$.  Thus the induced map from $H^0(X, \sh{R}_i \otimes \sh{R}_{n-i}^{\sigma^i}) \to H^0(X, \sh{R}_i \sh{R}_{n-i}^{\sigma^i})$ is surjective.  Therefore, for all $n \geq n_0$, the natural map
\[ R_i \otimes T_{n-i} \to H^0(X, \sh{R}_i \sh{R}_{n-i}^{\sigma^i})z^n\]
is surjective.  Applying \eqref{widget}, we see that for $n \geq n_0$, 
\[ T_n = \sum_{i=1}^k R_i T_{n-i}.\]
By induction,  $T$ is generated as a left $R$-module by $T_{\leq n_0}$.
\end{proof}

The next step in proving transversality of the data $\surfD$ is to show that $\sh{ADC}$ is cosupported on points with dense orbits.  This is the key result of the current section, and is surprisingly delicate to prove.  

\begin{proposition}\label{newprop}
Assume Assumption-Notation~\ref{assnot6}.  Then all points in $Z \cup \Lambda \cup \Lambda' $ have dense $\sigma$-orbits.
\end{proposition}

\begin{proof}
Suppose that there is a point in $Z \cup \Lambda \cup \Lambda'$ without a dense orbit.  We claim that $R$ is not noetherian.

Let $C$ be the Zariski closure of the orbits of  all points without dense orbits in $\Supp(\Lambda \cup Z \cup \Lambda')$.  Then $C$ is a reduced but not necessarily irreducible curve on $X$.   By passing to a Veronese, we may assume that all irreducible components of $C$ are $\sigma$-invariant.  

Let $\sh{A}' \supseteq \sh{I}_{\Omega} \sh{A}$ and $\sh{C}' \supseteq \sh{C}$ be maximal so that
\[ \sh{A}'\sh{C}' \subseteq \sh{D}.\]
 For all $n \geq 1$, let $\Pi_n$ be the closed subscheme of $X$ defined by 
$\sh{A}' \sh{D}^{\sigma} \cdots \sh{D}^{\sigma^{n-1}} (\sh{C}')^{\sigma^n}$.   For any irreducible component $\Gamma$ of $C$,  the intersections $W_n \cap \Gamma$  are supported on nonsingular points of $X$ of infinite order, which are therefore also nonsingular points of $\Gamma$; thus $\deg_{\Gamma}(\Pi_n \cap \Gamma)$ and $\deg_{\Gamma}(Z \cap \Gamma)$ make sense.  Because of the maximality of $\sh{A}'$ and $\sh{C}'$, we have  
\beq\label{joehill}
\deg_{\Gamma}(\Pi_n \cap \Gamma) = n \deg_{\Gamma} (Z \cap \Gamma)
\eeq
for all $n \geq 1$.

Fix $0 \neq f \in \bbar{R}_1$, and let $F = \div_X(f) + \Delta_1$.  As 
\[f \in H^0(X, \sh{R}_1) \subseteq H^0(X, \sh{A}'(\sh{C}')^{\sigma} \Lsh),\]
we have that $F \cap \Gamma \supseteq \Pi_1 \cap \Gamma$.  Thus 
\[ \deg_{\Gamma}(\Lsh|_{\Gamma}) = \Delta_1 . \Gamma = F.\Gamma \geq \deg_{\Gamma}(\Pi_1 \cap \Gamma) = \deg_{\Gamma}(Z \cap \Gamma).\]

We first suppose that this inequality is strict.  Thus suppose there is an irreducible component $\Gamma$ of $C$ so that
\beq\label{case1}
 \deg_{\Gamma}(\Lsh |_{\Gamma}) > \deg_{\Gamma} (Z \cap \Gamma).
 \eeq

 Suppose also that $Z \cap \Gamma \neq \emptyset$.   For any $p \in Z\cap \Gamma $, consider the closed subscheme $Z_p = Z|_{\{p\}}$ of $Z$ supported at $p$.  For any such $p$, there is an integer $d_p > 0$ such that  the Zariski closure of ${\{ \sigma^n(Z_p)\}}_{n \in \ZZ}$ is $d_p \Gamma $.   Let 
\[d = \min \{ d_p \st p \in Z \cap \Gamma\},\]
and let $x \in Z$ be a point with $d_x = d$.

Let ${G}$ be  the $(d+1)$-uple curve defined by the Weil divisor $(d+1) \Gamma$.   There is a natural map 
\[ \phi:  B(X, \Lsh, \sigma) \to B({G}, \Lsh|_{{G}}, \sigma|_G).\]
  Let $S := \phi(R)$.  That is, 
  \[ S = \bigoplus_{n \geq 0} \biggl( \frac{\bbar{R}_n}{ H^0(X, \sh{I}_{{G}} \Lsh_n) \cap \bbar{R}_n} \biggr) z^n = \bigoplus_{ n \geq 0} \bbar{S}_n z^n\]
  We claim that $S$ is not noetherian.  
  
Let $\sh{M}_n := \Lsh_n|_{{G}}$.  
 For all $n$, let $\sh{S}_n$ be the image of $\sh{R}_n \otimes \struct_{{G}}$ under the natural map 
 \[ \sh{R}_n \otimes \struct_{{G}} \to \sh{M}_n.\]
 The sections in $\bbar{S}_n$ generate the subsheaf $\sh{S}_n$ of $\sh{M}_n$.  
  Let 
  \[ k = \deg_{\Gamma}(Z \cap \Gamma) = \deg_{\Gamma} (\Pi_1 \cap {\Gamma}),\]
  and let 
  \[ \ell = \deg_{\Gamma} (\Lsh |_{\Gamma}).\]
  By \eqref{case1}, $\ell > k > 0$.
    
Let $\sh{S}'_n := \sh{R}_n \cdot \struct_{\Gamma}$ be the image of the natural map $\sh{R}_n \otimes \struct_{\Gamma} \to \Lsh_n |_{\Gamma}$.
   One can easily see that the formula for $\sh{T}_n$ gives a constant $c \geq 0$ so that 
      \[ \deg_{\Gamma}(\Omega \cap \Gamma + W_n \cap \Gamma) = nk + c\]
      for all $n \gg 0$.  
Thus 
  \[ \deg_{\Gamma} (\sh{S}'_n) = n(\ell-k)-c\]
  for all $n \gg 0$.  
 
We will work with the nonreduced scheme ${G}$ carefully.   Fix $n$; let $V_n$ be the subscheme of $X$ defined by $\sh{I}_{\Omega} \sh{I}_{W_n}$.  Let 
$P$ be the scheme-theoretic intersection $ {G} \cap V_n$.  Let 
$ \Supp P = \{p_1, \ldots, p_r \}$.  Recall that $\Gamma$ is nonsingular at all $p_i$.  Therefore, 
$ P \cap \Gamma$ is the divisor $m_1 p_1 + \cdots + m_r p_r$ on $\Gamma$ 
for some integers $m_i \geq 1$.  

 If $f \in \struct_{X, p_i}$, let $\bbar{f}$ be its  image in $\struct_{\Gamma, p_i}$.  Then for $i =1, \ldots, r$, there are elements
$f_i \in (\sh{R}_n \Lsh_n^{-1})_{p_i} \subseteq \struct_{X, p_i}$ so that the valuation of $ \bbar{f_i}$ in the discrete valuation ring $\struct_{\Gamma, p_i}$ is $m_i>0$.   In particular, the image of $f_i$ in $\struct_{{G}, p_i}$ is a non-zerodivisor.  
By taking the locally free rank 1 ideal sheaf on ${G}$ generated by the images of the $f_i$ in $\struct_{{G}, p_i}$, we obtain an invertible ideal sheaf $\sh{N}_n$ on ${G}$.  The sheaf $\sh{N}_n$  defines a locally principal subscheme $Q$ of ${G}$, with $Q \subseteq P$, so that the scheme-theoretic intersections $Q \cap \Gamma$ and $ P \cap \Gamma$ are equal.   

 Let $\sh{N}'_n := \sh{N}_n \otimes \sh{M}_n$; then $\sh{N}'_n$ is an invertible subsheaf of $\sh{S}_n$ with 
\[\deg_{\Gamma}(\sh{N}'_n |_{\Gamma}) = \deg_{\Gamma}( \sh{S}'_n).\]
Thus 
\[ \deg_{\Gamma}(\sh{N}'_n|_{\Gamma}) = n(l-k)-c\]
 for $n \gg 0$.

 As 
\[\lim_{n \to \infty} n(l-k)-c = \infty,\]
by \cite[Corollary~3.14]{RS-0}
 the sequence of bimodules $\{ (\sh{N}'_n)_{\sigma^n}\}_{n \geq 0}$ is a left and right ample sequence on ${G}$.  
Since the cokernel of 
 \[   \sh{N}'_n \to  \sh{S}_n\]
 is supported  on  a set of dimension 0,   by \cite[Corollary~3.7]{S-surfprop}, $\{(\sh{S}_n)_{\sigma^n}\}$ is also a left and right ample sequence on ${G}$.  
 
 Let 
 \[ A := \bigoplus_{n \geq 0} H^0({G}, \sh{S}_n) z^n.\]
   By Lemma~\ref{lem-fg}, $A$ is finitely generated as a left and right $S$-module.  Thus it suffices to prove that  $A$ is not noetherian.  
   
   Let $\sh{J}$ be the ($\sigma$-invariant) ideal sheaf of $d \Gamma $ on ${G}$. 
   Let $J$ be the ideal 
    \[ J := \bigoplus_{n \geq 0} H^0({G}, \sh{J} \sh{M}_n \cap \sh{S}_n)z^n\]
 of $T$.    Let $\sh{E}$ be the subsheaf $ \sh{D} \struct_{G}$ of $\struct_{G}$.  
 
 There are integers $a, n_0 \geq 0$ so that 
 \[ (\sh{S}_n)_{\sigma^{-(n-a)}(x)} = (\sh{E}^{\sigma^{n-a}}\sh{M}_n)_{\sigma^{-(n-a)}(x)}\]
 for all $n \geq n_0$;  in fact, we may take $a=1$.  As $d\Gamma$ is the Zariski closure of $\{ \sigma^n(Z_x)\}$, we have 
containments
 \[ (\sh{JM}_n)_{\sigma^{-(n-a)}(x)} \subseteq (\sh{S}_n)_{\sigma^{-(n-a)}(x)} \subseteq (\sh{M}_n)_{\sigma^{-(n-a)}(x)}.\]
 Thus for any $m \geq n_0$ and $n \geq 1$, we have
 \[   J_n A_m\subseteq 
 H^0({G},  \sh{J} \sh{M}_{n}(\sh{E}^{\sigma^{m-a}} \sh{M}_m)^{\sigma^n})z^{m+n} = 
 H^0({G}, \sh{J} \sh{E}^{\sigma^{n+m-a}}   \sh{M}_{n+m}) z^{m+n}.\]

 The kernel and cokernel of 
 \[ \sh{J} \otimes \sh{S}_n \to \sh{J} \sh{M}_n \cap \sh{S}_n\]
 are supported on sets of dimension 0, and $\{ (\sh{S}_n)_{\sigma^n}\}$ is a left and right ample sequence on ${G}$.  Thus by \cite[Corollary~3.7]{S-surfprop}, there is $n_1 $ so that   the sheaf $\sh{J M}_n \cap \sh{S}_n$ is globally generated for $n \geq n_1$.    We may assume that $n_1 \geq n_0$.  
     Then  for any 
$ m \geq n_1$ and $n \geq m+n_0$, we have
\[ ((J_{\leq m}) \cdot A)_n \subseteq H^0({G},   \sh{J} \sh{E}^{\sigma^{n-a}} \sh{M}_{n}) z^{n} \subsetneqq J_{n}.\]
Thus $J$ is not finitely generated as a 
  right ideal of $A$.
 
If \eqref{case1} holds and $Z \cap \Gamma = \emptyset$, then there is some $p \in (\Lambda \cup \Lambda')\cap \Gamma$.  Suppose that $p \in \Lambda'$.  Let $d $ be such that $d \Gamma \supseteq (\Lambda')_p$.
Let $G := (d+1)\Gamma$ and let $\sh{S}_n$ be the image of $\sh{R}_n \otimes \struct_{G} \to \Lsh_n |_{G}$.  Let 
\[ A := \bigoplus_{n \geq 0} H^0(G, \sh{S}_n)z^n.\]
As above, it suffices to prove that $A$ is not noetherian.  If $\sh{J}$ is the ideal sheaf of $d\Gamma$ on $G$, then it follows as above that the right ideal
\[ J := \bigoplus_{n \geq 0} H^0\bigl(G, \sh{J} (\Lsh_n |_{G}) \cap \sh{S}_n\bigr)\]
of $A$ is not finitely generated.
The case that $p \in \Lambda$ is similar.

 It remains to consider the case that  
 \[
 \deg_{C}(Z \cap  C) = \deg_{C} (\Lsh |_{C}).
 \]
Let $\Gamma$ be an irreducible component of $C$ and let $n \geq 1$.  Since
\[
 n \deg_{\Gamma}(\Lsh |_{\Gamma}) \geq \deg_{\Gamma}(\Omega \cap {\Gamma}) + \deg_{\Gamma}(W_n \cap {\Gamma}) \geq \deg_{\Gamma}(\Pi_n \cap {\Gamma}) = n \deg_{\Gamma}(Z \cap {\Gamma}),
\]
  we have 
  \beq\label{case2}
 \deg_{\Gamma}(\Pi_n \cap  \Gamma) = n \deg_{\Gamma} (Z \cap \Gamma) = \deg_{\Gamma} (\Lsh_n |_{\Gamma}) = \deg_{\Gamma}(\Omega \cap {\Gamma}) + \deg_{\Gamma}(W_n \cap {\Gamma}) > 0.
 \eeq
Choose some $p \in Z \cap \Gamma$.  For all $i$, let 
\[ p_i := \sigma^{-i}(p).\]
By reindexing the orbit of $p$ if necessary, we may assume that  $p_i \not\in Z$ for $i <0$.  Our assumption on the defining data for  $\sh{R}  $ implies that  if $i \geq 2$, then $p_i \not \in \Omega \cup W_1$.   Let $\struct := \struct_{X, p}$.   As usual, we will identify all $\struct_{X, p_i}$ with $\struct$.  Note that $\struct$ is a regular local ring of dimension 2, since $X$ is normal by assumption and the orbit of $p$ is infinite.  
Let $\mf{d}$ be the central stalk of $R$ at $O(p)$; that is, there are integers $b$, which we assume to be at least 1, and $N$, which we assume to be at least $2b$, so that for $n \geq N$ and $b \leq i \leq n-b$ we have that
\[ (\sh{R}_n)_{p_i} = \mf{d}.\]

The point $p$ has infinite order on $\Gamma$, and so  $\Gamma $ is  nonsingular at $p$.  Let $y$ be the local equation of $\Gamma$ in $\struct$.  Thus there is some  $x \in \struct$ so that $x$ and $y$ generate the maximal ideal $\mf{m}$ of $\struct$. 

 For all $i \geq 1$, let
\[ J(i) := \bigoplus_{n \geq 1} (\bbar{R}_n \cap H^0(X, \sh{I}_{i \Gamma} \Lsh_n))z^n \subseteq R.\]
Let $d = \min\{ i \st y^i \in \mf{d} \}$.  
Let \[ A: = R/J(d+1).\]  
We will show that $A$ is not noetherian.  
We first show that $\GK  A \geq 2$.  

Now, $\bbar{R}_N \not \subseteq H^0(X, \sh{I}_{\Gamma} \Lsh_N)$, since $\Gamma$ is not contained in the gap divisor $\Omega$.  Let $f \in \bbar{R}_N \smallsetminus H^0(X, \sh{I}_{\Gamma} \Lsh_N)$ and let $F = \div_X(f) + \Delta_N$.  The germ of $F$ at $p_b$ is in $\mf{d} \smallsetminus y \struct$ and is thus equal to $xr + ys$ for some $r \in \struct \smallsetminus y \struct$ and $s \in \struct$.
Recall 
\[ \deg_{\Gamma}(F \cap \Gamma) = \deg_{\Gamma} (\Delta_N \cap \Gamma) =  \deg_{\Gamma}(\Pi_N \cap \Gamma)\]
by \eqref{case2}.  Since  
\[f \in H^0(X,  \sh{A}' \sh{D}^{\sigma} \cdots \sh{D}^{\sigma^{N-1}}(\sh{C}')^{\sigma^N} \Lsh_N)\]
and $b \geq 1$,
 we see that $F$ does not vanish at $p_{b+Nj}$ unless $j = 0$.

  Since $y^d \in \mf{d}$,  there is some $h \in \bbar{R}_N$ so that the germ of $H := \div_X(h) + \Delta_n$ at $p_b$ is equal to $y^d$.    Thus we have $H = d \Gamma + C$, where $C(p_b) \neq 0$.

Let $ m \geq 2$.  For $i = 1 \ldots m-1$, define $\gamma_i \in \bbar{R}_{Nm}$ by 
\[ 
\gamma_i := f f^{\sigma^N} \cdots f^{\sigma^{N(i-1)}} h^{\sigma^{Ni}} f^{\sigma^{N(i+1)}} \cdots f^{\sigma^{N(m-1)}}.
\]
Then 
\begin{multline*}
\div_X(\gamma_i) + \Delta_{Nm} = \\
  F + \sigma^{-N} (F) +  \cdots + \sigma^{-N(i-1)}(F) + \sigma^{-Ni}(H) + 
 \sigma^{-N(i+1)} (F) + \cdots + \sigma^{-N(m-1)}(F) \\
  = d \Gamma + F +   \cdots + \sigma^{-N(i-1)}(F) + \sigma^{-Ni}(C) + 
 \sigma^{-N(i+1)} (F) + \cdots + \sigma^{-N(m-1)}(F).
 \end{multline*}
 
Fix $1 \leq i \leq m-1$.   The local equation of $\div_X(\gamma_i) + \Delta_{Nm}$ at $p_{b+Ni}$ is equal to $y^d \eta$ for some local unit $\eta$.  On the other hand, if $j \neq i$, the local equation of 
 $\div_X(\gamma_j) + \Delta_{Nm}$ at $p_{b+Ni}$ is equal to $(xr + ys) y^d \beta$ for some $0 \neq \beta \in \struct$.  In particular, we see that modulo $y^{d+1}$, the set 
 \[\{ \gamma_j \st j \neq i\}\]
 does not generate $\gamma_i$.  Thus the elements
\[ \{ \gamma_i z^{Nm} + J(d+1) \st 1 \leq i \leq m-1 \} \subseteq A_{Nm}\]
are linearly independent, and $\dim_{\kk}(A_{Nm}) \geq m$.  Therefore, $\GK A \geq 2$.

Suppose that $A$ is noetherian.     Let $\sh{S}_n$ be the image of the natural map $\sh{R}_n \otimes \struct_{\Gamma} \to \Lsh_n |_{\Gamma}$.  Let
 \[ B := R/J(1) \subseteq \bigoplus_{n \geq 0} H^0(\Gamma, \sh{S}_n) z^n.\]
Recall from \eqref{case2} that 
 \[\deg_{\Gamma} (W_n \cap \Gamma)+ \deg_{\Gamma}(\Omega \cap \Gamma) = \deg_{\Gamma}(\Lsh_n |_{\Gamma})\]
 for all $n \geq 1$.  
  That is, 
 \[
 \deg_{\Gamma}(\sh{S}_n ) = \deg_{\Gamma}(\Lsh_n |_{\Gamma}) - \deg_{\Gamma}(\Omega \cap \Gamma + W_n \cap \Gamma ) = 0
 \] 
 for all $n \geq 1$, and $\dim H^0(\Gamma, \sh{S}_n) \leq 1$ for all $n$.  Therefore,  $\GK B \leq 1 $.

 Each  $J(i)/J(i+1)$ is a finitely generated $A$-module.  The $A$-action on $J(i)/J(i+1)$ factors through $B$.  Thus each $J(i)/J(i+1)$ is also a finitely generated $B$-module and thus has GK-dimension $\leq 1$.  The module $A_A$ has  a finite filtration by modules of the form $J(i)/J(i+1)$;  therefore  $\GK A \leq 1$.  This gives a  contradiction.   

Thus every point in $Z \cup \Lambda \cup \Lambda'$ must have a dense orbit.
 \end{proof}

\begin{theorem}\label{thm-data-CT}
Assume Assumption-Notation~\ref{assnot6}.  Then the surface data 
\[\surfD = \surfdata\]
 is transverse, and $T(\surfD)$ is a finitely generated left and right $R$-module.  Further, $T(\surfD)$ and $\sh{T}(\surfD)$ are noetherian.
\end{theorem}
\begin{proof}
We have seen in  Lemma~\ref{lem-omega-transverse} that  $\{\sigma^n \Omega\}$ is critically transverse and in Proposition~\ref{newprop} that  all points in $Z \cup \Lambda \cup \Lambda'$ have dense $\sigma$-orbits.   
By Lemma~\ref{lem-ample} the sequence of bimodules $\{ (\sh{R}_n)_{\sigma^n} \}$ is left and right ample.  Thus by Lemma~\ref{lem-fg}, $T:=T(\surfD)$ is a finitely generated left and right $R$-module.  Thus $T$ is noetherian.  By Theorem~\ref{thm-surfprop-main} the data $\surfD$ is transverse and $\sh{T}(\surfD)$ is a noetherian bimodule algebra.
\end{proof}

\begin{corollary}\label{cor-normalize}
Let $R$ be a birationally commutative projective  surface.  Then there is transverse surface data 
\[ \surfD = \surfdata,\]
where $X$ is normal and $\sigma$ is numerically trivial, and an integer $\ell$ so that
\[ R \ver{\ell} \subseteq T(\surfD)\]
and $T (\surfD)$ is a finitely generated  left and right module over $R \ver{\ell}$.
\end{corollary}
\begin{proof}
By Theorem~\ref{thm-data}, there are a positive integer $\ell$ and   \adm\ surface data $\surfD = \surfdata$ so that  
\[\sh{R}(X)\ver{\ell} = \sh{T} (\surfD).\]
By possibly increasing $\ell$, we may also assume that $\bbar{R}_1$ generates $K = \kk(X)$.  Thus Assumption-Notation~\ref{assnot6} holds for $R\ver{\ell}$.  By Theorem~\ref{thm-data-CT}, $T(\surfD)$ is a finitely generated left and right $R \ver{\ell}$-module, and  the data $\surfD$ is in fact transverse.
 \end{proof}

\section{Constructing $Y$} \label{FINDY}

Let us review our progress towards proving 
  Theorem~\ref{thm-surfclass}.   In Theorem~\ref{thm-data}, we constructed surface data $\surfD = \surfdata$ for an appropriate bimodule algebra $\sh{R}$ associated to $R$; in Corollary~\ref{cor-normalize} we showed that this data is actually transverse, and that $T (\surfD)$ is a finite left and right module over some $R \ver{k}$.  In some sense, we may think of $T(\surfD)$ as a normalization of $R$, or more properly of $R\ver{k}$; note that the  variety $X$  given in Theorem~\ref{thm-data} is normal.  

Of course, there is no guarantee that $R$ is really associated to a normal variety.  Modifying $X$ to find the true surface associated to  $R$ turns out to be  quite technical.  In this section, we construct the scheme $Y$ on which $R$ actually lives by carefully studying the rational maps on $X$ defined by the sections in $\bbar{R}_n$.

We will assume that we are in the situation of Assumption-Notation~\ref{assnot6}.    By Theorem~\ref{thm-data-CT} the data $\surfD$ is transverse.   

\begin{notation}\label{not-FINDY}
We establish notation that we will use throughout the  section.  
For any $n \geq 1$, $\bbar{R}_n$ defines a rational map 
\[ \xymatrix{
X \ar@{-->}[r]^{\beta_n} & \PP^N } \]
that is birational onto its image.   Let $Y_n$ be the closure of the image of $X$; we also use $\beta_n$ to denote the induced birational map from $X$ to $Y_n$.  If we let $\alpha'_n:  X'_n \to X$ be  the blowup of $X$ at the base  locus $W_n$ of $\bbar{R}_n$, then by \cite[Example~II.7.17.3]{Ha} there is a birational morphism $\gamma'_n: X'_n \to Y_n$ such that the diagram
\[ \xymatrix{
X'_n \ar[d]_{\alpha'_n} \ar[rd]^{\gamma'_n} & \\
X \ar@{-->}[r]_{\beta_n}		& Y_n }\]
commutes.   
Let $\lambda_n:X_n \to X'_n$ be the normalization of $X'_n$, and let $\alpha_n := \alpha'_n \lambda_n$ and $\gamma_n := \gamma'_n \lambda_n$.  Thus we have
\[ \xymatrix{
X_n \ar[d]_{\alpha_n} \ar[rd]^{\gamma_n} & \\
X \ar@{-->}[r]_{\beta_n}		& Y_n }\]
for all $n \geq 1$.  Note, that as $X$ is normal, $\alpha_n^{-1}$ is defined at all points in $X \smallsetminus W_n$.

Let $W$ be the set-theoretic cosupport of $\sh{ADC}$ and let 
\[\mathbb{W} := \bigcup_{p \in W} O(p).\]
That is, $\mathbb{W}$ is  the union of the finitely many (dense) orbits that meet some $W_n$.  Let $U := X \smallsetminus \mathbb{W}$.   Note that all $\beta_n$ are defined at all points in $U$.
For any $n \geq 1$, let $E_n$ be the exceptional locus of $\alpha_n$.  Then $\alpha_n$ induces an isomorphism from $X_n \smallsetminus E_n \to X \smallsetminus W_n$.  Let $U_n := \alpha_n^{-1}(U)$.  We caution that $U_n$ is {\em not} $X_n \smallsetminus E_n$.

If $N > n \geq 1$,   let $\pi^N_n: Y_N \to Y_n$
 be the birational map induced from the multiplication $\bbar{R}_n (\bbar{R}_{N-n})^{\sigma^n} \subseteq \bbar{R}_N$ and Lemma~\ref{lem-definemap}, with $E = \sigma^{-n}(\Omega)$.   That is, the diagram of birational maps
 \[ \xymatrix{ X \ar@{-->}[r]^{\beta_N} \ar@{-->}[rd]_{\beta_n} &  Y_N \ar@{-->}[d]^{\pi^N_n} \\
  & Y_n}\]
  commutes, and for any $x \in U \smallsetminus \sigma^{-n}(\Omega)$, $\pi^N_n$ is defined at $\beta_N(x)$.  
Likewise, 
the multiplication $\bbar{R}_{N-n} (\bbar{R}_n)^{\sigma^{N-n}} \subseteq \bbar{R}_N$ gives a commuting diagram of birational maps 
 \[ \xymatrix{ 
 X \ar@{-->}[r]^{\beta_N} \ar@{-->}[rd]_{\beta_n \circ \sigma^{N-n}} &  Y_N \ar@{-->}[d]^{\rho^N_n} \\
  & Y_n.}\]
  The map $\rho^N_n$ is defined at $\beta_N(x)$ if $x \in U \smallsetminus \sigma^{-(N-n)}(\Omega)$.
  \end{notation}

We record for future reference an elementary lemma on birational maps.

\begin{lemma}\label{lem-INVERSE}
Let $\beta: X \to Y$ be a birational map of projective varieties that is defined  and is a local isomorphism at $x \in X$; let $y := \beta(x)$.  Then $\beta^{-1}$ is defined at $y$; in particular, if $x' \in X$ with $\beta(x') = y$, then $x' = x$.
\end{lemma}
\begin{proof}
This is almost tautological.  The fact that $\beta$ induces an isomorphism between the local rings $\struct_{Y,y}$ and $\struct_{X, x}$ means that there are open neighborhoods $x \in V \subseteq X$ and $y \in V' \subseteq Y$ so that $\beta$ restricts to an isomorphism between $V$ and $V'$.  This means that $\beta^{-1}$ gives a well-defined map 
\[ V' \to V \subseteq X.\]
This precisely says that the birational map $\beta^{-1}: Y \to X$ is defined at $y$.
\end{proof}

For all $n   \geq 1$, let $A_n$ be the set of points $p \in U$ such that $\beta_n$ is not a local isomorphism at $p$; that is, the set of $p$ such that the induced map from $\struct_{Y_n, \beta_n(p)} \to \struct_{X, p}$ is not an isomorphism.  Write $A_n$ as the disjoint union $A_n = C_n \sqcup Q_n \sqcup P_n$, where $C_n$ is the intersection of a curve in $X$ with $U$, $Q_n$ is 0-dimensional and supported on points of infinite order under $\sigma$, and $P_n$ is 0-dimensional and supported on points of finite order.   By assumption on the cardinality of $\kk$, any curve in $X$ must meet $U$ in uncountably many points, and so the sets $C$, $P$, and $Q$ are well-defined.

\begin{proposition}\label{prop-stablesubscheme}
There is some $m_1$ such that $A_{m_1}$ is $\sigma$-invariant  and $A_n = A_{m_1}$ for all $n \geq m_1$; further, $\bbar{C}_{m_1} \subseteq U$ and $Q_{m_1} = \emptyset$. 
\end{proposition}
\begin{proof}
 Let $y \in U$ and let $N > n \geq 1$.  If $\pi^N_n$ is defined at $\beta_N(y)$ and $\beta_n$ is a local isomorphism at $y$, then from the inclusions
 \[ \struct_{Y_n, \beta_n(y)} \subseteq \struct_{Y_N, \beta_N(y)} \subseteq \struct_{X,y}\]
  clearly $\beta_N$ is a local isomorphism at $y$.  As $\pi^N_n$ is defined on $\beta_n(U \smallsetminus \sigma^{-n} \Omega)$, we see that  
\[ A_N \subseteq A_n \cup \sigma^{-n} \Omega.\]
Making the same argument with the map $\rho^N_n$, we obtain that
\[ A_{N} \subseteq \sigma^{-(N-n)}(A_n) \cup \sigma^{-(N-n)} \Omega.\]
Thus 
\beq \label{GEO}
A_{n+m} \subseteq \sigma^{-m} (\Omega \cup A_n) \cap (A_n \cup \sigma^{-n} \Omega)
\eeq
for any $n, m \geq 1$.  
In particular, 
\[ C_{n+m} \subseteq \sigma^{-m} (\Omega \cup C_n) \cap (C_n \cup \sigma^{-n} \Omega)\]
for all $n, m \geq 1$.

Fix $n \geq 1$.  Then 
 $\sigma^{-n} \Omega \cap \sigma^{-m} \Omega$ is finite for $m \gg 0$.  Further,  $C_n \cap \sigma^{-m} \Omega$ and $\sigma^{-m} C_n \cap \sigma^{-n} \Omega$ are finite for $m \gg 0$. 
Thus for $m \gg 0$, we have
\[ C_{n+m} \subseteq C_n \cap \sigma^{-m} C_n.\]
Either $C_{n+m} = 0$ for $m \gg 0$, or $C_n$ contains a $\sigma$-invariant curve $C'$ so that $C_{n+m} = C'$ for $m \gg 0$.  In either case, 
 there is some $n_1$ such that if $n \geq n_1$, then $C_n = C_{n_1}$ and  $\sigma(C_{n_1}) = C_{n_1}$.  Let $C := C_{n_1}$.    

Now $C$ is a curve in $U$; let $\bbar{C}$ be its closure in $X$.  Then $\bbar{C}$ is also $\sigma$-stable, and  since all orbits in $\mathbb{W}$ are Zariski-dense in $X$, we have  $\mathbb{W} \cap \bbar{C} = \emptyset$.  Thus $C = \bbar{C} \subseteq U$.

For $n, m  \geq n_1$ we have
\begin{align*}
 Q_{n+m} \subseteq A_{n+m} & \subseteq \sigma^{-m}(\Omega \cup C \cup Q_n \cup P_n) \cap (C \cup Q_n \cup P_n \cup \sigma^{-n}\Omega) \\
& = C \cup \bigl( \sigma^{-m}(\Omega \cup Q_n \cup P_n) \cap (Q_n \cup P_n \cup \sigma^{-n} \Omega) \bigr).
\end{align*}
As 
\[Q_n \cap C = Q_n \cap \sigma^{-m}(C) = \emptyset\]
 for $n \geq n_1$, and $Q_n \cap \sigma^k(P_m) = \emptyset$ for all $n, m, k$ by definition,
we see that 
\[ Q_{n+m} \subseteq  \sigma^{-m} (\Omega \cup Q_n) \cap (Q_n \cup \sigma^{-n} \Omega)\]
for $n \geq n_1$ and $m \geq 1$.  

Choose $k$ such that 
\beq \label{choosek}
\Omega \cap \sigma^{-1} \Omega \cap \cdots \cap \sigma^{-k} \Omega = \emptyset.
\eeq
  Such $k$ exists because, by transversality of the data $\surfD$,  $\Omega$ contains no forward $\sigma$-orbits. Choose $n_2 \geq n_1$ such that if $n \geq n_2$, then we have for $i = 0, \ldots, k$ that
\[ 
Q_{n_1 +i} \cap \sigma^{-(n-(n_1+i))} (Q_{n_1+i} \cup \Omega) = \emptyset.\]
We may do this because each finite set $Q_{n_1+i}$ is supported on infinite orbits and  $\{ \sigma^n\Omega\}$ is critically transverse.

  We have  for $i = 0, \ldots, k$ and for $r \geq n_2$ that 
\begin{align*}
Q_r & \subseteq \sigma^{-(r-(n_1+i))}(Q_{n_1+i} \cup \Omega) \cap (Q_{n_1+i} \cup \sigma^{-(n_1+i)} \Omega) \\
	& \subseteq \sigma^{-(n_1+i)} \Omega
	\end{align*}
and so $Q_r = \emptyset$ for $r \geq n_2$ by \eqref{choosek}.

Finally,  $\Omega$ does not contain any points of finite order, and by construction $P_n$ is disjoint from $C$ and from $Q_n$.  Thus \eqref{GEO} implies that 
\[ P_{n+m} \subseteq P_n \cap \sigma^{-m} P_n\]
 for all $n, m \geq n_1$.   Thus there is some $n_3 \geq n_2$ such that $ P_{n_3}$ is $\sigma$-invariant, and  $P_n = P_{n_3}$ if if $n \geq n_3$.   The result is proven for $m_1 := n_3$.  
\end{proof}

\begin{notation}\label{not-FINDY2}
Let $m_1$ be the integer given by Proposition~\ref{prop-stablesubscheme}.  Let $A := A_{m_1}$, $C := C_{m_1}$, and $P := P_{m_1}$.   Recall that $Q_{n} = \emptyset$ for all $n \geq m_1$.  
\end{notation}

\begin{corollary}\label{cor-U-finite}
 For $n \geq m_1$, the only curves in $X_n$ that are contracted by $\gamma_n$ are contained in the exceptional locus $E_n$ of $\alpha_n$.  In particular, the map $\gamma_n$ is finite at all points of $U_n$, and $\beta_n$ is finite at all points of $U$.
\end{corollary}
\begin{proof}  
Suppose that $n \geq m_1$ and that 
$\gamma_n$ contracts some irreducible curve $\Gamma$ that is not contained in $E_n$.  By assumption on the cardinality of $\kk$, $\Gamma$ meets $U_n$.  By construction,  $\Gamma \cap U_n\subseteq \alpha_n^{-1}C_n  = \alpha_n^{-1}C$.  Now, as $\alpha_n$ is an isomorphism away from $E_n$, the curve $\alpha_n^{-1}(C)$ is closed in $X_n$.  
Thus  
\[ \Gamma = \bbar{\Gamma \cap U_n}  \subseteq \bbar{\alpha_n^{-1}(C)} = \alpha_n^{-1} C.\]  
This means that $\alpha_n(\Gamma) \subseteq C \subset U$, so $\alpha_n(\Gamma)$ is disjoint from $\Bs(\bbar{R}_n)$.   As $\beta_n$ contracts the curve $\alpha_n(\Gamma)$,  Lemma~\ref{lem-contractedcurves} implies that $\Delta_n . \alpha_n(\Gamma) = 0$.  This contradicts the ampleness of $\Delta_n$ by the Nakai-Moishezon criterion (\cite[Theorem~V.I.10]{Ha}; see \cite[Theorem~1.2.23]{Laz} for a reference that includes singular surfaces).   As $\alpha_n$ is a local isomorphism at all points in $U_n$, the statement on $\beta_n$ follows immediately.
\end{proof}

We recall some terminology from commutative algebra.  Let $R$ be a commutative noetherian $\kk$-algebra, and let $T$ be its normalization.  Recall that the {\em $S_2$-ification} of $R$ is the unique minimal $\kk$-algebra $S \subseteq T$ such that $R \subseteq S$ and $S$ satisfies Serre's condition $S_2$.  More explicitly, 
\[ S = \bigcap_{\substack{ P \in \Spec R \\ 
					\rm{ht}\, P = 1 
			}}
					R_{P}.\]
See \cite[Definition~2.2.2(ii)]{Kollar1985} and subsequent discussion. 

We give a lemma on the domain of definition of birational maps of $S_2$-ifications.

\begin{lemma}\label{lem-S2}
Let $T$ be a normal commutative domain that is a  finitely generated $\kk$-algebra, and let $R, R' \subseteq T$ be finitely generated subalgebras so that $T$ is the normalization of both $R$ and $R'$.  Let $S$, respectively $S'$, be the $S_2$-ification of $R$, respectively $R'$.  Suppose that the induced birational map 
\[ \xymatrix{ \pi: \Spec R \ar@{-->}[r] & \Spec R' }\]
is defined away from a locus of codimension 2.  Then the induced birational map
\[ \xymatrix{ \zeta: \Spec S \ar@{-->}[r] & \Spec S' }\]is defined everywhere; that is, $S' \subseteq S$.
\end{lemma}
\begin{proof}
Because $\pi$ is defined in codimension 2, for every height 1 prime $P$ of $R$, $\pi$ is defined at the generic point of $V(P) \subseteq \Spec R$.   That is, for every height 1 $P$, we have $R' \subseteq R_P$.  Thus
\[ R' \subseteq  \bigcap_{\substack{ P\in \Spec R \\ 
					\rm{ht}\, P= 1 
			}}
					R_{P} = S .\]
By minimality of the $S_2$-ification, $S' \subseteq S$.  
\end{proof}

\begin{proposition}\label{prop-stableY} If $n \gg 0$, then the rational maps $\pi^{n+1}_n, \rho^{n+1}_n: Y_{n+1} \to Y_n$  are local isomorphisms everywhere on the image of $U$; in particular, they are defined everywhere on the image of $U$.
\end{proposition}

\begin{proof}
We continue to let $m_1$, $A$, $C$, and $P$ be as in  Notation~\ref{not-FINDY2}.  In particular, if $n \geq m_1$ then $\gamma_n:  X_n \to Y_n$ is a local isomorphism at all points in $U_n \smallsetminus \alpha_n^{-1}(C\cup P)$.  

Let $n \geq m_1$.  Recall that if $x \in U \smallsetminus \sigma^{-n}(\Omega)$, then $\pi^{n+1}_n$ is defined at $\beta_{n+1}(x)$. As $\sigma^{-n}(\Omega)$ contains no points of finite order, $\pi^{n+1}_n$ is defined at all points in $\beta_{n+1}(P)$.   We saw in Corollary~\ref{cor-U-finite} that $\beta_m$ is finite at all $p \in U$, and in particular, at all $p \in P$.  By finiteness of the integral closure, there is some $m_2 \geq m_1$ so that if $n \geq m_2$, then $\pi^{n+1}_n$ is a local isomorphism at all points in $\beta_{n+1}(P)$.  

Let $n \geq m_2$.     Now, let $x \in U \smallsetminus A$.
  Then $\beta_{n+1}$ and $\beta_n$ are local isomorphisms at $x$, and thus by Lemma~\ref{lem-INVERSE}, $\beta_{n+1}^{-1}$ is defined at $\beta_{n+1}(x)$.  Thus
$\pi^{n+1}_n = \beta_n \beta_{n+1}^{-1}$ is defined and is a local isomorphism at $\beta_{n+1}(x)$.

The only points where $\pi^{n+1}_n$ may not be defined thus lie in $\beta_{n+1}(C \cap \sigma^{-n}(\Omega))$.  
If $x \in U \smallsetminus (\sigma^{-n}(\Omega) \cap C)$, then  $\pi^{n+1}_n$ is defined and is a local isomorphism at $\beta_{n+1}(x)$.   

The intersection $\sigma^{-n} (\Omega) \cap C$ is finite by transversality of $\Omega$, since $C$ is $\sigma$-invariant.  
%
Thus  $\pi^{n+1}_n$ is defined at  the generic point of each component of of $\beta_{n+1} (C)$.  As $\gamma_n$ is finite at each point of $U_n$, by finiteness of the integral closure there is $m_3 \geq m_2$ such that if $n \geq m_3$, then $\pi^{n+1}_n$ is an isomorphism at the generic point of each component of $\beta_{n+1}(C)$.  

For each $n \geq m_3$, let $\delta_n : Z_n \to Y_n$ be the projective variety obtained by taking the $S_2$-ification of $Y_n$ at all points in  $\beta_n(C)$.  Since   $X_n$ is normal, the birational morphism $\gamma_n: X_n \to Y_n$ factors through $Z_n$.  There are thus a  birational morphism $\epsilon_n: X_n \to Z_n$ and a birational map $\eta_n: X \to Z_n$ so that the diagram
\[ \xymatrix{
 & X_n \ar[ld]_{\alpha_n} \ar[d]^{\epsilon_n} \ar[rd]^{\gamma_n}  && \\
X \ar@{-->}[r]^{\eta_n} \ar @/_10pt/ @{-->}[rr]_{\beta_n}	& Z_n \ar[r]^{\delta_n}	& Y_n}\]
commutes.  In particular, $\eta_n$ is defined at all points of $U$.

Let  $\zeta^{n+1}_n: Z_{n+1} \to Z_n$ be the  induced birational map such that the diagram
\beq \label{S2} 
\xymatrix{
X \ar @{-->}[r]^{\eta_{n+1}} \ar @{-->}[rd]_{\eta_n} \ar@/^18pt/@{-->}[rr]^{\beta_{n+1}}
&  Z_{n+1} \ar[r]^{\delta_{n+1}} \ar@{-->}[d]^{\zeta^{n+1}_n} 	& Y_{n+1} \ar@{-->}[d]^{\pi^{n+1}_n} \\
 &  Z_n \ar[r]_{\delta_n}	& Y_n }
 \eeq
 commutes.  
 We claim that for $n \geq m_3$, that  the rational map $\zeta^{n+1}_n$ is defined at all points of $\eta_{n+1} (U)$.  
 
 Let $p \in U \smallsetminus (C\cup P)$.  Then $\beta_{n+1}$ and therefore $\eta_{n+1}$ is a local isomorphism at $p$, so, using Lemma~\ref{lem-INVERSE},  $\zeta^{n+1}_n = \eta_n \eta_{n+1}^{-1}$ is defined at $\eta_{n+1}(p)$.
 
 If $p \in P$ and $\beta_{n+1}(p) \not \in \beta_{n+1}(C)$, then by our
 choice of $n$, the map $\pi^{n+1}_n$ is a local isomorphism at $\beta_{n+1}(p)$.  Thus $\beta_{n}(p) \not \in \beta_n(C)$.  By construction $\delta_n$ is a local isomorphism at  $\eta_n(p)$ and so $\zeta^{n+1}_n = \delta_n^{-1} \pi^{n+1}_n \delta_{n+1}$ is defined at $\eta_{n+1}(p)$.

 Now let  $p \in C$.   We have seen that $\pi^{n+1}_n$ is defined on $\beta_{n+1}(U)$, except at a 0-dimensional locus contained in $\beta_{n+1}(C)$.  It follows from Lemma~\ref{lem-S2} that $\zeta^{n+1}_n$ is defined at $\eta_{n+1}(p)$.   This completes the proof of the claim.
 
 %
 
Using  finiteness of the integral closure again, we may choose $m_4 \geq m_3$ so that $\zeta^{n+1}_n$ is a local isomorphism at all points of $\eta_{n+1}(U)$   for $n \geq m_4$.  For $n \geq m_4$, let 
\[F_n  = \{ p \in C \st \delta_n \mbox{ is not a local isomorphism at } \eta_n(p) \}.\]
  As $\delta_n$ is finite, the  set $F_n$ is finite.  Arguing as in the proof of Proposition~\ref{prop-stablesubscheme}, and using the maps $\zeta^{n+m}_n$, we have that
\[F_{n+m} \subseteq \bigl(F_n \cup \sigma^{-n}( \Omega \cap C)\bigr) \cap  \sigma^{-m}\bigl( F_n \cup  (\Omega \cap C)\bigr) \]
if $n \geq m_4$.  
Since $\Omega \cap C$ consists of finitely many points of infinite order, for $n, m \gg 0$ we have that
\[
\sigma^{-n} (\Omega \cap C) \cap \sigma^{-m}(F_n \cup (\Omega \cap C)) = F_n \cap \sigma^{-m}(\Omega \cap C)  = \emptyset.\]
Thus for $m, n \gg 0$, we have that 
\[ F_{n+m} \subseteq F_n \cap \sigma^{-m} F_m,\]
and so for $n \gg 0$, we have that $F_n = F_{n+1}$ is $\sigma$-invariant.  In particular,  $F_n \cap \sigma^{-n} \Omega = \emptyset$.  This means that $\delta_{n+1}$ is a  local isomorphism at all points of $ \eta_{n+1}(\sigma^{-n} \Omega \cap U)$.   By  \eqref{S2}, $\pi^{n+1}_n = \delta_n \zeta^{n+1}_n \delta_{n+1}^{-1}$ is defined everywhere in $\beta_{n+1}(\sigma^{-n} \Omega \cap U)$.  Thus $\pi^{n+1}_n$ is defined everywhere in $\beta_{n+1}(U)$ for $ n \geq m_4$.

The argument that for $n \gg 0$, $\rho^{n+1}_n$ is defined everywhere on $\beta_{n+1}(U)$ is completely symmetric.  By finiteness of the integral closure, we see that for $n \gg 0$ both $\pi^{n+1}_n$ and $\rho^{n+1}_n$ are local isomorphisms at every point of $\beta_{n+1}(U)$.
\end{proof}

We establish some more notation, which we will use in the next few results.

\begin{notation}\label{not-maps}
Assume Assumption-Notation~\ref{assnot6}.  
Let $m$ be such that for $n \geq m-1$, the rational maps $\pi^{n+1}_n$ and $\rho^{n+1}_n$ are defined and are local isomorphisms at every point in $\gamma_{n+1}(U)$.  We call $Y_m$ a {\em stable scheme} for $R$.   Let $C \subset U$ be the $\sigma$-invariant curve where $\gamma_m$ is not a local isomorphism, and let $P \subset U$ be the $\sigma$-invariant 0-dimensional subscheme where $\gamma_m$ is not a local isomorphism.  
Let $F $ be the  $\sigma$-invariant subset of $C$ that maps onto points where $Y_m$ does not satisfy $S_2$.  

For all $n \geq 2$, we define a birational automorphism $\tau_n$ of $Y_n$ by setting 
\[\tau_n := (\pi^n_{n-1})^{-1} \rho^n_{n-1}.\]
  By construction,  $\tau_n \gamma_n = \gamma_n \sigma $ as birational maps.   Proposition~\ref{prop-stableY} implies that for $n \geq m$, $\tau_n$ is an automorphism of $\gamma_n(U_n)$.
  
  For all $n, m \geq 1$, we  define birational maps 
\[ p^{n+m}_n, r^{n+m}_n:  X_{n+m} \to X_n,\] 
where $p^{n+m}_n := (\alpha_n)^{-1} \alpha_{n+m}$ and
$r^{n+m}_n := \alpha_n^{-1} \sigma^m \alpha_{n+m}$.
By construction, we have 
\[ \pi^{n+m}_n \gamma_{n+m} = \gamma_n p^{n+m}_n\]
and
\[ \rho^{n+m}_n \gamma_{n+m} = \gamma_n r^{n+m}_n\]
as birational maps from $X_{n+m}$ to $Y_n$.  
\end{notation}

  Recall that $W$ is the cosupport of $\sh{ADC}$, that $\mathbb{W} = \bigcup_{p \in W} O(p)$ and that $U = X \smallsetminus \mathbb{W}$.  
The map $\beta_m$ is defined and finite at every point of $U$.  Heuristically, it may fold $U$ along $C$ or pinch $U$ into a cusp at some point of $P$.  At points of $F$, $\beta_m$ does some additional cusping, since there $Y_m$ fails $S_2$.
We will construct a finite morphism $\theta: X \to Y$ by cusping and folding $U$ along $C$, $P$, and $F$, and gluing in the points of $\mathbb{W}$.

There is one technicality still to dispose of: in order to glue as described, we need the sets
\[ \gamma_m(U_m) = \beta_m(U)\]
and 
\[ \gamma_m (\alpha_m^{-1}(\mathbb{W}))\]
to be disjoint, at least for large $m$.  Proving this is the content of the next few results.

Ideally, this result would be a consequence of some sort of ampleness statement, since we are trying to show, roughly speaking, that the $\gamma_m$ separate points of $X_m$ for $m$ large.  Unfortunately, we have not been able to prove this directly.  Instead, we  analyze in detail   
 how our various birational maps affect $\mathbb{W}$.  To do this, we establish  more notation.    

\begin{notation}\label{not-nice}
Fix $w \in W$.   For any $n \in \ZZ$,  let $w_n := \sigma^{-n}(w)$.  Let $\struct = \struct_{X, w}$ and identify all $\struct_{X, w_n}$ with $\struct$ as usual.  Let $\mf{a}:= (\sh{I}_{\Omega}\sh{A})_w$, let $\mf{c}:= \sh{C}_w$, and let $\mf{d}:= \sh{D}_w$.  

Let $\mf{m}^n_i$ be the germ of $\sh{R}_n \Lsh_n^{-1}$ at $w_i$, regarded as an ideal in $\struct$, as usual. By assumption, we have:
 \begin{itemize}
\item $\mf{m}^n_0 =  \mf{a}$;
\item if $1 \leq i \leq n-1$ then $\mf{m}^n_{i} = \mf{d}$;
\item $\mf{m}^n_{n} = \mf{c}$; and
\item if $i < 0$ or $i > n$ then $\mf{m}^n_{i+j} = \struct$
\end{itemize}
for all $n \geq 1$.

We define subschemes 
\[\Ea := \alpha_2^{-1}(w_0),\]
\[\Ed := \alpha_2^{-1}(w_1),\]
and
\[\Ec := \alpha_2^{-1}(w_2)\]
of $X_2$.   
That is,  $\Ed$ is the exceptional locus obtained by blowing up the ideal $\mf{d}$ and normalizing, and similarly for $\Ea$ and $\Ec$.  

For any $ n \geq 1$ and $i \in \ZZ$, define $E^n_i \subseteq Y_n$ by  
\beq\label{def-Eni}
E^n_i := \gamma_n(\alpha_n^{-1} (w_i)).
\eeq
\end{notation}

\begin{lemma}\label{lem-nice}
Let $w \in W$.  
Then for all $n, m \geq 1$ the map $p^{n+m}_n$ is defined at all points in 
$\alpha_{n+m}^{-1} (O(w) \smallsetminus \{ w_{n}\})$, and is a local isomorphism
at all points in $\alpha_{n+m}^{-1}(O(w) \smallsetminus \{w_n, w_{n+1}, \ldots, w_{n+m}\})$.
Likewise, $ r^{n+m}_n$ is defined 
at all points in $\alpha_{n+m}^{-1}(O(w) \smallsetminus \{w_m\})$ and is a local isomorphism at all points in $\alpha_{n+m}^{-1} (O(w) \smallsetminus \{ w_0, \ldots, w_m \})$.
\end{lemma}
\begin{proof}
Fix $n, m \geq 1$ and let $w \in W$.  
By definition, $\Supp W_n \subseteq \{ w_0, \ldots, w_n\}$, and so the map $\alpha_n$ is a local isomorphism at all points in $\alpha_n^{-1} (O(w) \smallsetminus \{ w_0, \ldots, w_n\})$.  Furthermore, as $\mf{m}^n_0 = \mf{a}$, we have that $\alpha_n^{-1}(w_0) \cong \Ea$.  Likewise, $\alpha_n^{-1}(w_n) \cong \Ec$, and for $1 \leq i \leq n-1$, $ \alpha_n^{-1} (w_i) \cong \Ed$.  

Therefore, if $i< 0$ or $i > n+m$, then $p^{n+m}_n = \alpha_n^{-1} \alpha_{n+m}$ is defined and is a local isomorphism at the point $\alpha_{n+m}^{-1}(w_i)$.  For $0 \leq i \leq n-1$, the stalks 
$\mf{m}^{n+m}_i$ and $\mf{m}^n_i$ 
 are isomorphic.  Thus  $\alpha_n^{-1} \alpha_{n+m}$ extends to a map that is defined and a local isomorphism at all points of $\alpha_{n+m}^{-1}(w_i)$.   For $n+1 \leq i \leq n+m$, $\alpha_n^{-1}$ is defined at $w_i$, so $p^{n+1}_n$ is defined on $\alpha_{n+m}^{-1}(w_i)$, although it is not necessarily a local isomorphism.

We repeat this analysis for the maps $r^{n+m}_n$.  If $i < 0$ or $i > n+m$, then $\alpha_{n+m}$ is a local isomorphism at the point $\alpha_{n+m}^{-1}(w_i)$, and $\alpha_n^{-1}$ is defined (and is thus a local isomorphism) at $w_{i-m} = \sigma^m(w_i)$.    Thus $r^{n+m}_n = \alpha_n^{-1} \sigma^m \alpha_{n+m}$ is a local isomorphism at $\alpha_{n+m}^{-1}(w_i)$.   If $m+1 \leq i \leq n+m$, then $\alpha_{n+m}^{-1} (w_i) \cong \alpha_n^{-1} (w_{i-m})$ and $r^{n+m}_n$ extends to  a local isomorphism at all points in $\alpha_{n+m}^{-1}(w_i)$.  Finally, if $0 \leq i \leq m-1$, then $\alpha_n^{-1}$ is defined at $w_{i-m}$ and so $r^{n+m}_n$ is defined on $\alpha_{n+m}^{-1} (w_i)$.
\end{proof}

\begin{lemma}\label{lem-W}
 Let $w \in W$.

$(1)$  For $m, n \geq 1$, if $i \neq n$ then $\pi^{n+m}_n$ is defined at all points in 
$E^{n+m}_i = \gamma_{n+m}(\alpha_{n+m}^{-1} (w_i))$, and $\pi^{n+m}_n(E^{n+m}_i) = E^n_i$.

$(2)$  For $m, n \geq 1$, if $i \neq m$  then $\rho^{n+m}_n$ is defined at all points in 
$E^{n+m}_i = \gamma_{n+m}(\alpha_{n+m}^{-1} (w_i))$, and $\rho^{n+m}_n(E^{n+m}_i) = E^n_{i-m}$.
\end{lemma}

\begin{proof}
(1)  On $X_{n+m}$, the rational functions in $\bbar{R}_{n+m}$ define the morphism 
\[\gamma_{n+m}: X_{n+m} \to  Y_{n+m}.\]
  The rational map  induced by $\bbar{R}_n$ is easily seen to be 
\[ \gamma_n p^{n+m}_n = \pi^{n+m}_n \gamma_{n+m}: X_{n+m}  \dra Y_n,\]
 and the rational map induced by  $\bbar{R}_m^{\sigma^n}$ is
\[ \gamma_m r^{n+m}_m = \rho^{n+m}_m \gamma _{n+m}: X_{n+m} \dra Y_m.\]
If $i \neq n$, then $p^{n+m}_n$ and $r^{n+m}_m$ are defined at all points in $\alpha_{n+m}^{-1}(w_i)$.  

We wish to apply Lemma~\ref{lem-definemap}.  To do so, we must calculate the divisors and base loci on $X_{n+m}$ associated to  the vector spaces $\bbar{R}_{n+m}$, $\bbar{R}_n$, and $\bbar{R}_m^{\sigma^n}$.   

For $0 \leq i \leq m+n$, let $F_i$ be the effective exceptional Weil divisor $\alpha_{n+m}^{-1}(w_i)$.  Now, $\sh{I}_{F_i}$ is the expansion of $\mf{m}^{n+m}_i$ to $X_{n+m}$.  By \cite[Proposition~7.1]{Ha},
 the expansion of $\mf{m}^{n+m}_i$ to $X'_{n+m}$ is Cartier; the ideal sheaf $\sh{I}_{F_i}$ is its pullback to $X_{n+m}$ and is thus also Cartier.  

By Lemma~\ref{lem-baseideal}, we have
\[ D^{X_{n+m}}(\bbar{R}_{n+m}) = \alpha_{n+m}^* D_{n+m}  - F_0 - \cdots - F_{n+m}.\]

Let $G_3 :=  D^{X_{n+m}}(\bbar{R}_{n+m})$.    Let 
\[ G_1: = \alpha_{n+m}^* D_{n} - F_0 - \cdots - F_{n-1},\]
and let
\[ G_2 := \alpha_{n+m}^* \sigma^{-n}D_{m} - F_{n+1} \cdots - F_{n+m}.\]
By assumption, 
$G_1 - D^{X_{n+m}}(\bbar{R}_n)$ and $G_2 - D^{X_{n+m}}(\bbar{R}_m^{\sigma^n})$ are both effective and supported on $F_n$.  That is, the base locus of  the rational functions  in $\bbar{R}_n$ with respect to the Cartier divisor $G_1$ is contained in  $F_n$.   Likewise, the base locus of the rational functions in $\bbar{R}_m^{\sigma^n}$ with respect to the Cartier divisor $G_2$ is also contained in $F_n$.

We now apply Lemma~\ref{lem-definemap} to the multiplication $\bbar{R}_n \bbar{R}_m^{\sigma^n} \subseteq \bbar{R}_{n+m}$.  We have that 
\[G_3 - G_1 - G_2 = -F_n + \alpha^*_{n+m} \sigma^{-n} \Omega.\]
Now, $\sigma^{-n} \Omega \cap O(w) \subseteq \{w_n\}$.  Thus by  Lemma~\ref{lem-definemap}, the rational map
\[ \xymatrix{\pi^{n+m}_n = (\pi^{n+m}_n \gamma_{n+m} ) \gamma_{n+m}^{-1}: Y_{n+m} \ar@{-->}[r] & Y_n}\]
is defined at every point of $\gamma_{n+m}(F_i)$ for every $i \neq n$.  That is, if $i \neq n$, then $\pi^{n+m}_n$ is defined at all points of $E^{n+m}_i$.    That  $\pi^{n+m}_n(E^{n+m}_i) = E^n_i$ is immediate.

The proof of (2) is symmetric:  we use the multiplication $\bbar{R}_m \bbar{R}_n^{\sigma^m} \subseteq \bbar{R}_{n+m}$.
\end{proof}

We will use the combinatorial Lemma~\ref{lem-M} to study  $\pi^{n+m}_n$ and $\rho^{n+m}_n$  for large $n$.

\begin{defn}
Let $X$ be a scheme and let $E \subset X$ be a proper closed subscheme of  $X$.  A {\em birational  image of $E$} is an open neighborhood $U$ of $E$ in $X$, together with a proper birational morphism
\[ \phi: U \to V \]
to some scheme $V$, such that $\phi|_{U \smallsetminus E}$ is an isomorphism onto its image.  

We put a preorder $\succcurlyeq$ on the collection of birational local images of $E$ as follows.  Let $\phi_1: U_1 \to V_1$ and $\phi_2:  U_2 \to V_2$ be two birational local images of $E$.  Then we say that $\phi_1$ {\em covers} $\phi_2$ and write $\phi_1 \succcurlyeq \phi_2$ if there is an open neighborhood $V'$ of $\phi_1(E)$ and a morphism $\psi:  V' \to V_2$ so that $\psi \phi_1 = \phi_2$ as birational maps from $U$ to $V_2$.  If $\phi_1 \succcurlyeq \phi_2 \succcurlyeq \phi_1$, we say that $\phi_1 \approx \phi_2$. The  {\em  local images of $E$}   are the $\approx$-congruence classes of birational  images of $E$; we denote the partial order on  local images of $E$ induced by $\succcurlyeq$ by $\geq$.  
\end{defn}

\begin{lemma}\label{lem-ACC}
Let $E$ be a proper closed subscheme of $X$ and let $\mb{M}$ be a collection of  local images of $E$.  Then $(\mb{M}, \geq)$ satisfies the ascending chain condition.
\end{lemma}

\begin{proof}
Let $\{ \phi_n:  X \dra V_n\}$ be an increasing sequence of representatives of  local images of $E$.  Then there are birational maps 
$\psi_n: V_n \dra V_{n-1}$, defined on a neighborhood of $\phi_n(E)$,  so that the diagram
\[ \xymatrix{
X \ar@{-->}[r]^{\phi_n} \ar@{-->}[rd]_{\phi_{n-1}}	& V_n \ar@{-->}[d]^{\psi_n} \\
& V_{n-1}
}\]
commutes for all $n$.  In particular, we have
\[ \dim(\phi_1(E)) \leq \dim(\phi_2(E)) \leq \cdots \leq \dim(E).\]
Therefore $\dim \phi_n(E)$ is equal to some constant $d$  for $n \geq n_0$.

Let $n \geq n_0$ and choose $h \leq\dim E$.  Let 
$W^h_n $ be the set of $ x \in E $ so that there is an irreducible component of $\phi_n^{-1} \phi_n(x)$ of dimension $\geq h$.
Now, 
\[
 \phi^{-1}_{n+1} \phi_{n+1}(x) \subseteq \phi^{-1}_{n+1} \psi_{n+1}^{-1} \psi_{n+1} \phi_{n+1}(x) 
 = \phi^{-1}_{n+1} \psi^{-1}_{n+1} \phi_n(x) = \phi^{-1}_n \phi_n(x).
\]
Thus $W^h_{n+1} \subseteq W^h_n$.  There is thus some $n_1 \geq n_0$ so that for $n \geq n_1$, we have
$W^h_n = W^h_{n+1}$ for all $0 \leq h \leq \dim E$.  We claim that for $n > n_1$, $\psi_n$ is finite.  

To see this, suppose the claim is false.   Then there is a positive-dimensional irreducible subscheme $C$ in $\phi_{n}(E)$ that $\psi_n$ contracts to a point.  Let $e := \dim \phi^{-1}_n(C) - \dim C$ and let $f := \dim \phi^{-1}_n(C)$; note that $f > e$.  Now, a generic point in $\phi^{-1}_{n}(C)$ is in $W^e_{n}$ but not in $W^f_{n}$.  However, we have
$\phi^{-1}_{n}(C) \subseteq W^f_{n-1}$.  But $n-1 \geq n_1$ so $W^f_n = W^f_{n-1}$.  This gives a contradiction.

Thus for $n \geq n_0$ there are open neighborhoods $V'_n$ of $\phi_n(E)$ so that $\psi_n(V'_n) \subseteq V'_{n-1}$ and so that $\psi_n |_{V'_n}$ is finite and birational.   By finiteness of the integral closure, for $n \gg 0 $ the map $\psi_n$ is a  local isomorphism at all points of $\phi_n(E)$; that is, for $n \gg 0 $ we have $\phi_n \approx \phi_{n+1}$ as claimed.
\end{proof}

\begin{lemma}\label{lem-W2}
 Let $w \in W$.  Define $E^n_i$ as in \eqref{def-Eni}.  
Then 
there are  integers $n_1$ and $b \geq 1$ so that: 

$(1)$  If $n \geq n_1$ and $i \leq n-b$ then for all $m \geq 1$, $\pi^{n+m}_n$ is a local isomorphism at all points of $E^{n+m}_i$.

$(2)$ If $n \geq n_1$  then for all $m \geq 1$ and $i \geq m+b$, $\rho^{n+m}_n$ is a local isomorphism at all points of $E^{n+m}_i$.
\end{lemma}

\begin{proof}
It suffices to prove the lemma for $m = 1$.  Our key claim is that the $\{ E^n_i \st 1 \leq i \leq n-1\}$ are a set of birational local images of $\Ed$ that satisfy the hypotheses of Lemma~\ref{lem-M}.  

More precisely, fix $i \geq 1$ and $n \geq i+1$.  By Lemma~\ref{lem-nice}, the birational map
\[ \alpha_n^{-1} \sigma^{-(i-1)} \alpha_2 = (p^n_{i+1})^{-1} (r^{i+1}_2)^{-1}: X_2 \dra X_n \]
is defined and is a local isomorphism in a neighborhood of $\Ed$.  
Let
\[ a^n_i := \gamma_n \alpha_n^{-1} \sigma^{-(i-1)} \alpha_2:  X_2 \dra Y_n.\]
This maps $\Ed$ onto $E^n_i$ and is a birational local image of $\Ed$.  

Consider the diagram
\[ \xymatrix{
X_2 \ar@{-->}[r]
\ar@{-->}[rd]
 \ar@/^1pc/@{-->}[rr]^{a^{n+1}_i} \ar@/_1pc/@{-->}[rrdd]_{a^n_i}
& X_{n+1} \ar[r]_{\gamma_{n+1}} \ar[d]^>>>>{p^{n+1}_{n}}	& Y_{n+1} \ar@{-->}[dd]^{\pi^{n+1}_{n}} \\
& X_{n} \ar[rd]^{\gamma_{n}}	&& \\
&& Y_{n}. }\]
By Lemma~\ref{lem-W}, $\pi^{n+1}_n$ is defined on $E^{n+1}_i$, so   $a^{n+1}_i \succcurlyeq a^{n}_i$.   Likewise,  consider the diagram
\[ \xymatrix{
X_2 \ar@{-->}[rr]^{a^{n+1}_{i+1}} \ar@{-->}[rrd]_{a^{n}_{i}}	&& Y_{n+1} \ar@{-->}[d]^{\rho^{n+1}_{n}} \\
&& Y_{n}
}\]
By Lemma~\ref{lem-W}, $\rho^{n+1}_n$ is defined on $E^{n+1}_{i+1}$.  Thus  $a^{n+1}_{i+1} \succcurlyeq a^{n}_{i}$.  

Let $\mb{M}: = \{a^n_i \st 1 \leq i \leq n-1\}$.  By Lemma~\ref{lem-ACC}, $(\mb{M}, \succcurlyeq)$ satisfies the ascending chain condition.  Applying Lemma~\ref{lem-M}, we see that there are some $N, b$ so that if $n \geq N$ and $1 \leq j \leq n-b$, then $a^{n}_i \approx a^{n-1}_i $ (as local images of $\Ed$), and if 
$b \leq j \leq n-1$ then $a^n_i \approx a^{n-1}_{i-1}$.  That is, if $1 \leq j \leq n-b$, then $\pi^n_{n-1}$ is a local isomorphism at all points of $E^n_j$, and if $b \leq j \leq n-1$ then $\rho^n_{n-1}$ is a local isomorphism at $E^n_j$. 

Now let $i \leq 0$.  Note that $a^n_1$ is defined on $\alpha_2^{-1}(w_i)$ for all $n \geq 2$.    Further, if $a^2_1$ is a local isomorphism at $\alpha_2^{-1}(w_i)$, then so are all $a^n_1$.  There are finitely many  $i \leq 0$ where $a^2_1$ is not a local isomorphism at $a_2^{-1}(w_i)$.  Since $\pi^n_{n-1} a^n_1 = a_1^{n-1}$, we see that $a^n_1 \succcurlyeq a^{n-1}_1$ as  local images of $\alpha_2^{-1}(w_i)$. Thus by increasing $N$ if necessary, we may ensure that $\pi^{n+1}_{n}$ is a local isomorphism at $E^{n+1}_j$ for all $n \geq N$ and $j \leq n-b$.    Likewise, and again increasing $N$ if needed, we may ensure that $\rho^{n+1}_{n}$ is a local isomorphism at $E^{n+1}_j$ for all $n \geq N$ and $j \geq b$.
\end{proof}

We are finally ready to prove:

\begin{proposition}\label{prop-noprob-W}
Assume Assumption-Notation~\ref{assnot6} and Notation~\ref{not-maps}.  For all $n \gg 0$, the sets $\gamma_n(\alpha_n^{-1}(\mathbb{W}))$ and $\gamma_n(U_n)$ are disjoint.
\end{proposition}
\begin{proof}
 Let $w \in W$, and adopt Notation~\ref{not-nice}.  
Let $n_1$ and $b$ be the integers constructed in Lemma~\ref{lem-W2}; let $N \geq \max\{n_1, 2b\}$ be such that for $n \geq N$, $\tau_n$ is an automorphism of $\gamma_n(U_n)$.   This exists by Proposition~\ref{prop-stableY}.  

Suppose there is some $e \in E_N$ and $u \in U_N$ such that $\gamma_N(e) = \gamma_N(u) =x$.  As $e$ and $u$ are in different connected components of $\gamma_N^{-1}(x)$, clearly $x$ is of finite order, say $k$,  under $\tau$.  Let $i$ be such that $\alpha_N(e) = w_i$.  

First suppose that $i \leq N-b$.  If $i \geq 0$, let $n := N + (i+1)k$; if $i<0$, let $n := N + k$.  As $i < N< n$, $(p^n_N)^{-1}$ is defined at $e$; let $e' := (p^{n}_N)^{-1}(e)$.   Let $u' := (p^n_N)^{-1}(u)$.  Then
\[  \pi^n_N \gamma_n (u') = \gamma_N p^n_N(u')  = x = \gamma_N p^n_N(e') = \pi^n_N \gamma_n(e').\]
Note that $i< N$, so $\pi^n_N$ is defined at $\gamma_n(e')$.  
By the choice of $N$ and $n$, $\pi^n_N$ is one-to-one on $\gamma_n(U_n) \cup E^n_i$ and so 
\[ \gamma_n(e') = \gamma_n(u').\]

But now, as $(\tau_N)^k(x) = x$, we have
\[ x = \pi^n_N\gamma_n(e') = \pi^n_N (\tau_N)^{n-N} \gamma_n(e') =  \rho^n_N \gamma_n(e') = 
\gamma_N r^n_N(e').\]
Our assumption on $i$ ensures that  $n - N \neq i$ and so $r^n_N(e')$ is well-defined.   As 
\[ \alpha_N r^n_N(e') = w_{i-(n-N)},\]
we see that $r^n_N(e') \not\in \{e', u\}$.  
We have produced a new point in $\gamma_N^{-1}(x)$.  Continuing, we may produce infinitely many such points, which is impossible.  Thus $i > N-b$.  

Arguing symmetrically, we obtain that $i < b$.  Since $N \geq 2b$, we see that no such $e$ can exist.
\end{proof}

We now construct $Y$.  

\begin{theorem}\label{thm-stableY}
Assume Assumption-Notation~\ref{assnot6}.  Then there are a projective variety $Y$ and a finite birational morphism $\theta: X\to Y$  such that for all $n \gg 0$ the rational map from $Y$ to $\PP^{N_n}$ induced by the rational functions in $\bbar{R}_n$ is a closed immersion at every point of $Y \smallsetminus \theta(\mathbb{W})$. 
Further, there  are a numerically trivial  automorphism $\phi$ of $Y$ such that $\theta \sigma = \phi \theta$, an ample and $\phi$-ample invertible sheaf $\sh{M}$ on $Y$ so that $\theta^* \sh{M} = \Lsh$, and a locally principal subscheme $\Phi$ of $Y$ so that $\Omega = \theta^* \Phi$.  For this $Y$, $\sh{M}$, and $\phi$, we have $R \subseteq B(Y, \sh{M}, \phi)$.   For $n \gg 0$,  the rational functions in $\bbar{R}_n$ correspond to sections of the invertible sheaf
$ \sh{I}_{\Phi}  \sh{M} \sh{M}^{\phi} \cdots \sh{M}^{\phi^{n-1}}$, and  their 
 base locus  is equal (set-theoretically) to $\theta(W_n)$.    
\end{theorem}
\begin{proof}
  We continue to use Notation~\ref{not-FINDY} and  Notation~\ref{not-maps}, so $m$ is such that $Y_m$ is stable, and $C \cup P$ is the subset of $U$ on which $\beta_m$ is not  a local isomorphism.  By Proposition~\ref{prop-noprob-W}, by increasing $m$ if necessary we may assume also that
   \beq\label{noprobW}
    \gamma_m(\alpha_m^{-1}(\mathbb{W})) \cap \gamma_m(U_m) = \emptyset.
    \eeq
   Let $\tau$ be the birational automorphism $\tau_m$ of $Y_m$.

 Let $H \subseteq \mathbb{W}$ be the set
\[\{ x \in \mathbb{W} \st \mbox{ either $\beta_m$ is undefined at $x$ or $\beta_m$ is not a local isomorphism at $x$} \}.\]
  Thus $H = \{h_1, \ldots, h_s\}$ is the finite  set of ``bad points'' of $\beta_m$ that do not lie on $C \cup P$.   Let $G := \alpha_m^{-1}(H)$.

We claim that the sets $\beta_m(U \smallsetminus (C \cup P))$, $\beta_m(C \cup P) $, $\gamma_m(G) = \gamma_m \alpha_m^{-1}(H)$, and $\gamma_m\alpha_m^{-1} (\mathbb{W} \smallsetminus H)$ are pairwise disjoint.  To see this, recall that $\beta_m$ is a local isomorphism at all points of $X \smallsetminus (C \cup P \cup H)$.  
Thus if $x \in U \smallsetminus (C \cup P)$, then $\beta_m^{-1}$ is defined at $\beta_m(x)$.   As $\alpha_m^{-1}$ is defined at $x$, if $x' \in X_m$ with $\gamma_m(x') = \beta_m(x)$, then $x'  = \alpha_m^{-1}(x)$.     Thus $\beta_m(U \smallsetminus (C \cup P))$ is disjoint from the other three sets.  That $\beta_m(C \cup P)$ is disjoint from the other sets follows;  recall that $\beta_m(U) \cap \gamma_m \alpha_m^{-1}(\mathbb{W}) = \emptyset$.

If $x \in \alpha_m^{-1} (\mathbb{W} \smallsetminus H)$ and $x' \in \alpha_m^{-1}(\mathbb{W})$ with $\gamma_m(x) = \gamma_m(x')$, then note that $\beta_m^{-1}$ is defined at $\gamma_m(x)$.  Therefore, 
\[ \alpha_m(x') = \beta_m^{-1} \gamma_m(x) = \alpha_m(x)\]
and $x' \not \in \alpha_m^{-1}(H)$.  This completes the proof of the claim.

To construct $Y$,  let
\[ V_1 := X \smallsetminus (C \cup P)\]
and let
\[ V_2 := Y_m \smallsetminus \gamma_m(G).\]
Let $V_{12} := V_1 \cap \alpha_m \gamma_m^{-1}(V_2)$, and let $V_{21} = V_2 \cap \gamma_m \alpha_m^{-1} (V_1)$.  
By the claim just previous,  $V_{12} = V_1  \smallsetminus H$  and $V_{21} = V_2 \smallsetminus \beta_m(C \cup P)$.   Further, $\beta_m(V_{12}) = V_{21}$; note that    $\beta_m$ is defined and is a local isomorphism at all $x \in V_{12}$. 


As $\beta_m$ defines a bijection between $V_{12}$ and $V_{21}$ that is a local isomorphism at each point, it is an isomorphism between $V_{12}$ and $V_{21}$.  By  \cite[Example~2.3.5]{Ha} there is a scheme $Y$ given by glueing $V_{1}$ and $V_2$ along the isomorphism $\beta_m: V_{12} \to V_{21}$.   For $i = 1, 2$ let $\psi_i$ be the induced map from $V_i$ to $Y$.

We now construct the automorphism $\phi$ of $Y$.  Let 
\[ V_{22} := V_2  \smallsetminus \gamma_m(\alpha_m^{-1} (\sigma^{-1}(H))).\]
We define morphisms
\[ \phi_1 := \psi_1 \sigma: V_1 \to Y,\]
\[ \phi_{21} := \psi_1 \sigma \beta_m^{-1}: V_{21} \to Y,\]
and
\[ \phi_{22} := \psi_2 \tau: V_{22} \to Y.\]

We check that $\phi_1$, $\phi_{21}$, and $\phi_{22}$ are well-defined; that is, that they are in fact morphisms.   First, $V_1$ is $\sigma$-invariant by construction, so $\sigma(V_1) \subseteq V_1$ and  $\phi_1$ is well-defined.  
Since $\beta_m^{-1}(V_{21}) = V_{12} \subseteq V_1$, $\phi_{21}$ is also well-defined.   Now, if $y \in V_{22} \cap \gamma_m(U_m)$, then, using \eqref{noprobW}, we have that $\tau(y) \in \gamma_m(U_m) \subseteq V_2$ and so $\phi_{22}$ is defined at $y$.  
Finally, if $y \in V_{22} \cap \gamma_m \alpha_m^{-1} (\mathbb{W})$, then  $\beta_m^{-1}$ is defined at $y$. Let $x = \beta_m^{-1}(y) \in \mathbb{W} \smallsetminus H \smallsetminus \sigma^{-1}(H)$. As $\sigma(x) \not \in H$, the map  $\tau = \beta_m \sigma \beta_m^{-1}$ is defined at $y$.  Further, $\beta_m$ is a local isomorphism at $\sigma(x)$, and so $\beta_m \sigma(x) \not \in \gamma_m(G)$ and $\tau(y) \in V_2$.    Thus $\psi_2$ is defined at $\tau(y)$.  

We next claim that $V_{21} \cup V_{22} = V_2$.  To see this, let $y \in V_2 \smallsetminus V_{22} = \gamma_m \alpha_m^{-1} (\sigma^{-1} H) \cap V_2$.  Then there is $x \in \sigma^{-1} (H)$ so that $y \in \gamma_m\alpha_m^{-1}(x)$; as $y \in V_2$, therefore  $x \not \in H$.  
As $x $ is certainly not in $C \cup P$, we see that $x \in V_{12}$, and $\beta_m(x) = y \in V_{21}$.


The diagram 
\[ \xymatrix{ 
& V_2 \ar@{-->}[ld]_{\sigma \beta_m^{-1}} \ar@{-->}[rd]^{\tau} &  \\
X \ar@{-->}[rr]^{\beta_m} \ar@{-->}[rd]_{\psi_1} && Y_m \ar@{-->}[ld]^{\psi_2} \\
& Y & }\]
of rational maps commutes by construction.  Note that the left side of this diagram gives $\phi_{21}$ and the right side gives $\phi_{22}$, considered as rational maps from $V_2$ to $Y$.  Thus $\phi_{21}$ and $\phi_{22}$ agree where both are defined; in particular, they agree on $V_{21} \cap V_{22}$.   By \cite[page~88]{Ha}, the morphisms $\phi_{21}$ and $\phi_{22}$ glue to give a birational morphism $\phi_2: V_2 \to Y$.  It is clear that  $\phi_1 = \phi_2 \beta_m$ on $V_{12}$, and so $\phi_1$ and $\phi_2$ glue via $\beta_m: V_{12} \to V_{21}$ to give a morphism  $\phi: Y \to Y$.  As $\phi$ is a local isomorphism at every point of $Y$, it is an automorphism of $Y$ by Lemma~\ref{lem-INVERSE}.   

Now let $V_3 := X \smallsetminus H$.  Note that $\beta_m$ is defined on $V_3$, and $\beta_m(V_3) = V_2$, by \eqref{noprobW}.  Define 
\[\psi_3 := \psi_2 \beta_m: V_3 \to Y.\]
  Now, $V_3 \cup V_1 = X$, and $V_3 \cap V_1 = V_{12}$.  By construction, $\psi_3 = \psi_1$ on $V_{12}$.   Thus we may glue $\psi_1$ and $\psi_3$  to obtain a morphism 
\[\theta: X \to Y.\]
   Clearly $\theta \sigma= \phi \theta$.  Furthermore, as both $\psi_3$ and $\psi_1$ are finite maps, $\theta$ is finite.   

Clearly $Y$ is integral.  We claim that $Y$ is also separated.  To see this, consider the diagonal $\Delta_Y = \{(y, y)\} \subseteq Y \times Y$.  This is the image of the diagonal $ \Delta_X \subseteq X\times X$ under the finite morphism $\theta \times \theta$.  As $X$ is separated, $\Delta_X$ is closed.  By  \cite[Exercise~3.5]{Ha} the finite morphism $\theta \times \theta$  is closed.  Thus  $\Delta_Y$ is also closed, and so $Y$ is separated.  Thus $Y$ is a variety.  

For all $n  \geq 1$, the rational functions in $\bbar{R}_n$ induce a rational map $\mu_n: Y \to \PP^{N_n}$ for appropriate $N_n$.   By construction, for $n \geq m$,  the indeterminacy locus of $\mu_n$ is equal to $\theta(W_n)$.  In particular,  it is contained in $\theta(\mathbb{W})$ and so supported at smooth points of $Y$.   Further, 
note that
if $n \geq m$ and  $x \in \theta(U)$, that locally at $x$ the rational map $\mu_n$ factors through the local isomorphism 
\[ \xymatrix{ (\pi^n_m)^{-1} \psi_2^{-1}: Y \ar@{-->}[r] &  Y_n \subseteq \PP^{N_n}. }\]
Thus $\mu_n$ is locally a closed immersion at any point of $Y \smallsetminus \theta(\mathbb{W})$.

By resolving the indeterminacy locus of $\mu_n$, we obtain a variety $Y'_n$, a  morphism $\xi_n: Y'_n \to Y$ and a morphism $\nu_n: Y'_n \to  \PP^{N_n}$ so that the diagram
\[ \xymatrix{
Y'_n \ar[rd]^{\nu_n} \ar[d]_{\xi_n} \\
Y \ar@{-->}[r]^{\mu_n} & \PP^{N_n} }
\]
commutes.    For all $n$, let $\sh{N}_n := \nu_n^* \struct(1)$ and let 
\[ \sh{K}_n := \bigl( (\xi_n)_* \sh{N}_n\bigr)^{**}.\]
Away from the indeterminacy locus of $\mu_n$, $\sh{K}_n $ is isomorphic to $ \mu_n^* O(1)$ and is invertible.  
 As any rank 1 reflexive module over a regular local ring is invertible,  $\sh{K}_n$ is invertible on the indeterminacy locus of $\mu_n$ as well, and therefore is an invertible sheaf on $Y$ for all $n \geq m$.  Thus $\bbar{R}_n \subseteq H^0(Y, \sh{K}_n)$, and the (set-theoretic) base locus of the sections $\bbar{R}_n$  of $\sh{K}_n$ is precisely $\theta(W_n)$ for $n \geq m$.  

For $n \geq m$, consider the Weil  divisor corresponding to  the invertible sheaf $\theta^* \sh{K}_n$ on $X$.  Away from the finitely many points in $ W_n$, this is equal to $\Delta_n - \Omega$.   As $X$ is smooth at all points of $ W_n$, by extending this equality to all of $X$, we obtain that
\[ \sh{I}_{\Omega}  \Lsh_n =   \struct_X(\Delta_n - \Omega) = \theta^* \sh{K}_n\]
for $n \geq m$.  

Let $\sh{M}:= \bigl( \sh{K}_m (\sh{K}_{m+1})^{-1} \bigr)^{\phi^{-m}}$, and let 
$ \sh{Z} := \sh{K}_m^{-1} \otimes \sh{M} \cdots \otimes \sh{M}^{\phi^{m-1}}$.  Then 
\[ \theta^* \sh{M} = \bigl( (\sh{I}_{\Omega} \Lsh_m) (\sh{I}_{\Omega} \Lsh_{m+1})^{-1} \bigr)^{\sigma^{-m}} \cong \bigl( \Lsh^{\sigma^m} \bigr)^{\sigma^{-m}} \cong \Lsh,\]
and 
\[\theta^* \sh{Z} = (\sh{I}_{\Omega} \Lsh_m)^{-1} \Lsh_m \cong (\sh{I}_{\Omega})^{-1}.\]
As $\theta^* \sh{Z}$ is an effective Cartier divisor, so is $\sh{Z}$; that is, $\sh{Z}^{-1}$ is an ideal sheaf defining a locally principal curve on $Y$.  We will denote this curve by $\Phi$; by construction, $\theta^* \Phi = \Omega$.  Note that $\sh{K}_n \cong \sh{I}_{\Phi} \sh{M}_n$.   
 
Recall that  $\Lsh$ is ample.  As $\theta$ is finite, $\sh{M}$ is ample by 
 \cite[Proposition 2.6.2]{EGA-III.1}.  Thus $Y$ carries an ample line bundle and so is projective.  The numeric action of $\phi$ is clearly still trivial, and so $\sh{M}$ is also $\phi$-ample by \cite[Theorem~1.7]{AV}.  We have $\bbar{R}_n \subseteq H^0(Y, \sh{K}_n) \subseteq H^0(Y, \sh{M}_n)$.
\end{proof}

We remark that the fact that $Y$ is a projective variety may also be deduced from \cite[Proposition~7.4]{RS}.

\section{Surjectivity in large degree}\label{SURJ}

 Let $\surfD = \surfdata$ be transverse surface data, and suppose that $R \subseteq T(\surfD)$ is a graded ring.   In this section, we digress for a moment to  establish sufficient conditions for $R$   and $T(\surfD)$ to be equal in large degree.   Our methods involve reducing the question to one involving subrings of twisted homogeneous coordinate rings of $\sigma$-invariant curves on $X$.  We wish to use the results of \cite{AS} on subrings of idealizers on curves; however, as those were proved only for reduced and irreducible curves, we repeat the proofs here in a more general context.  
  
\begin{theorem}\label{thm-surj}
Suppose that the surface data 
\[\surfD = \surfdata\]
 is transverse.  Let 
$T := T (\surfD)$
 and let 
 $\sh{T} :=\sh{T}( \surfD).$
   Let $R$ be a subalgebra of $T$  with $R_1 \neq 0$, and fix $0 \neq z \in R_1$.  Let $\bbar{R}_n := R_n z^{-1}$ and let $\sh{R}_n(X) := \bbar{R}_n \cdot \struct_X$.  
   
   Suppose  that  $\sh{R}_n(X) = \sh{T}_n$ for $n \gg 0$.   Let $W$ be the cosupport of $\sh{ADC}$ and let 
    \[\mathbb{W} := \bigcup_{p \in W} O(p).\]
  Further assume that for all $n \gg 0$, the rational map defined on $X$ by the rational functions in $\bbar{R}_n$ is birational onto its image and is a closed immersion at each point in $X \smallsetminus \mathbb{W}$.  Then $R_n = T_n$ for $n \gg 0$.  
\end{theorem}

We will prove Theorem~\ref{thm-surj} in several steps.  
We first establish some notation.   If ${\Gamma}$ is a $\sigma$-invariant proper subscheme of $X$, then $\sigma$ restricts to an automorphism of ${\Gamma}$, which we also denote by $\sigma$.  For any such ${\Gamma}$,  let $B_{\Gamma} = B({\Gamma}, \Lsh|_{\Gamma}, \sigma)$.  We may consider $T$ and $R$ to be subrings of $B(X, \Lsh, \sigma)$;  we will let $T_{\Gamma}$, respectively $R_{\Gamma}$, be the image of $T$, respectively $R$, under the natural map from $B(X, \Lsh, \sigma)$ to $B_{\Gamma}$.

\begin{proof}[Proof of Theorem~\ref{thm-surj}]
By Lemma~\ref{lem-ample}, the sequence of bimodules  $\{(\sh{R}_n)_{\sigma^n}\}$ is left and right ample; thus by Lemma~\ref{lem-fg}, $T$ is a finitely generated left and right $R$-module.  
Let $J_l := \rm{l. ann}_R(T/R)$ and let $J_r := \rm{r. ann}_R (T/R)$.  Note that $J_l$ is a graded right ideal of $T$ and that $J_r$ is a graded left ideal of $T$.  Our assumptions imply that $R$ and $T$ have the same graded quotient ring, and thus $J_l$ and $J_r$ are nonzero.  Let $K := J_r J_l$.  Then $K \neq 0$ is a nonzero graded ideal both of $R$ and of $T$.  Note also that by Theorem~\ref{thm-surfprop-main}, both $\sh{T}$ and $T$ are  noetherian. 

By \cite[Proposition~4.10]{S-surfprop}, there is a $\sigma$-invariant ideal sheaf $\sh{K}$ on $X$ such that for $n \gg 0$, we have  $K_n = H^0(X, \sh{K} \sh{R}_n) z^n$.    Let $\Gamma$ be the $\sigma$-invariant closed subscheme defined by $\sh{K}$; then $\dim \Gamma \leq 1$.      By transversality of the defining data for $\sh{R}$,  the $\sigma$-invariant subscheme $\Gamma $ is disjoint from $\mathbb{W}$, and  $\Omega \cap \Gamma $ consists of points of infinite order.  Let $\sh{J}$ be the ideal sheaf on $\Gamma $ of the scheme-theoretic intersection $\Omega \cap \Gamma $.  Since $\Omega$ is locally principal, 
 $\sh{R}_n |_{\Gamma} = \sh{J}( \Lsh_n |_{\Gamma })$ for $n \geq 1$.  

Note that $R/K$ and $R_{\Gamma}$ are equal in large degree, and $T/K$ and $T_{\Gamma}$ are equal in large degree.  Note also that as for $n \gg 0 $ the rational functions in $\bbar{R}_n$ define a closed immersion at all points of $X \smallsetminus \mathbb{W}$, that their restrictions to $\Gamma \subseteq X \smallsetminus \mathbb{W}$ also define a closed immersion for $n \gg 0$.

We claim that $R_{\Gamma}$ and $T_{\Gamma}$ are equal in large degree.  Before proving this claim, we give a lemma generalizing a result of Artin and Stafford.

\begin{lemma}\label{lem-AS4.6}
{\em (cf. \cite[Lemma~4.6]{AS})}
  Suppose, in addition, that there are no proper $\sigma$-invariant subschemes  $Y$ of $\Gamma$ so that $(T_Y)/(R_Y)$ is infinite-dimensional, and that there are $\sigma$-invariant ideal sheaves $\sh{I}_1, \ldots, \sh{I}_{\ell} \subseteq \struct_{\Gamma}$ so that $\sh{I}_1 \sh{I}_2 \cdots \sh{I}_{\ell} = 0 $ on $\Gamma$.  Then $T_{\Gamma}/R_{\Gamma}$ is finite-dimensional.
\end{lemma}

\begin{proof}
The proof is similar to the proof of \cite[Lemma~4.6]{AS}; we give it in detail because some of the details are different in our  more general context.

Suppose, in contrast, that $T_{\Gamma}/R_{\Gamma}$ is infinite-dimensional.  We first note that if $J$ is a nonzero graded  ideal of $T_{\Gamma}$, then there is a graded ideal $J' \supseteq K$ of $T$ so that $J = J' /K$ in large degree.  By \cite[Proposition~4.10]{S-surfprop}, in large degree $J'$ consists of sections of $\sh{R}_n$ that vanish on some $\sigma$-invariant proper subscheme $Y$ of $\Gamma$, and so (in large degree)  $T_{\Gamma}/J = T_Y$.  As $R_Y$ and $T_Y$ are equal in large degree but $R_{\Gamma}$ and $T_{\Gamma}$ are not equal in large degree, $J \not\subseteq R_{\Gamma}$.

Thus  $R_{\Gamma}$ and $T_{\Gamma}$ have no nonzero ideals in common.  
By induction, we may assume that $\ell = 2$.  Let $Z_1$ and $Z_2$, respectively, be the subschemes of $Y$ defined by $\sh{I}_1$ and $\sh{I}_2$, respectively.   Let $\sh{M} := \Lsh|_{\Gamma}$ and let $B := B(\Gamma, \sh{M}, \sigma)$.  For $i = 1, 2$, let 
\[ K_i := \bigoplus_{n\geq 0} H^0(\Gamma, \sh{I}_i \sh{M}_n)z^n \subseteq B,\]
 and let $M_i := K_i \cap T_{\Gamma}$.  Note that the $M_i$ are two-sided ideals of $T_{\Gamma}$.  As $\sh{I}_1$ is an $\struct_{Z_2}$-module, the right and left actions of $T_{\Gamma}$ on $M_1$ factor through $T_2 := T_{\Gamma}/M_2$.

Now, $T_{\Gamma}$ is noetherian and so $M_1$ is a finitely generated left and right $T_2$-module.  Let $R_2 := (R_{\Gamma} + M_2)/M_2 \subseteq T_2$.  By hypothesis, $R_2$ and $T_2$ are equal in large degree.  
Thus $R_2$ is noetherian, and  both $M_1$ and $N := R_{\Gamma} \cap M_1$ are finitely generated left and right $R_2$-modules.  
Let $N' := T_2 N T_2 \subseteq M_1$.  This is a finitely generated $R_2$-module; since $R_2$ and $T_2$ are equal in large degree, $N'$ and $N$ are also equal in large degree.

There is thus some $n_0$ so that $N_{\geq n_0} = N'_{\geq n_0}$ is a left and right $T_2$-module.  That is, $N_{\geq n_0}$ is an ideal of $T_{\Gamma}$.  As $N_{\geq n_0} \subseteq R_{\Gamma}$ and $R_{\Gamma}$ and $T_{\Gamma}$ have no nonzero ideals in common, $N_{\geq n_0} = 0$.

Since $(R_{\Gamma})_{\geq n_0} \cap M_1 = 0$, we have an injection $(R_{\Gamma})_{\geq n_0} \hookrightarrow T_1$.  
  This implies that the map defined by the sections $(\bbar{R}_{\Gamma})_n$ of $\sh{M}_n$ factors through $Z_1$ for $n \ge n_0$, and so is not an embedding.  This gives a contradiction.
\end{proof}

We return to the proof of Theorem~\ref{thm-surj}.  We show that $R_{\Gamma}$ and $T_{\Gamma}$ are equal in large degree.  
By noetherian induction on $\Gamma$, we may assume  for any proper $\sigma$-invariant closed subscheme $Y \subseteq \Gamma$ that $R_Y$ has finite codimension in $T_Y$.  

We first suppose that $\Gamma$ is reducible.  Let $k$ be such that $\sigma^k$ fixes all irreducible components of $\Gamma$.      The hypotheses of Lemma~\ref{lem-AS4.6} thus hold for $R_{\Gamma}\ver{k}$ and $T_{\Gamma}\ver{k}$.  Applying  Lemma~\ref{lem-AS4.6}, we see that $T_{\Gamma}\ver{k}/R_{\Gamma}\ver{k}$ is finite-dimensional.

We show that  this implies that $T_{\Gamma}/R_{\Gamma}$ is finite-dimensional.    Let $\sh{F}_n := \sh{R}_n|_{\Gamma}$, and let 
\[ \sh{F} := \bigoplus_{n \geq 0} ( \sh{F}_{n})_{\sigma^n}.\]  
The noetherian property of $\sh{T}$ descends to the $\struct_{\Gamma}$-bimodule algebra $\sh{F}$, and so $\sh{F}$ and its Veronese $\sh{F}\ver{k}$ are noetherian.   
 As the restriction of an ample sequence to a $\sigma$-invariant subscheme, the sequence of bimodules $\{ (\sh{F}_{rk})_{\sigma^{rk}}\}_{r \geq 0}$ is left and right ample. 
 
  Recall that $\sh{J} = \sh{I}_{\Omega} \struct_{\Gamma} \subseteq \struct_{\Gamma}$.  Let $\sh{M} := \Lsh|_{\Gamma}$. 
Fix $0  \leq i \leq k-1$, and let 
\[ \sh{P}  := \bigoplus_{n \geq 0} \sh{J} \sh{M}_{i+nk} = \bigoplus_{n \geq 0} \sh{F}_{i+nk}.\]
Let $\sh{N}:= (\Lsh_i^{\sigma^{-i}})|_{\Gamma}$.  
   The sheaf $\sh{JN}^{-1}$ on $\Gamma$ is  invertible, and so the submodule lattices of the right $\sh{F}\ver{k}$-modules  $\sh{P}$ and 
   \[(\sh{JN}^{-1})_{\sigma^{-i}} \otimes  \sh{P} \cong \bigoplus_{n \geq 0} (\sh{J} \sh{J}^{\sigma^{-i}} \sh{M}_{nk})_{\sigma^{nk}} \subseteq \sh{F}\ver{k}\]
    are isomorphic.   In particular, $\sh{P}$ is a coherent right $\sh{F}\ver{k}$-module.
   
   Fix $n_0$ so that if $n \geq i+n_0k$, then $(\bbar{R}_{\Gamma})_n$ generates $\sh{P}_n$.  Let $n_1 \geq n_0$ be such that 
   \[ \sh{F}_{i+n_0k} \oplus \cdots \oplus \sh{F}_{i+n_1k}\]
   generates $\sh{P}_{\geq i+n_0k}$ as a right $\sh{F}\ver{k}$-module.     Then for $r \geq n_1$, we have
  \begin{multline*}
(T_{\Gamma})_{i+rk} = H^0(\Gamma, \sh{JM}_{i+rk})z^{i+rk} = H^0 \bigl(\Gamma, \sum_{j = n_0}^{n_1} \sh{F}_{i+jk}  \sh{F}^{\sigma^{i+jk}}_{(r-j)k}\bigr) z^{i+rk} \\\
   =  \sum_{j=n_0}^{n_1}  H^0\bigl(\Gamma, \sh{F}_{i+jk} \sh{F}^{\sigma^{i+jk}}_{(r-j)k}\bigr) z^{i+rk}.
  \end{multline*}
    By Lemma~\ref{lem-surjective}, for fixed $\ell$ and for $r \gg 0$, we have 
     \[ (R_{\Gamma})_{\ell} (T_{\Gamma})_{rk }  =  H^0\bigl(\Gamma, \sh{F}_{\ell} \sh{F}_{rk}^{\sigma^{\ell}}\bigr) z^{\ell+rk}.\]  
    Recall that $R_{\Gamma}\ver{k}$ and $T_{\Gamma}\ver{k}$ are equal in large degree. Thus, by taking $r \gg 0$, we obtain that
     \begin{multline*}
       (R_{\Gamma})_{i + rk} \subseteq (T_{\Gamma})_{i + rk} =  \sum_{j=n_0}^{n_1} (R_{\Gamma})_{i+jk} (T_{\Gamma})_{(r-j)k} = \\
     \sum_{j=n_0}^{n_1} (R_{\Gamma})_{i+jk} (R_{\Gamma})_{(r-j)k} \subseteq (R_{\Gamma})_{i + rk}.
     \end{multline*}
     Since this holds for $ 0 \leq i \leq k-1$, $R_{\Gamma}$ has finite codimension in $T_{\Gamma}$.

  Now suppose that $\Gamma$ is irreducible but not reduced.  Then the nilradical $\sh{N}$ of $\struct_{\Gamma}$ is a $\sigma$-invariant nilpotent ideal sheaf on $\Gamma$; so the hypotheses of Lemma~\ref{lem-AS4.6} hold for $R_{\Gamma}$ and $T_{\Gamma}$, with $\sh{I}_1 = \sh{I}_2 = \cdots = \sh{I}_{\ell} = \sh{N}$.  We see again that $T_{\Gamma}/R_{\Gamma}$ is finite-dimensional.
    
  Thus we have reduced to considering the case that $\Gamma$ is reduced and irreducible.  Now, if $\Omega \cap \Gamma = \emptyset$, then $R_{\Gamma}$ and $T_{\Gamma}$ are equal in large degree by \cite[Theorem~4.1]{AS}; in particular, this holds if $\Gamma$ is a point.  If $\Omega \cap \Gamma$ is nonempty, and $\Gamma$ is  a reduced and irreducible curve, then $T_{\Gamma}/R_{\Gamma}$ is finite-dimensional by  \cite[Proposition~5.4]{AS}.

We have thus shown that there is an ideal $K$ of $T$ that is contained in $R$ and so that $(R/K)_n = (R_{\Gamma})_n = (T_{\Gamma})_n = (T/K)_n$  for $n \gg 0$.   Thus  $R_n = T_n$ for $n \gg 0$.
\end{proof}

\section{The proof of the main theorem}\label{FINAL}

In this final section, we prove Theorem~\ref{ithm-surfclass}, Theorem~\ref{ithm-main} and Theorem~\ref{ithm-chi2}.  We also discuss  possible extensions to birationally commutative graded algebras of higher GK-dimension.

We restate our main result.
\begin{theorem}\label{thm-surfclass}
Let $\kk$ be an uncountable algebraically closed field and let $R$ be a   birationally commutative projective surface over $\kk$.  Then  there is  transverse surface data $\mb{E}=(Y, \sh{M}, \phi, \sh{E}, \sh{F}, \sh{G}, \Phi)$, where $\phi$ is numerically trivial, and a positive integer $N$ so that $R \ver{N} \cong T(\mb{E})$. 
\end{theorem}
\begin{proof}
 By Corollary~\ref{cor-normalize}, there are a positive integer $\ell$ and transverse surface data
 \[ \surfD = \surfdata\]
 so that $X$ is normal, $\sigma$ is numerically trivial, $\bbar{R}_{\ell}$ generates $K$, and 
 \[ \sh{R}(X) \ver{\ell} \cong \sh{T}(\surfD).\]  
Note that Assumption-Notation~\ref{assnot6} holds for $R \ver{\ell}$.  Let $Z_n$ be the base locus of $\bbar{R}_{n\ell}$.
 
  Let $W$ be the cosupport of $\sh{ADC}$.  
Let
 \[\mathbb{W} := \bigcup_{p \in W} O(p).\]  
 By Theorem~\ref{thm-stableY} there are a projective variety $Y$, a numerically trivial  automorphism $\phi$ of $Y$, an ample invertible sheaf $\sh{M}$ on $Y$,  a locally principal subscheme $\Phi$ of $Y$, and a finite birational morphism $\theta: X \to Y$ so that for $n \gg 0$ the rational functions $\bbar{R}_{n\ell}$ induce a closed immersion into projective space at every point of $Y \smallsetminus \theta(\mathbb{W})$ and so that $\theta \sigma = \phi \theta$,  $\theta^* \sh{M} = \Lsh$, and $\theta^* \Phi = \Omega$; further, set-theoretically the base locus of the sections $\bbar{R}_{n\ell} \subseteq H^0(Y, \sh{I}_{\Phi} \sh{M}_n)$ is equal to $\theta(Z_{n})$.  

Let $\Lambda$, $Z$, and $\Lambda'$ be the subschemes of $X$ defined respectively by $\sh{A}$, $\sh{D}$, and $\sh{C}$.  Note that $\theta$ is a local isomorphism at all points in $\Lambda \cup Z \cup \Lambda'$.  Let $\sh{E}$ be the ideal sheaf of $\theta(\Lambda)$, let $\sh{F}$ be the ideal sheaf of $\theta(Z)$, and let $\sh{G}$ be the ideal sheaf of $\theta(\Lambda')$.  

The ideal sheaves $\sh{E}$, $\sh{F}$, and $\sh{G}$ on $Y$ pull back to $\sh{A}$, $\sh{D}$, and $\sh{C}$ respectively.  Furthermore, by working locally at each point of $\theta(W)$, we see that $\sh{I}_{\Phi} \sh{E}\sh{G} \subseteq \sh{F}$.   As distinct points in the cosupport of $\sh{D}$ have distinct $\sigma$-orbits, distinct points in the cosupport of $\sh{F}$ have distinct $\phi$-orbits.   

Let 
\[ \surfD' := (Y, \sh{M}, \phi, \sh{E}, \sh{F}, \sh{G},  \Phi).\]
By construction, $\sh{R}(Y)\ver{\ell} \cong \sh{T}(\surfD')$.  
We check that the data $\surfD'$ is transverse.    It follows from the corresponding statements for $X$ that 
$\{ \phi^n(\theta(Z))\}_{n \in \ZZ}$, $\{ \phi^n(\theta(\Lambda))\}_{n \geq 0}$, and $\{\phi^n(\theta(\Lambda'))\}_{n \leq 0}$ are critically transverse, and that for any reduced and irreducible $\Gamma \subset Y$, the set
$\{ n \st \phi^n(\Gamma) \subseteq \Phi\}$ is finite.  Since $\Phi$ is locally principal, by \cite[Lemma~3.1]{S-surfprop}, $\{\phi^n \Phi\}$ is critically transverse.  Thus $\mb{D}'$ is transverse surface data.

We have seen that for $n \gg 0$, the sections in $\bbar{R}_{n\ell}$ define a closed immersion at all points of $Y \smallsetminus \theta(\mathbb{W})$.  
 Theorem~\ref{thm-surj} now implies that  there is some $k \geq 1$ so that 
\[ R_{n\ell} = T_n\]
for $n \geq k$.  
Thus if $\mathbb{E}$ is the surface data given by Lemma~\ref{lem-Veronese-data}
so that 
\[ \sh{T} (\mathbb{D}')\ver{k} = \sh{T}(\mb{E}),\]
 then $\mb{E}$ is transverse and 
\[ R\ver{k\ell} = T (\mathbb{E}),\]
as claimed.
\end{proof}

Theorem~\ref{ithm-surfclass} follows directly from Theorem~\ref{thm-surfclass}.

\begin{proof}[Proof of  Theorem~\ref{ithm-main}]
 By Theorem~\ref{thm-surfclass}, we may replace $R$ with a Veronese subring so that there is transverse surface data $\surfD = \surfdata$ with $R \cong T(\surfD)$.  Let 
$\sh{C}':= (\sh{D} : \sh{I}_{\Omega} \sh{A})$.
Let 
\[\mb{F}:=(X, \Lsh, \sigma, \sh{A}, \sh{D}, \sh{C}', \Omega)\]
 and let $S:= T(\mb{F})$.  Note that $R_{\geq 1}$ is a left ideal of $S$.  
By \cite[Corollary~4.7]{S-surfprop}, $R$ and $\I_S(R_{\geq 1})$ are equal in large degree.

Now let 
$\sh{A}':= (\sh{D} :  \sh{C}')$.
Let 
\[\mb{G}:= (X, \Lsh, \sigma, \sh{A}', \sh{D}, \sh{C}', \emptyset)\]
 and let $T:=T(\mb{G})$.  Then $S_{\geq 1}$ is a right ideal of $T$, and applying \cite[Corollary~4.7]{S-surfprop} again, we see that $S$ and $\I_T(S_{\geq 1})$ are equal in large degree.  The pair $(\sh{A}', \sh{C}')$ is maximal with respect to 
\[ \sh{A}'\sh{C}' \subseteq \sh{D}.\]
If $\sh{D} = \struct_X$, then $T$ is a twisted homogeneous coordinate ring.  If $\sh{D}$ is nontrivial, then $T$ is an ADC ring, as defined in the introduction, and $R$ is in large degree equal to a left idealizer inside a right idealizer in $T$.
\end{proof}

We prove one corollary of Theorem~\ref{thm-surfclass}.  Recall that a noetherian connected $\NN$-graded algebra $R$ satisfies left (right) $\chi_j$ if, for all $i \leq j$, and for all finitely generated left (right) graded $R$-modules $N$, the groups
\[ \lim_{n \to \infty} \Ext^i_R(R/R_{\geq n}, N)\]
are finite-dimensional.  If $R$ satisfies left and right $\chi_{j}$ for all $j$, then  $R$ satisfies $\chi$.

\begin{corollary}\label{cor-chi}
Let $R$ be  a birationally commutative projective surface.  The following are equivalent:
\begin{enumerate}
\item $R$ satisfies left $\chi_2$;
\item $R$ satisfies right $\chi_2$;
\item $R$ satisfies $\chi$;
\item  Some Veronese $R\ver{k}$ is a twisted homogeneous coordinate ring $B(X, \Lsh, \sigma)$, where $\Lsh$ is $\sigma$-ample.
\end{enumerate}
\end{corollary}
\begin{proof}
For any $k \geq 1$, $R$ is a finitely generated left and right  $R\ver{k}$-module, by \cite[Lemma~4.10(iii)]{AS}.  Thus left or right $\chi_i$ hold for $R$ if and only if they hold for $R \ver{k}$.  By Theorem~\ref{thm-surfclass}, for some $k$ we have  $R\ver{k} \cong T(\surfD)$, where $\surfD$ is transverse surface data.  The four conditions are equivalent for $T(\surfD)$ by \cite[Theorem~5.9]{S-surfprop}.  
\end{proof}

Theorem~\ref{ithm-chi2} follows from Corollary~\ref{cor-chi}.

To end this chapter, we make a few remarks on a possible extension of Theorem~\ref{thm-surfclass} to rings of higher GK-dimension.  One class of rings one would like to understand are what we will call for purposes of discussion {\em GK 5 birationally commutative projective surfaces}:    that is, connected $\NN$-graded noetherian domains $R$  whose graded quotient ring is of the form
\[ K [z, z^{-1}; \sigma] \]
for some field $K$ of transcendence degree 2 and geometric, but non-quasi-trivial automorphism $\sigma$ of $K$.  (Recall  from Theorem~\ref{thm-R-model} that  such a ring must have GK-dimension 5.)  The results of \cite{RS} on algebras generated in degree 1 apply  also to GK 5 birationally commutative projective surfaces.   We conjecture that Theorem~\ref{thm-surfclass}  extends to algebras of this type that are not generated in degree 1.  To make this extension would require modifying proofs in the current work that use the hypothesis that $\sigma$ is quasi-trivial.  

This hypothesis is often convenient, and is sometimes essential to our proofs. We have not been able to construct alternate proofs for all results that rely on quasi-triviality.  It appears, however, that one can use the Kodaira classification of surfaces to show that there are limited possibilities for defining data of GK 5 surfaces.   This work is ongoing.

There are also GK 4 finitely $\NN$-graded subrings of  $K[z, z^{-1}; \sigma]$ for appropriate  $K$ of transcendence degree 2 and $\sigma \in \Aut_{\kk}(K)$.  By Theorem~\ref{thm-R-model}, for these algebras $\sigma$ is not geometric.  Rogalski and Stafford  conjecture \cite{RS} that these GK 4 algebras are never noetherian.  Proving or disproving this conjecture is an interesting and delicate question; if some prove to be noetherian, one would presumably hope to classify them, as well.

Finally, one can ask whether some version of Theorem~\ref{thm-surfclass} works for higher-dimensional varieties.  To begin to answer this question, it clearly makes sense to work over $\CC$ and to work with automorphisms $\sigma$ that have a projective model $X$ on which they  are quasi-trivial.  By \cite[Lemma~2.23]{DQZhang2008}, this is equivalent to assuming that some power of $\sigma$ lies in $\Aut^o(X)$.   

To try to extend Theorem~\ref{thm-surfclass} to this context, one has first to understand na\"ive blowups at higher-dimensional subvarieties.  What transversality properties of the movement of a surface under an automorphism of a 4-fold ensure that the corresponding na\"ive blowup is noetherian?  Assuming that these questions are solved, one would then need to modify the many proofs in this paper that depend on the fact that $X$ is a surface; for example, we constantly use the fact that by definition, the base locus of the sections in $\bbar{R}_n$ must be 0-dimensional, since the corresponding invertible sheaf is defined as a minimal object.  Similar problems would arise in using the techniques of \cite{RS}, since their proofs rely on  Castelnuovo-Mumford regularity techniques that require low dimension of the subschemes one is blowing up.  So far little is known about these interesting questions.

\bibliographystyle{amsalpha}
\bibliography{biblio}

\end{document}